\documentclass[plain]{sn-jnl}
\jyear{2021}%

\usepackage{amsmath}
\usepackage{amssymb}
\usepackage{subcaption}
\usepackage[numbers]{natbib}

\theoremstyle{thmstyleone}%
\newtheorem{theorem}{Theorem}%  meant for continuous numbers
\newtheorem{lemma}{Lemma}
\newtheorem{corollary}{Corollary}
\newtheorem{result}{Main result}

\theoremstyle{thmstyletwo}%
\newtheorem{example}{Example}%
\newtheorem{remark}{Remark}%

\theoremstyle{thmstylethree}%
\newtheorem{assumption}{Assumption}

\begin{document}

\title[A framework for randomized time-splitting in LQ optimal control]{A framework for randomized time-splitting in linear-quadratic optimal control}


\author[1]{\fnm{D.W.M.} \sur{Veldman}}\email{daniel.veldman@math.fau.de}

\author[1,2,3]{\fnm{E.} \sur{Zuazua}}\email{enrique.zuazua@fau.de}

\affil[1]{\orgdiv{Chair in Dynamics, Control, and Numerics (Alexander-von-Humboldt Professorship)}, \orgname{Friedrich-Alexander Universit\"at (FAU) Erlangen-Neuremberg}, \orgaddress{\street{Cauerstrasse 11}, \postcode{91052}, \city{Erlangen}, \country{Germany}}}

\affil[2]{\orgdiv{Departamento do Matem\'aticas}, \orgname{Universidad Autonoma de Madrid}, \orgaddress{\street{Ciudad Universitaria de Cantoblanco}, \postcode{28049}, \city{Madrid}, \country{Spain}}}

\affil[3]{\orgdiv{Chair of Computational Mathematics}, \orgname{Fundaci\'on Deusto}, \orgaddress{\street{Av.\ de las Universidades 24}, \postcode{48007}, \city{Bilbao}, \state{Basque-Country}, \country{Spain}}}


\abstract{Inspired by the successes of stochastic algorithms in the training of deep neural networks and the simulation of interacting particle systems, we propose and analyze a framework for randomized time-splitting in linear-quadratic optimal control. In our proposed framework, the linear dynamics of the original problem is replaced by a randomized dynamics. To obtain the randomized dynamics, the system matrix is split into simpler submatrices and the time interval of interest is split into subintervals. The randomized dynamics is then found by selecting randomly one or more submatrices in each subinterval. 

We show that the dynamics, the minimal values of the cost functional, and the optimal control obtained with the proposed randomized time-splitting method converge in expectation to their analogues in the original problem when the time grid is refined. The derived convergence rates are validated in several numerical experiments. Our numerical results also indicate that the proposed method can lead to a reduction in computational cost for the simulation and optimal control of large-scale linear dynamical systems. }

\keywords{Random Batch Method, Operator Splitting, Optimal Control, Model Predictive Control}


\pacs[MSC Classification]{65C99, 49M99, 65L20, 37M05}

\maketitle

\section{Introduction} \label{sec:intro}
Solving an optimal control problem for a large-scale dynamical system can be computationally demanding. This problem appears in numerous applications. One example is Model Predictive Control (MPC), which requires the solution of several optimal control problems on a receding time horizon \cite{grune2017, ko2020}. Another example is the training of Deep Neural Networks (DNNs), which can be approached as an optimal control problem for a large-scale nonlinear dynamical system, see, e.g., \cite{E2017, benning2019, esteve2021, esteve2021sparse, ruizbalet2021}. 
Because the computational cost for gradient-based deterministic optimization algorithms explodes on large training data sets, Neural Networks (NNs) are typically trained using stochastic optimization algorithms such as stochastic gradient descent or stochastic (mini-)batch methods, see, e.g., \cite{bottou2018}. In such methods, the update direction for the parameters of the NN is not computed based on the complete training data set, but on a subset of the available training data that is chosen randomly in each iteration. It can be shown that such methods converge in expectation to a (local) minimum of the considered cost functional, see, e.g., \cite{bottou2018}. 

These successes inspired the development of Random Batch Methods (RBMs) for the simulation of interacting particle systems \cite{jin2020,li2020,jin2020convergence}. Because the number of interactions between $N$ particles is of order $N^2$, the forward simulation of a system with a large number of particles is computationally demanding. A RBM reduces the required computational cost by reducing the number of considered interactions as follows. First, the considered time interval is divided into a number of subintervals of length $\leq h$. In each subinterval, particles are grouped in randomly chosen batches (of at least two particles) and only the interactions between particles in the same batch are considered. The number of considered interactions now grows as $PN$, where $P$ is the size of the considered batches, and a significant reduction in computational time can be achieved when $P \ll N$. It can be shown that the expected error introduced by this process is proportional to $\sqrt{h}$, where $h$ denotes (an upper bound on) the length of the considered time intervals, see \cite{jin2020}. 

The computation of optimal controls for interacting particle systems is even more computationally demanding than the forward simulation because it requires several simulations of the forward dynamics and the associated adjoint problem, see, e.g., \cite{lee1967}. Because the optimal control for the RBM-approximated dynamics can be computed significantly faster than the control for the original dynamics, it has been proposed in \cite{ko2020} to control the original system with the controls optimized for the RBM dynamics. The numerical experiments in \cite{ko2020} indeed indicate that this approach can lead to a reasonably good approximation of the control for the original system. In \cite{ko2020}, the control of the original dynamics with the RBM-optimal controls is combined with an MPC strategy, which creates additional robustness against the errors introduced by the RBM-approximation. However, even for the simplest case that does not consider the combination with MPC, a formal proof that the optimal control computed for the RBM-approximated dynamics indeed converges to the optimal control for the original system for $h \rightarrow 0$ was not given. 

In this paper, we study, motivated by the ideas from \cite{ko2020}, the classical linear-quadratic (LQ) optimal control problem constrained by randomized dynamics. Extensions of these results to a nonlinear setting are not only of interest for the control of interacting particle systems as considered in \cite{ko2020}, but have also applications in the training of certain DNNs which can be viewed as (the time discretization) of an optimal control problem, see, e.g., \cite{E2017, benning2019, esteve2021, esteve2021sparse, ruizbalet2021}. The results for the LQ problem in this paper form a starting point for the study of these more involved problem settings. 

%In the considered linear setting, the achievable reduction in computational cost is similar to the reduction that can be obtained for classical operator splitting algorithms for discretized Partial Differential Equations (PDEs), see, e.g., \cite{quarteroni1994}, which date back to \cite{douglas1955, peaceman1955, douglas1956}. These splittings can also be incorporated in our proposed frame work. 
%In these methods,  a discretization of the operator is viewed as the sum of certain simpler operators. These simpler operators are then used sequentially in every time step, which can lead to a significant reduction in computational cost. Using the same splittings in our randomized framework leads to a similar reduction in computational cost. %The main difference with the approaches proposed in \cite{douglas1955, peaceman1955, douglas1956} and our approach is then that the order in which the simpler operators are used is determined randomly in our framework. 

In this paper, we propose a framework for the simulation and optimal control of large-scale linear dynamical systems. In our proposed framework, the system matrix is split into submatrices and the time interval of interest is split into subintervals of length $\leq h$. The randomized dynamics is then found based on the randomly selected submatrices in each subinterval. Similarly as in \cite{jin2020,li2020,jin2020convergence}, we show that the randomized dynamics converges to the dynamics of the original system at a rate $\sqrt{h}$. The main contributions of this paper concern the LQ optimal control problem in which the original dynamics is replaced by these randomized dynamics. In particular, we show that the minimal values of the cost functional and the corresponding optimal controls for the RBM-dynamics converge (in $L^2$ and in expectation) to their analogues for the original dynamics when $h \rightarrow 0$. The found convergence rates are validated by several numerical examples. Numerical results also indicate that the proposed method can lead to a reduction in computational cost. 

The remainder of this paper is structured as follows. Section \ref{sec:method} contains a precise description of our proposed stochastic simulation method and a summary of the main results of the paper. Section \ref{sec:convergence} contains the detailed proofs of the convergence of the proposed method. The proposed method and the obtained convergence results are illustrated by several numerical examples in Section \ref{sec:examples}. The conclusions and discussions are presented in Section \ref{sec:conc}. 

\section{Proposed method and main results} \label{sec:method}

\subsection{Proposed method} \label{ssec:method}
We consider the evolution of a large-scale Linear Time Invariant (LTI) dynamical system of the form
\begin{equation}
\dot{x}(t) = Ax(t) + Bu(t), \qquad \qquad x(0) = x_0, \label{eq:dyn_x}
\end{equation}
where the state $x(t)$ evolves in $\mathbb{R}^N$, the control $u(t)$ evolves in $\mathbb{R}^q$, $A \in \mathbb{R}^{N \times N}$ is the system matrix, $B \in \mathbb{R}^{N \times q}$ is the input matrix, and $x_0 \in \mathbb{R}^N$ is the initial condition. 

A typical problem associated to the dynamics \eqref{eq:dyn_x} is to find the optimal control $u^*(t)$ that minimizes the quadratic cost functional
\begin{equation}
J(u) = \frac{1}{2}\int_0^T \left( (x(t)-x_d(t))^\top Q (x(t)-x_d(t))+u(t)^\top R u(t)  \right) \ \mathrm{d}t,  \label{eq:J}
\end{equation}
where the given target trajectory $x_d(t)$ evolves in $\mathbb{R}^N$, the weighting matrix $Q \in \mathbb{R}^{N \times N}$ is symmetric and positive semi-definite, and the weighting matrix $R\in \mathbb{R}^{q \times q}$ is symmetric and positive definite. It is well known that the optimal control $u^*(t)$ exists and that it is unique, see, e.g., \cite{minoux1986, kirk2004}.

\begin{remark} \label{rem:phi}
When the state-dimension $N$ is large, the optimal control $u^*(t)$ is typically computed using a gradient-based algorithm in which the gradient of $J(u)$ is computed from the adjoint state $\varphi(t)$ that satisfies (see, e.g., \cite{kirk2004})
\begin{equation}
-\dot{\varphi}(t) = A^\top \varphi(t) + Q (x(t) - x_d(t)), \qquad \qquad \varphi(T) = 0, \label{eq:dyn_phi}
\end{equation}
where $x(t)$ is the solution of \eqref{eq:dyn_x}. Note that the adjoint state $\varphi(t)$ is computed by integrating \eqref{eq:dyn_phi} backward in time starting from the final condition $\varphi(T) = 0$. The gradient of the cost functional $J(u)$ is then obtained as
\begin{equation}
\left(\nabla J(u)\right)(t) = B^\top \varphi(t) + R u(t). 
\end{equation}
\end{remark}

In our proposed randomized time-splitting method, the matrix $A$ is written as the sum of $M$ submatrices $A_m$ 
\begin{equation}
A = \sum_{m=1}^M A_m. \label{eq:Aform}
\end{equation}
Typically, the submatrices $A_m$ will be more sparse than the original matrix $A$. For ease of presentation, the results in this paper are presented under the following assumption. 

\begin{assumption} \label{ass:dissipativity}
The submatrices $A_m$ in \eqref{eq:Aform} are dissipative, i.e.\ 
$\langle x, A_m x \rangle \leq 0$ for all $x \in \mathbb{R}^N$ and all $m \in \{1,2,\ldots ,M \}$.
\end{assumption}

\begin{remark}
Note that there always exists a constant $a > 0$ such that the matrices $A_m - aI$ are dissipative for $m \in \{1,2,\ldots ,M \}$. Assumption \ref{ass:dissipativity} is therefore not essential for the convergence of the proposed method, but without Assumption \ref{ass:dissipativity} the error estimates are less clean and grow exponentially in time. This idea is made more precise in Remark \ref{rem:dissipativity} in Section \ref{ssec:forward_dynamics_deterministic}.
\end{remark}

We then choose a temporal grid in the time interval $[0,T]$
\begin{equation}
0 = t_0 < t_1 < t_2 < \ldots < t_{K-1} < t_K = T, 
\end{equation}
and denote
\begin{equation}
h_k = t_k - t_{k-1}, \qquad \qquad h = \max_{k \in \{1,2,\ldots, K\}} h_k. 
\end{equation}
In each of the $K$ subintervals $[t_{k-1}, t_k)$, we randomly select a subset of indices in $\{1,2,\ldots , M \}$. The idea of the proposed method is to consider a linear combination of the submatrices $A_m$ with the indices that have been selected for each time interval. This can lead to a significant reduction in computational time when the submatrices $A_m$ are well-chosen and only a small number of submatrices $A_m$ are selected in each time interval. 

To make this idea more precise, we enumerate all of the $2^M$ subsets of $\{1,2, \ldots, M \}$ as $S_1, S_2, \ldots S_{2^M}$. Note that one of the subsets $S_{\omega}$ will be the empty set. To every subset $S_{\omega}$ ($\omega \in \Omega := \{1, 2, \ldots , 2^M \}$) we then assign a probability $p_{\omega}$ with which this subset is selected. This probability is the same in each of the time intervals $[t_{k-1},t_k)$. Because we select only one subset $S_{\omega}$ in each time interval, the probabilities $p_{\omega}$ should satisfy
\begin{equation}
\sum_{\omega = 1}^{2^M} p_{\omega} = 1. \label{eq:sump}
\end{equation}
From the chosen probabilities $p_{\omega}$, we then compute the probability $\pi_m$ that an index $m \in \{1,2, \ldots, M \}$ is an element of the selected subset
\begin{equation}
\pi_m = \sum_{\omega \in \Omega_m} p_\omega, \qquad \qquad \Omega_m = \{ \omega \in \{1,2, \ldots, 2^M \} \mid m \in S_{\omega} \}. \label{eq:pim}
\end{equation}
Observe that $\Omega_m$ is the set of the indices $\omega$ of the sets $S_{\omega}$ that contain the index $m$. 
We need the following (weak) assumption on the selected probabilities $p_{\omega}$. 
\begin{assumption} \label{ass:positivity}
The probabilities $p_{\omega}$ ($\omega \in \{1,2, \ldots ,2^M \}$) are assigned such that 
\begin{itemize}
\item \eqref{eq:sump} is satisfied and
\item the probabilities $\pi_m$ defined in \eqref{eq:pim} are positive for all $m \in \{1, ,2, \ldots , M \}$. 
\end{itemize}
\end{assumption}

In each of the $K$ time intervals $[t_{k-1},t_k)$, we then randomly select an index $\omega_k \in \{1,2, \ldots , 2^M \}$ according to the chosen probabilities $p_\omega$ (and independently of the other indices $\omega_1, \omega_2, \ldots \omega_{k-1}, \omega_{k+1}, \omega_{k+1}, \ldots , \omega_K$). The selected indices form a vector
\begin{equation}
\boldsymbol \omega := (\omega_1, \omega_2, \ldots , \omega_K) \in \{1,2,\ldots ,2^M \}^K =: \Omega^K. 
\label{eq:def_Omega}
\end{equation}
For the selected $\boldsymbol \omega\in \Omega^K$, we then define a piece-wise constant matrix $t \mapsto \mathcal{A}_h(\boldsymbol \omega,t)$
\begin{equation}
\mathcal{A}_h(\boldsymbol \omega, t) = \sum_{m \in S_{\omega_k}} \frac{A_m}{\pi_m}, \qquad \qquad t \in [t_{k-1}, t_k). \label{eq:defAt}
\end{equation}
The scaling by $1/\pi_m$ assures that the expected value of $\mathcal{A}_h$ is $A$ because
\begin{equation}
\sum_{\omega=1}^{2^M} \sum_{m \in S_{\omega}} \frac{A_m}{\pi_m} p_\omega 
= \sum_{m=1}^M \sum_{\omega \in \Omega_m} \frac{A_m}{\pi_m} p_\omega 
= \sum_{m=1}^M \frac{A_m}{\pi_m}\pi_m 
= \sum_{m=1}^M A_m = A, \label{eq:expectationA_step1}
\end{equation}
where the first identity follows after interchanging the two summations using the definition of $\Omega_m$ in \eqref{eq:pim}, the second from the definition of $\pi_m$ in \eqref{eq:pim}, and the last identity from the decomposition of $A$ in \eqref{eq:Aform}. 

\begin{example} \label{ex:A12_1}
In the simplest situation, we decompose the original matrix $A$ into $M = 2$ matrices as $A = A_1 + A_2$. We then need to assign $2^M = 4$ probabilities $p_\ell$ to the subsets $S_1 = \{ 1 \}$, $S_2 = \{ 2 \}$, $S_3 = \{ 1,2 \}$, and $S_4 = \emptyset$. In this example, we choose $p_1 = p_2 = \tfrac{1}{2}$ and $p_3 = p_4 = 0$. This choice indeed satisfies Assumption \ref{ass:positivity} because $\pi_1 = p_1 + p_3 = \tfrac{1}{2} > 0$ and $\pi_2 = p_2 + p_3 = \tfrac{1}{2} > 0$. The matrix $\mathcal{A}_h(\boldsymbol \omega,t)$ is thus either equal to $2A_1$ with probability $p_1 = \tfrac{1}{2}$ or equal to $2A_2$ with probability $p_2 = \tfrac{1}{2}$. The expected value of $\mathcal{A}_h$ is then indeed $\tfrac{1}{2} 2 A_1 + \tfrac{1}{2} 2 A_2 = A_1 + A_2 = A$.
\end{example}

To reduce the computational cost for solving \eqref{eq:dyn_x}, the matrix $A$ is replaced by a $\mathcal{A}_h(\boldsymbol \omega,t)$ in the RBM. For the selected vector of indices $\boldsymbol \omega \in \Omega^K$, we thus obtain a solution $t \mapsto x_h(\boldsymbol \omega,t)$
\begin{equation}
\dot{x}_h(\boldsymbol \omega,t) = \mathcal{A}_h(\boldsymbol \omega,t) x_h(\boldsymbol \omega,t) + B u(t), \qquad \qquad x_h(\boldsymbol \omega,0) = x_0. \label{eq:dyn_xtilde}
\end{equation} 
The main contribution of this paper concerns the optimal controls computed based on the RBM-dynamics \eqref{eq:dyn_xtilde}. In particular, we consider the minimization of the functional
\begin{equation}
J_h(\boldsymbol \omega, u) = \frac{1}{2}\int_0^T \left( (x_h(\boldsymbol \omega,t)-x_d(t))^\top Q (x_h(\boldsymbol \omega,t)-x_d(t))+u(t)^\top R u(t)  \right) \ \mathrm{d}t, \label{eq:Jtilde}
\end{equation}
over all $u \in L^2(0,T ; \mathbb{R}^q)$ subject to the dynamics \eqref{eq:dyn_xtilde}. The minimizer of $J_h(\boldsymbol \omega, \cdot)$ depends on the selected indices $\boldsymbol \omega \in \Omega^K$ and is denoted by $u^*_h(\boldsymbol \omega,t)$. Because $R$ is positive definite, the minimizer $u^*_h(\boldsymbol \omega,t)$ exists and is unique. As we will show in \eqref{eq:bound_uhopt_step1}--\eqref{eq:bound_uhopt} in Section \ref{ssec:preliminaries}, the minimizers $u^*_h(\boldsymbol \omega,t)$ are uniformly bounded because $R$ is positive definite.

\begin{remark} \label{rem:phih}
Similarly as for the original cost functional $J(u)$ in \eqref{eq:J}, we can compute the optimal control $u_h(\boldsymbol \omega,t)$ that minimizes $J_h(\boldsymbol \omega,u)$ by a gradient-based algorithm. We can again compute the gradient of $J_h(\boldsymbol \omega,u)$ from the adjoint state $\varphi_h(\boldsymbol \omega,t)$ which satisfies 
\begin{equation}
-\dot{\varphi}_h(\boldsymbol \omega,t) = \left( \mathcal{A}_h(\boldsymbol \omega,t) \right)^\top \varphi_h(\boldsymbol \omega,t) + Q(x_h(\boldsymbol \omega,t) - x_d(t)), \qquad \varphi_h(\boldsymbol \omega,T) = 0. \label{eq:dyn_phih}
\end{equation}
The gradient of $J_h(\boldsymbol \omega,u)$ is then obtained as
\begin{equation}
\nabla J_h(\boldsymbol \omega,u) = B^\top \varphi_h(\boldsymbol \omega,t) + Ru(t). 
\end{equation}
Note that when the randomized dynamics for $x_h(\boldsymbol \omega,t)$ in \eqref{eq:dyn_xtilde} can be solved faster than the original dynamics for $x(t)$ in \eqref{eq:dyn_x}, the same reduction in computational cost is typically also obtained for the randomized adjoint equation \eqref{eq:dyn_phih} compared to the original adjoint equation \eqref{eq:dyn_phi}. Because the computation of the optimal control $u^*(t)$ (resp.\ $u_h^*(\boldsymbol \omega,t)$) requires several evaluations of the forward dynamics \eqref{eq:dyn_x} (resp.\ \eqref{eq:dyn_xtilde}) and the adjoint equation \eqref{eq:dyn_phi} (resp.\ \eqref{eq:dyn_phih}), it is natural to expect the same relative speed-up for $u_h^*(\boldsymbol \omega,t)$ (compared to $u_h^*(t)$) as for $x_h(\boldsymbol \omega,t)$ (compared to $x(t)$). This idea is confirmed by the numerical experiments in Section \ref{sec:examples}. 
\end{remark}

To conclude this subsection, we summarize the proposed approach to approximate the solution $x(t)$ of \eqref{eq:dyn_x} for a given control $u(t)$ and/or the optimal control $u^*(t)$ that minimizes $J(\cdot)$ in \eqref{eq:J} subject to \eqref{eq:dyn_x} in Algorithm \ref{alg:RBM}. The accuracy of the obtained solutions $x_h(\boldsymbol \omega,t)$ and/or $u^*_h(\boldsymbol \omega,t)$ depends on the chosen submatrices $A_m$ in Step 1, the chosen probabilities $p_\omega$ in Step 2, and the chosen time grid $t_0, t_1, \ldots, t_K$ in Step 3. This dependence is captured by the error estimates in the next subsection. 

It should be emphasized that we do not have that $\mathbb{E}[x_h(t)] = x(t)$ (for a fixed control $u(t)$) or that $\mathbb{E}[u^*_h(t)] = u^*(t)$. Repeating Step 4 in Algorithm \ref{alg:RBM} for different realizations of $\boldsymbol \omega$ and averaging the obtained results leads to better approximations of $\mathbb{E}[x_h(t)]$ and/or $\mathbb{E}[u^*_h(t)]$ and can therefore only improve the approximation of $x(t)$ and $u^*(t)$ to some extend. A better way to increase the accuracy of the proposed method is to repeat Algorithm \ref{alg:RBM} for a choice of submatrices $A_m$, probabilities $p_{\omega}$, and a time grid $t_0, t_1, \ldots, t_K$ that reduce the error estimates presented in the next subsection.

\begin{algorithm}
\hrule
\caption{The proposed randomized time-splitting method}
\label{alg:RBM}
\hrule
\textbf{Step 1} Decompose the matrix $A$ into $M$ submatrices $A_m$ as in \eqref{eq:Aform}, preferably such that Assumption \ref{ass:dissipativity} is satisfied. 

\textbf{Step 2} Enumerate the $2^M$ subsets of $\{1,2, \ldots , M \}$ and assign probabilities $p_1$, $p_2$, \ldots, $p_{2^M}$ such that Assumption \ref{ass:positivity} is satisfied. 

\textbf{Step 3} Divide the considered time interval $[0,T]$ into $K$ subintervals $[t_{k-1}, t_k)$ and choose an index $\omega_k$ according to the selected probabilities in Step 2 for each subinterval. Store the selected indices in a vector $\boldsymbol \omega = (\omega_1, \omega_2, \ldots \omega_K)$. 

\textbf{Step 4} Compute the solution $x_h(\boldsymbol \omega,t)$ of the dynamics \eqref{eq:dyn_xtilde} for a certain given control $u(t)$ and/or compute the control $u^*(\boldsymbol \omega,t)$ that minimizes $J_h(\boldsymbol \omega,\cdot)$ in \eqref{eq:Jtilde} subject to the dynamics \eqref{eq:dyn_xtilde}. 

\hrule 

\end{algorithm}

\begin{remark}
The presented framework is somewhat different from the problem setting considered in previous publications on RBMs for interacting particle systems, see, e.g., \cite{jin2020,li2020,jin2020convergence,ko2020}. Appendix \ref{app:interacting_particles} shows how these RBMs can be accommodated in our proposed framework. 
\end{remark}

\subsection{Main results}
The main results of this paper concern the effect of replacing the system matrix $A$ in the original LQ optimal control problem \eqref{eq:dyn_x}--\eqref{eq:J} by the randomized matrix $\mathcal{A}_h(\boldsymbol \omega,t)$ defined in \eqref{eq:defAt}. Clearly, the deviation of the randomized matrix $\mathcal{A}_h(\boldsymbol \omega,t)$ from the original matrix $A$ will influence the accuracy of the obtained results. The deviation of $\mathcal{A}_h(\boldsymbol \omega,t)$ from $A$ is measured by
\begin{equation}
\mathrm{Var}[\mathcal{A}] := \sum_{\omega=1}^{2^M} \left\| \sum_{m \in S_{\omega}} \frac{A_m}{\pi_m} - A \right\|^2 p_\omega, \label{eq:def_varA}
\end{equation}
where $\| \cdot \|$ denotes the operator norm. The quantity $\mathrm{Var}[\mathcal{A}]$ is thus the average squared distance of $\mathcal{A}_h(\boldsymbol \omega,t)$ from $A$, weighted with the probabilities $p_1$, $p_2$, \ldots, $p_{2^M}$ with which different values of $\mathcal{A}_h(\boldsymbol \omega,t)$ occur. Naturally, the error estimates below show that reducing $\mathrm{Var}[\mathcal{A}]$ will also reduce the errors introduced by the proposed randomized time-splitting method.

\setcounter{example}{0}
\begin{example}[continued] \label{ex:A12_2}
We again consider the situation from Example \ref{ex:A12_1} in which $A$ is decomposed into $M=2$ submatrices as $A = A_1 + A_2$ and $\mathcal{A}_h(\boldsymbol \omega,t)$ is either $2A_1$ or $2A_2$, both with probability $\tfrac{1}{2}$. We now compute the variance $\mathrm{Var}[\mathcal{A}]$ according to \eqref{eq:def_varA} and find
\begin{equation}
\mathrm{Var}[\mathcal{A}] = \| 2A_1 - A \|^2 p_1 + \| 2A_2 - A \|^2 p_2 = \| A_1 - A_2 \|^2.
\end{equation}
\end{example}

Examples \ref{ex:A123} and \ref{ex:popt} at the end of this subsection further illustrate how $\mathrm{Var}[\mathcal{A}]$ depends on the decomposition of $A$ into submatrices $A_m$ and the selected probabilities $p_{\omega}$.

\begin{remark} \label{rem:def_varAW}
When $A$ in an approximation of an unbounded operator as in the examples in Section \ref{sec:examples}, it is natural to introduce an additional (invertible) weighting matrix $W$ and compute
\begin{equation}
\mathrm{Var}_W[\mathcal{A}] := \sum_{\ell=1}^{2^M} \left\| \left( \sum_{m \in S_\ell} \frac{A_m}{\pi_m} - A \right)W \right\|^2 p_\ell.  \label{eq:def_varAW}
\end{equation}
Clearly, we want to choose $W$ such that $AW$ and the matrices $A_mW$ can be considered as approximations of bounded operators. In that case, $\mathrm{Var}_W[\mathcal{A}]$ is also an approximation of a finite quantity. A natural choice is $W = (A - \lambda I)^{-1}$ for some $\lambda$ in the resolvent of $A$. %Extending the results in this paper to a meaningful infinite-dimensional setting is a challenging topic for future research. 
\end{remark}

The first main result of this paper is an estimate for the difference 
\begin{equation}
e_h(\boldsymbol \omega,t) = x_h(\boldsymbol \omega,t) - x(t) \label{eq:res_def_e1}
\end{equation}
between the RBM-dynamics \eqref{eq:dyn_xtilde} and the original dynamics \eqref{eq:dyn_x}. 

\begin{result} \label{res:este0}
Assume that Assumptions \ref{ass:dissipativity} and \ref{ass:positivity} hold and that the input $u(t)$ in \eqref{eq:dyn_x} is the same as in the input $u(t)$ in \eqref{eq:dyn_xtilde}, then
\begin{equation}
\mathbb{E}[\lvert e_h(t) \rvert^2] \leq C_{[A,B,x_0,T,u]} h \mathrm{Var}[\mathcal{A}], \label{eq:res_este0}
\end{equation}
for all $t \in [0,T]$. 
\end{result}

The first main result follows directly from Theorem \ref{thm:este} in Subsection \ref{ssec:forward_dynamics_deterministic}. 

The expectation operator $\mathbb{E}$ is taken with respect to all possible outcomes $\boldsymbol \omega \in \Omega^K$. A precise definition will be given in Section \ref{ssec:preliminaries}. 
The constant $C_{[A,B,x_0,T,u]}$ can be taken as $(\| A\|T^2 + 2T)( \lvert x_0 \rvert+\lvert Bu \rvert_{L^ 1(0,T;\mathbb{R}^N)})^2$. The estimate thus only depends on the used submatrices $A_m$, the probabilities $p_{\omega}$, and the used temporal grid $t_0,t_1,\ldots , t_K$ through $h\mathrm{Var}[\mathcal{A}]$ defined in \eqref{eq:def_varA}. The proof of Main result \ref{res:este0} is inspired by the proofs of convergence of the RBM in \cite{jin2020, jin2020convergence}. 

The estimate \eqref{eq:res_este0} shows that the expected squared error is proportional to the temporal grid spacing $h$. We can thus make the expected squared error in the forward dynamics arbitrary small by reducing the grid spacing. Note that Markov's inequality, see, e.g., \cite{rohatgi2015}, shows that
\begin{equation}
\mathbb{P}[\lvert e_h(\boldsymbol \omega,t) \rvert > \varepsilon] = 
\mathbb{P}[\lvert e_h(\boldsymbol \omega,t) \rvert^2 > \varepsilon^2] < 
\frac{\mathbb{E}[\lvert e_h(t) \rvert^2]}{\varepsilon^2}.
\end{equation}
The probability that we select an $\boldsymbol \omega \in \Omega^K$ for which $\lvert e_h(\boldsymbol \omega,t) \rvert$ exceeds any given treshold $\varepsilon > 0$ is thus controlled by $\mathbb{E}[\lvert e_h(t) \rvert^2]$. According to \eqref{eq:res_este0}, we can make this probability as small as desired by choosing the temporal grid spacing $h$ small enough. However, one should also keep in mind that decreasing $h$ will increase the computational cost for the RBM-dynamics \eqref{eq:dyn_xtilde} and that the computational advantage of the RBM is lost when the required grid spacing is too small.

\setcounter{example}{0}
\begin{example}[continued] \label{ex:A12_3}
To illustrate why Main result \ref{res:este0} could be true, we again consider the situation from Example \ref{ex:A12_1} in which $A$ is decomposed as $A = A_1 + A_2$ and $\mathcal{A}_h(\boldsymbol \omega,t)$ is equal to $2A_1$ or $2A_2$, both with probability $\tfrac{1}{2}$. We additionally assume that $u(t) \equiv 0$, that the time grid $t_k = k T / K$ ($k \in \{0,1,2, \ldots , K \}$) is uniform with grid spacing $h = T/ K$, and that $A_1$ and $A_2$ commute. 
Because $u(t) = 0$, the solution of \eqref{eq:dyn_x} is $x(t) = e^{At}x_0$ and the solution of \eqref{eq:dyn_xtilde} is
\begin{align}
x_h(\boldsymbol \omega,T) = e^{2A_{\omega_K} h} \cdots  e^{2A_{\omega_2} h} e^{2A_{\omega_1} h} x_0 = e^{2A_1T_1(\boldsymbol \omega) + 2A_2 T_2(\boldsymbol \omega)} x_0. \label{eq:exA12comm_xh}
\end{align}
Here, $T_1(\boldsymbol \omega)$ and $T_2(\boldsymbol \omega)$ denote the times during which $A_1$ and $A_2$ are used, i.e.\
\begin{equation}
T_1(\boldsymbol \omega) = \frac{T}{K}\sum_{\ell=1}^K \chi_1(\omega_\ell), \qquad \qquad 
T_2(\boldsymbol \omega) = \frac{T}{K}\sum_{\ell=1}^K \chi_2(\omega_\ell),
\end{equation}
where the characteristic functions $\chi_1(\omega)$ and $\chi_2(\omega)$ are defined by the property that $\chi_i(\omega) = 1$ when $\omega = i$ and $\chi_i(\omega) = 0$ otherwise ($i \in \{ 1,2\}$). Note that the second identity in \eqref{eq:exA12comm_xh} uses that $A_1$ and $A_2$ commute.  Because $\mathbb{E}[\chi_1] = \mathbb{E}[\chi_2] = \tfrac{1}{2}$, it follows that $\mathbb{E}[T_1] = \mathbb{E}[T_2] = T/2$. When we now consider the limit $K \rightarrow \infty$ (so $h \rightarrow 0$), the law of large numbers states that $T_1$ and $T_2$ converge to $T/2$ (in probability). The RHS of \eqref{eq:exA12comm_xh} thus converges (in probability) to $e^{AT}x_0 = x(T)$ for $K \rightarrow \infty$. Note that the convergence in Main result \ref{res:este0} is in expectation, which is stronger than convergence in probability. 
\end{example}

We now present the main results aimed at the LQ optimal control problem constrained by randomized dynamics. Because the optimal control $u^*_h(\boldsymbol \omega,t)$ depends on the selected indices $\boldsymbol \omega$, we need the following result. The key difference with the first main result is that the input $u_h(\boldsymbol \omega,t)$ may now depend on the randomly selected indices $\boldsymbol \omega$. As will be explained at the start of Section \ref{sec:convergence}, this makes the arguments for the convergence of the RBM in \cite{jin2020, jin2020convergence} break down. 

Note that replacing $u(t)$ in \eqref{eq:dyn_x} and \eqref{eq:dyn_xtilde} by $u_h(\boldsymbol \omega,t)$ results in solutions $x(\boldsymbol \omega,t)$ and $x_h(\boldsymbol \omega,t)$ that now both depend on the selected indices $\boldsymbol \omega$. The second main result now gives a bound for the expected value of the difference
\begin{equation}
e_h(\boldsymbol \omega,t) = x_h(\boldsymbol \omega,t) - x(\boldsymbol \omega,t).
\end{equation}

\begin{result}
\label{res:este2} 
Consider any control $u_h : \Omega^K \rightarrow L^2(0,T ; \mathbb{R}^q)$. 
Assume that Assumptions \ref{ass:dissipativity} and \ref{ass:positivity} are satisfied and let $U$ be such that $\lvert u_h(\boldsymbol \omega) \rvert_{L^2(0,T; \mathbb{R}^q)} \leq U$ for all $\boldsymbol \omega \in \Omega^K$, then
\begin{equation}
\mathbb{E}[\lvert e_h(t) \rvert^2] \leq C_{[A,B,x_0,T,U]} h \mathrm{Var}[\mathcal{A}].
\end{equation}
\end{result}
The second result follows directly from Theorem \ref{thm:este2} in Subsection \ref{ssec:forward_dynamics_stochastic}. 

  Just as in the first main result, the expectation is taken over all possible values of $\boldsymbol \omega \in \Omega^K$ and the constant $C_{[A,B,x_0,T,U]}$ does not depend on the chosen submatrices $A_m$ in \eqref{eq:Aform}, the chosen probabilities $p_{\omega}$, and the used temporal grid.

Using this result, we can now obtain a no-gap result which shows that the minimal value of the cost functional $J_h(\boldsymbol \omega, u_h^*(\boldsymbol \omega))$ is (in expectation) close to the minimal value $J(u^*)$ in the original problem when $h \mathrm{Var}[\mathcal{A}]$ is small enough. 

\begin{result} \label{res:nogap}
Let $u^*(t)$ be the control that minimizes the cost functional $J(u)$ in \eqref{eq:J} and let $u_h^*(\boldsymbol \omega,t)$ be the control that minimizes the cost functional $J_h(\boldsymbol \omega, u)$ in \eqref{eq:Jtilde}. Then
\begin{equation}
\mathbb{E}[\lvert J_h(u^*_h) - J(u^*) \rvert] \leq C_{[A,B,x_0,Q,R,x_d,T]} \left(  \sqrt{h \mathrm{Var}[\mathcal{A}]} + h \mathrm{Var}[\mathcal{A}] \right). \label{eq:res_nogap}
\end{equation}
\end{result}

The third main result is identical to Theorem \ref{thm:nogap} in Subsection \ref{ssec:nogap}. 

For $h\mathrm{Var}[\mathcal{A}]$ small enough, Main result \ref{res:nogap} clearly implies that $\mathbb{E}[\lvert J_h(u^*_h) - J(u^*) \rvert] \leq C_{[A,B,x_0,Q,R,x_d,T]} \sqrt{h \mathrm{Var}[\mathcal{A}]}$, which is also the rate that is observed in numerical experiments. We keep the second term on the RHS of \eqref{eq:res_nogap} to assure that the estimate is valid for all values of $h\mathrm{Var}[\mathcal{A}]$, and not just for sufficiently small values of $h\mathrm{Var}[\mathcal{A}]$.

By Markov's inequality, this result thus implies that, for any $\varepsilon > 0$, the probability that $\lvert J(u^*_h(\boldsymbol \omega)) - J(u^*) \rvert > \varepsilon$ can be made arbitrarily small by reducing the temporal grid spacing $h$. 

The next main result shows that the optimal control for the RBM-problem $u_h^*(\boldsymbol \omega)$  also converges (in expectation) to the optimal control of the original problem $u^*$ when $h \rightarrow 0$. 

\begin{result} \label{res:controls}
Let $u_h^*(\boldsymbol \omega,t)$ be the minimizer of $J_h(\boldsymbol \omega,\cdot)$ in \eqref{eq:Jtilde} and $u^*(t)$ be the minimizer of $J$ in \eqref{eq:J}, then
\begin{equation}
\mathbb{E}[\lvert u_h^* - u^*\rvert_{L^2(0,T; \mathbb{R}^q)}^2] \leq C_{[A,B,x_0,Q,R,x_d,T]} h \mathrm{Var}[\mathcal{A}]. \label{eq:corr_controls}
\end{equation}
\end{result} 

The fourth main result follows directly from Theorem \ref{thm:controls} in Subsection \ref{ssec:controls}. 

The fourth main result justifies the use of the optimal control $u^*_h(\boldsymbol \omega)$, that is optimized for the RBM-dynamics to control the original dynamics, as proposed in \cite{ko2020}. An almost immediate corollary of Main result \ref{res:controls} is that the trajectories of the original dynamics \eqref{eq:dyn_x} resulting from the controls $u^*_h(\boldsymbol \omega,t)$ and $u^*(t)$ will also be close to each other, see Corollary \ref{corr:xstar} in Subsection \ref{ssec:controls}. This further justifies the strategy in \cite{ko2020}.

When the control $u^*_h(\boldsymbol \omega)$ is close to the control $u^*$ that is optimal for the original dynamics, the performance $J(u^*_h(\boldsymbol \omega))$ should also be close to the optimal performance $J(u^*)$. This idea is formalized by the fifth and last main result. 

\begin{result} \label{res:suboptimality}
Let $u^*(t)$ be the control that minimizes the cost functional $J(u)$ in \eqref{eq:J} and let $u_h^*(\boldsymbol \omega,t)$ be the control that minimizes the cost functional $J_h(\boldsymbol \omega, u)$ in \eqref{eq:Jtilde}. Then
\begin{equation}
\mathbb{E}[\lvert J(u^*_h) - J(u^*) \rvert ] \leq C_{[A,B,x_0,Q,R,x_d,T]}h \mathrm{Var}[\mathcal{A}]. \label{eq:res_suboptimality}
\end{equation}
\end{result}

The fifth main result is identical to Corollary \ref{corr:suboptimality} in Subsection \ref{ssec:controls}. Main result \ref{res:suboptimality} is proven as a corollary of Main result \ref{res:controls}/Theorem \ref{thm:controls}. 

The fifth main result is particularly important because it shows that the performance $J(u^*_h(\boldsymbol \omega))$ obtained with control $u^*_h(\boldsymbol \omega)$ optimized for the randomized dynamics is close to the optimal performance $J(u^*)$ when $h\mathrm{Var}[\mathcal{A}]$ is sufficiently small. This further motivates strategies in which the original system is controlled by a control $u^*_h(\boldsymbol \omega)$ that is optimized for the randomized dynamics, as was proposed in \cite{ko2020}. 
 
\subsection{Further examples for \texorpdfstring{$\mathrm{Var}[\mathcal{A}]$}{Var[A]} {\color{red} and computational cost}}
The quantity $\mathrm{Var}[\mathcal{A}]$ describes how the derived estimates depend on the decomposition of $A$ into submatrices and the selected probabilities $p_1, p_2, \ldots , p_{2^M}$. 
We therefore present two other examples that illustrate how $\mathrm{Var}[\mathcal{A}]$ depends on the decomposition of $A$ into submatrices $A_m$ and the selected probabilities $p_{\omega}$. 

\begin{example} \label{ex:A123}
We decompose the matrix $A$ into $M = 3$ parts $A = A_1 + A_2 + A_3$ and consider two choices for the probabilities $p_\omega$. In the first case, we only use one of the submatrices $A_m$ simultaneously. We thus assign probabilities $p_1 = p_2 = p_3 = \tfrac{1}{3}$ to the subsets $S_1 = \{ 1\}$, $S_2 = \{ 2 \}$, and $S_3 = \{ 3 \}$ and zero probability to the other 5 subsets of $\{1,2,3 \}$. We then have that $\pi_1 = \pi_2 = \pi_3 = \tfrac{1}{3}$ and the variance $\mathrm{Var}[\mathcal{A}]$ in \eqref{eq:def_varA} becomes
\begin{align}
&\mathrm{Var}[\mathcal{A}] = \| 3A_1 - A \|^2p_2 + \| 3A_2 - A \|^2p_3 + \| 3A_3 - A \|^2p_4 \nonumber \\
&= \tfrac{1}{3} \left( \| 2 A_1 - A_2 - A_3 \|^2 + \| 2 A_2 - A_1 - A_3 \|^2 + \| 2 A_3 - A_1 - A_2 \|^2\right). \label{eq:ex_A123_1}
\end{align}
In the second case, we always use two of the three submatrices $A_m$ simultaneously. We thus assign probabilities $p_4 = p_5 = p_6 = \tfrac{1}{3}$ to the subsets $S_4 = \{ 1,2\}$, $S_5 = \{ 2,3 \}$, and $S_6 = \{ 1,3 \}$ and zero probability to the other 5 subsets of $\{1,2,3 \}$. We then have that $\pi_1 = p_4 + p_6$, $\pi_2 = p_4 + p_5$, and $\pi_3 = p_5 + p_6$, so that $\pi_1 = \pi_2 = \pi_3 = \tfrac{2}{3}$. The variance $\mathrm{Var}[\mathcal{A}]$ in \eqref{eq:def_varA} becomes
\begin{align}
&\mathrm{Var}[\mathcal{A}] = \| \tfrac{3}{2}(A_1+A_2) - A \|^2p_5 + \| \tfrac{3}{2}(A_2+A_3) - A \|^2p_6 + \| \tfrac{3}{2}(A_1+A_3) - A \|^2p_7 \nonumber \\
&= \tfrac{1}{3} \left( \| \tfrac{1}{2} (A_1 + A_2) - A_3 \|^2 + \| \tfrac{1}{2} (A_2 + A_3) - A_1 \|^2 + \| \tfrac{1}{2} (A_1 + A_3) - A_2 \|^2\right). \label{eq:ex_A123_2} 
\end{align}
Observe that $\| \tfrac{1}{2} (A_1 + A_2) - A_3 \|^2 = \tfrac{1}{4} \| 2 A_3 - A_1 - A_2 \|^2$ and that similar expressions relate the other terms in \eqref{eq:ex_A123_1} and \eqref{eq:ex_A123_2}. The variance for the first case in \eqref{eq:ex_A123_1} is thus four times larger than the variance for the second case in \eqref{eq:ex_A123_2}. Increasing the overlap between the possible values of $\mathcal{A}_h(\boldsymbol \omega,t)$ thus reduces $\mathrm{Var}[\mathcal{A}]$ and will improve the accuracy of the proposed method. It is worth noting that similar observations have been made for domain decomposition methods, for which it is well-known that increasing the overlap between subdomains increases the convergence rate (see, e.g., \cite[Section 1.5]{dolean2015}). Note however that increasing the overlap will also reduce the sparsity of $\mathcal{A}_h(t)$ and thus also increase the computational cost. This will be illustrated further in Example \ref{ex:comp_cost} and the numerical examples in Section \ref{sec:examples}.
\end{example}

\begin{example} \label{ex:popt}
It is not always optimal to choose the probabilities uniform. To illustrate this, we assume $A=A_1+A_2$ has a block-diagonal decomposition
\begin{equation}
A = \begin{bmatrix}
A_{11} & 0 \\ 
0 & A_{22}
\end{bmatrix}, \qquad A_1 = \begin{bmatrix}
A_{11} & 0 \\ 0 & 0
\end{bmatrix}, \qquad A_2 = \begin{bmatrix}
0 & 0 \\
0 & A_{22}
\end{bmatrix}.
\end{equation}
It easy to verify that $\| \alpha A_1 + \beta A_2 \| = \max\{ \lvert \alpha \rvert \| A_1\|, \lvert \beta \rvert \| A_2\| \}$ for any $\alpha,\beta \in \mathbb{R}$. We assign the (at this point undetermined) probability $p_1 = p$ to the subset $S_1 = \{ 1\}$, the probability $p_2 = 1-p$ to the subset $S_2 = \{ 2\}$, and probabilities $p_3 = p_4 = 0$ to the subsets $S_3 = \emptyset$ and $S_4 = \{1,2 \}$. It follows that $\pi_1 = p$ and $\pi_2 = 1-p$ and that 
\begin{align}
\mathrm{Var}[\mathcal{A}] &= \| \tfrac{1}{p}A_1 - A \|^2 p + \| \tfrac{1}{1-p} A_2 - A \|^2 (1-p) \nonumber \\
&= \| \tfrac{1}{p}((1-p)A_1 - pA_2) \|^2 p + \| \tfrac{1}{1-p} (pA_2 - (1-p)A_1) \|^2 (1-p) \nonumber \\
&= \| (1-p)A_1 - pA_2 \|^2 \left(\tfrac{1}{p} + \tfrac{1}{1-p} \right)
= \left\| \sqrt{\tfrac{1-p}{p}}A_1 + \sqrt{\tfrac{p}{1-p}}A_2 \right\|^2 \nonumber \\
&= \left( \max\left\{\sqrt{\tfrac{1-p}{p}}\|A_1\|, \sqrt{\tfrac{p}{1-p}}\| A_2 \|  \right\} \right)^2. \label{eq:ex3_VarA}
\end{align}
It is now easy to see that $\mathrm{Var}[\mathcal{A}]$ is minimal when $\sqrt{\tfrac{1-p}{p}}\|A_1\| = \sqrt{\tfrac{p}{1-p}}\| A_2 \|$. Solving this equation for $p$, we find optimal probability
\begin{equation}
p^* = \frac{\| A_1\|}{\| A_1\| + \| A_2 \|}. \label{eq:ex3_pstar}
\end{equation}
We observe that the larger the submatrix $A_1$ is compared to $A_2$, the larger the probability $p$ with which the submatrix $A_1$ is selected should be. 
Inserting the optimal probability $p^*$ in \eqref{eq:ex3_pstar} into the expression for $\mathrm{Var}[\mathcal{A}]$, we find that
\begin{equation}
\mathrm{Var}[\mathcal{A}]^* = \| A_1 \| \| A_2 \|. 
\end{equation}
With uniform probabilities, i.e., with $p = 1/2$, $\mathrm{Var}[\mathcal{A}] = \max\{ \| A_1 \|^2, \| A_2 \|^2 \}$, see \eqref{eq:ex3_VarA}. When $\|A_1 \| / \| A_2 \| \gg 1$ or $\|A_1 \| / \| A_2 \| \ll 1$, using the optimal probability $p^*$ in \eqref{eq:ex3_pstar} can thus reduce $\mathrm{Var}[\mathcal{A}]$ significantly. 
\end{example}

{ %\color{red}
We conclude this section with two examples that illustrate the potential reduction in computational cost offered by the proposed randomized time-splitting method.

\begin{example} \label{ex:comp_cost}
Let $A \in \mathbb{R}^{N\times N}$ be a sparse symmetric negative semi-definite matrix with a bandwidth $b$, i.e.\ $[A]_{ij} = 0$ when $\lvert i-j \rvert > b$. Select $n_1,n_2,n_3 \in \{1,2,\ldots,N \}$ such that $n_1 + n_2 + n_3 = N + 2b$. It is then possible to split $A$ as $A = A_1 + A_2 + A_3$ with
\begin{equation}
A_1 = \begin{bmatrix}
A_{11} & 0 \\ 0 & 0
\end{bmatrix}, \qquad A_2 = \begin{bmatrix}
0_{n_1-b} & 0 & 0 \\
0 & A_{22} & 0 \\
0 & 0 & 0_{n_3-b}
\end{bmatrix}, \qquad A_3 = \begin{bmatrix}
0 & 0 \\ 0 & A_{33}
\end{bmatrix},
\end{equation}
where $A_{11} \in \mathbb{R}^{n_1 \times n_1}$, $A_{22} \in \mathbb{R}^{n_2 \times n_2}$, $A_{33} \in \mathbb{R}^{n_3 \times n_3}$, $0_n$ denotes an $n \times n$ zero matrix, and the $0$'s denote zero matrices of appropriate size. We assign probabilities $p_1 = p_2 = p_3 = \tfrac{1}{3}$ to the subsets $S_1 = \{ 1\}$, $S_2 = \{ 2 \}$, and $S_3 = \{ 3 \}$ and zero probability to the other 5 subsets of $\{1,2,3 \}$. The computational cost for one time step with the matrix $A_1$ is $O(n_1^r)$, where $r \in [1,3]$ is a certain power that depends on $b$, the time discretization scheme, and the method used to solve the resulting linear systems. In particular, $r = 1$ when $A$ is tridiagonal (i.e.\ when $b = 1$), $r = 3$ for an implicit time discretization scheme in which the resulting linear systems are solved by Gaussian elimination, and $r = 2$ for an implicit time discretization scheme in which the resulting linear systems are solved based on a precomputed Lower-Upper (LU) factorization. Similarly, the computational cost for one time step with the matrices $A_2$ or $A_3$ or with the full matrix $A$ is $O(n_2^r)$ or $O(n_3^r)$ or $O(N^r)$, respectively. The proposed randomized time-splitting scheme is therefore expected to reduce the computational cost for one forward simulation (on the same temporal grid) by a factor
\begin{equation}
\frac{p_1n_1^r + p_2n_2^r + p_3n_3^r}{N^r}. 
\end{equation}
When $b \ll N$, it is possible to choose $n_1 \approx n_2 \approx n_3 \approx N/3$, and the reduction in computational cost is then $\approx 1/3^r$. Note that the expected reduction in computational cost can only be observed when $n_1$, $n_2$, and $n_3$ are sufficiently large. As explained in Section \ref{sec:intro}, we expect that the computation of optimal controls is sped up by the same factor as the forward simulation. 

Similarly as in the second case in Example \ref{ex:A123}, we also consider the situation in which the overlap is increased. We thus assign probabilities $p_4 = p_5 = p_6 = \tfrac{1}{3}$ to the subsets $S_4 = \{ 1,2\}$, $S_5 = \{ 2,3 \}$, and $S_6 = \{ 1,3 \}$ and zero probability to the other 5 subsets of $\{1,2,3 \}$. The cost of doing one time step with the matrices $A_1 + A_2$, $A_2 + A_3$, or $A_1 + A_3$ is then proportional to $(n_1 + n_2 -b)^r$, $(n_2 + n_3 - b)^r$, or $(n_1 + n_3)^r$, respectively. When $b \ll N$ and $n_1 \approx n_2 \approx n_3 \approx N/3$ the proposed randomized time-splitting scheme thus reduces the expected computational cost by a factor $2^r / 3^r$. Increasing the overlap thus increases the expected computational cost of the randomized time splitting method by a factor $2^r$, but it also reduces $\mathrm{Var}[\mathcal{A}_h]$ by a factor $4$, see Example \ref{ex:A123}. Choosing the level of overlap is thus a trade-off between accuracy and computational cost. 
\end{example}

\begin{example}
When $A \in \mathbb{R}^{N \times N}$ is symmetric but not sparse, we can select $n_1, n_2, n_3 \in \{1,2, \ldots, N \}$ such that $n_1 + n_2 + n_3 = N$, and split $A$ as $A = A_1 + A_2 + \ldots + A_6$ with
\begin{align}
&A_1 = \begin{bmatrix}
A_{11} & 0 & 0 \\ 0 & 0 & 0 \\ 0 & 0 & 0
\end{bmatrix}, \qquad \quad 
& &A_2 = \begin{bmatrix}
0 & 0 & 0 \\ 0 & A_{22} & 0 \\ 0 & 0 & 0
\end{bmatrix}, \qquad \quad
& &A_3 = \begin{bmatrix}
0 & 0 & 0 \\ 0 & 0 & 0 \\ 0 & 0 & A_{33}
\end{bmatrix}, \nonumber \\
&A_4 = \begin{bmatrix}
0 & A_{12} & 0 \\ A_{21} & 0 & 0 \\ 0 & 0 & 0
\end{bmatrix}, \quad \quad 
& &A_5 = \begin{bmatrix}
0 & 0 & 0 \\ 0 & 0 & A_{23} \\ 0 & A_{32} & 0
\end{bmatrix}, \quad \quad
& &A_6 = \begin{bmatrix}
0 & 0 & A_{13} \\ 0 & 0 & 0 \\ A_{31} & 0 & 0
\end{bmatrix},
\end{align}
where $A_{11} \in \mathbb{R}^{n_1 \times n_1}$, $A_{22} \in \mathbb{R}^{n_2 \times n_2}$, and $A_{33} \in \mathbb{R}^{n_3 \times n_3}$. The cost for doing one time step with $A_1$, $A_2$, or $A_3$ is $O(n_1^r)$, $O(n_2^r)$, and $O(n_3^r)$, respectively, with $r$ as in Example \ref{ex:comp_cost}. Similarly, the cost for doing one time step with $A_4$, $A_5$, or $A_6$ is $O((n_1+n_2)^r)$, $O((n_2+n_3)^r)$, and $O((n_1+n_3)^r)$, respectively. When we assign probabilities $\tfrac{1}{6}$ to the six singleton subsets of $\{1,2,\ldots, 6 \}$ and zero probability to the other, the proposed randomized time-splitting scheme is expected to reduce the computational cost for one forward simulation (on the same temporal grid) by a factor
\begin{equation}
\frac{n_1^r + n_2^r + n_3^r + (n_1+n_2)^r + (n_2+n_3)^r + (n_1+n_3)^r}{6N^r} \approx \frac{1}{2}\left( \frac{1}{3^r} + \frac{2^r}{3^r} \right), 
\end{equation}
where the latter approximation holds when $n_1 \approx n_2 \approx n_3 \approx N/3$. 
\end{example}
}

\section{Convergence analysis} \label{sec:convergence}

The proof of convergence for the RBM optimal control problem is divided into several stages. 

In the first stage, we consider a control $u \in L^2(0,T; \mathbb{R}^q)$ that does not depend on the selected indices $\boldsymbol \omega$. We then show that the expected difference between the RBM-dynamics \eqref{eq:dyn_xtilde} and the original dynamics \eqref{eq:dyn_x} can be bounded in terms of $h \mathrm{Var}[\mathcal{A}]$ as in Main result \ref{res:este0}. The proof of this statement is inspired by the results for interacting particles systems in  \cite{jin2020,jin2020convergence}. 

Because we will also need to deal with the optimal control $u^*_h(\boldsymbol \omega,t)$ that minimizes $J_h(\boldsymbol \omega, \cdot)$, we consider a general family of controls $u_h(\boldsymbol \omega,t)$ (with $\boldsymbol \omega \in \Omega^K$) in the second stage. This is a nontrivial extension of the results in the previous stage because the crucial idea in the proof for the first stage and in \cite{jin2020,jin2020convergence} is that the solutions $x(t_{k-1})$ and $x_h(\boldsymbol \omega,t_{k-1})$ do not depend on $\omega_k$ (the index that is used in the time interval $[t_{k-1}, t_k)$). This is clearly no longer the case when we insert an input $u_h(\boldsymbol \omega,t)$ that depends on $\boldsymbol \omega$, so also on $\omega_k$, into the dynamics \eqref{eq:dyn_x} and \eqref{eq:dyn_xtilde}. This problem is particularly clear when we consider the family of optimal controls $u_h^*(\boldsymbol \omega)$ for which $u_h^*(\boldsymbol \omega,t_{k-1})$ will depend on the choices for the `future' indices $\omega_k, \omega_{k+1}, \ldots \omega_K$. 

In the third stage, we prove the no-gap condition presented in Main result \ref{res:nogap}. A crucial result for the proof is an auxiliary lemma (Lemma \ref{lem:dJuh}) that bounds the differences $J_h(\boldsymbol \omega,u) - J(u)$ and $J_h(\boldsymbol \omega,u_h(\boldsymbol \omega)) - J(u_h(\boldsymbol \omega))$ (in expectation). For controls $u$ that do not depend on $\omega$, a bound on $J_h(\boldsymbol \omega,u) - J(u)$ can be obtained directly from Main result \ref{res:este0}. For controls $u_h(\boldsymbol \omega)$ that do depend on $\boldsymbol \omega$, we need to use Main result \ref{res:este2} to find the bound on the expected difference $J_h(\boldsymbol \omega,u_h(\boldsymbol \omega)) - J(u_h(\boldsymbol \omega))$. For brevity, Lemma \ref{lem:dJuh} considers controls $u_h(\boldsymbol \omega)$ that depend on $\boldsymbol \omega$ (which of course also covers the case in which the control does not depend on $\boldsymbol \omega$). The no-gap condition (i.e., a bound on $J_h(u_h^*(\boldsymbol \omega)) - J(u^*)$) can then be obtained using classical arguments from the calculus of variations and Lemma \ref{lem:dJuh} applied to the optimal controls $u^*$ and $u_h^*(\boldsymbol \omega)$. 

In the fourth stage, we bound the difference between the RBM-optimal control $u^*_h(\boldsymbol \omega)$ and the control $u^*$ optimized for the original dynamics. To this end, we first bound the expected difference between the gradients of $J_h(\boldsymbol \omega, \cdot)$ and $J$. The bound on the difference in the optimal controls then follows from classical arguments based on the $\alpha$-convexity of the functional $J_h(\boldsymbol \omega,\cdot)$. Finally, the bound for the difference $J(u^*_h(\boldsymbol \omega)) - J(u^*)$ follows easily from the previously derived bound on the difference between the optimal controls $u^*_h(\boldsymbol \omega)$ and $u^*$. 

The four stages discussed above will be proved in detail in Subsections \ref{ssec:forward_dynamics_deterministic}--\ref{ssec:controls}. We first present some preliminaries in Subsection \ref{ssec:preliminaries}. 

\subsection{Preliminaries} \label{ssec:preliminaries}
We will use the following notation. The transpose of a real column vector $x$ is denoted by $x^\top$. Similarly, the transpose of a real matrix $A$ is denoted by $A^\top$. The entry in the $i$-th row and $j$-th column of $A$ is denoted by $[A]_{ij}$.  The standard Euclidean innerproduct of two vectors $x,y \in \mathbb{R}^N$ is denoted by $\langle x, y \rangle := x^\top y$. The corresponding norm is denoted by $\lvert x \rvert := \sqrt{x^\top x}$. The (operator) norm of a matrix $A \in \mathbb{R}^{N \times N}$ is denoted by
\begin{equation}
\| A \| := \max_{\lvert x \rvert =1} \lvert Ax \rvert. \label{eq:operator_norm}
\end{equation}
We use $C_{[a,b, \ldots, d]}$ to denote a constant that only depends on the parameters $a,b, \ldots , d$. The value of $C_{[a,b, \ldots, d]}$ may vary from line to line. The $L^p$-norm of a function in $u \in L^p(0,T; \mathbb{R}^q)$ (for $1 \leq p < \infty$ and $p = \infty$) is defined as
\begin{equation}
\lvert u \rvert_{L^p(0,T;\mathbb{R}^q)} := \sqrt[p]{\int_0^T \lvert u(t) \rvert^p \ \mathrm{d}t  }, \qquad \lvert u \rvert_{L^\infty(0,T;\mathbb{R}^q)} :=  \underset{t \in [0,T]}{\mathrm{ess\ sup}} \ \vert u(t) \rvert.
\end{equation}

We now set up the precise probabilistic setting for our problem. The set $\Omega^K$ defined in \eqref{eq:def_Omega} is the natural sample space for the considered problem. To turn $\Omega^K$ into a probability space, we assign  a probability $p(\boldsymbol \omega)$ to each $\boldsymbol \omega \in \Omega^K$ according to
\begin{equation}
p(\boldsymbol \omega) = p_{\omega_1} p_{\omega_2} \ldots p_{\omega_K}. 
\end{equation}
Note that we use here that each index $\omega_k$ is chosen independently from the other indices $\omega_1, \omega_2, \ldots, \omega_{k-1}, \omega_{k+1}, \omega_{k+1}, \ldots , \omega_K$. 

A random element on the sample space $\Omega^K$ is a function $X : \Omega^K \rightarrow V$ from the sample space $\Omega^K$ to a vector space $V$. When $V = \mathbb{R}$, $X : \Omega^K \rightarrow \mathbb{R}$ is also called a random variable. Note that we can embed $V$ into $V^{\Omega^K}$ by associating to each element $x \in V$ the constant function $X(\boldsymbol \omega) = x$ for all $\boldsymbol \omega \in {\Omega^K}$. Constant functions $X(\boldsymbol \omega) = x$ are called deterministic. Functions $X(\boldsymbol \omega)$ that are not deterministic are called stochastic. The expectation operator $\mathbb{E}$ assigns to a random variable $X \in V^{\Omega^K}$ an element of the vector space $V$ 
\begin{align}
\mathbb{E}[X] &= \sum_{\boldsymbol \omega \in \Omega^K} X(\boldsymbol \omega) p(\boldsymbol \omega) \nonumber \\
&= \sum_{\omega_1 =1}^{2^M} \sum_{\omega_2=1}^{2^M} \cdots \sum_{\omega_K =1}^{2^M}  X(\omega_1, \omega_2, \ldots, \omega_K) p_{\omega_1} p_{\omega_2} \cdots p_{\omega_K}. \label{eq:def_expectation}
\end{align}
It is immediate from this definition that $\mathbb{E}$ is linear. When $V = \mathbb{R}$, we also see that $\mathbb{E}[X] \geq 0$ when $X(\boldsymbol \omega) \geq 0$ for all $\boldsymbol \omega \in \Omega^K$. 

Several random elements appear in the randomized splitting method outlined in Subsection \ref{ssec:method}. One example is the matrix $\mathcal{A}_h(\boldsymbol \omega,t)$ defined in \eqref{eq:defAt}. When $t \in [t_{k-1}, t_k)$, $\mathcal{A}_h(\boldsymbol \omega, t)$ only depends on $\omega_k$. Therefore,  the definitions in \eqref{eq:def_expectation} and  \eqref{eq:defAt} show that (for $t \in [t_{k-1}, t_k)$)
\begin{align}
\mathbb{E}[\mathcal{A}_h(t)] &= \sum_{\boldsymbol \omega \in \Omega^K} \mathcal{A}_h(\boldsymbol \omega,t) p(\boldsymbol \omega) =  \sum_{\omega_1 =1}^{2^M} \sum_{\omega_2=1}^{2^M} \cdots \sum_{\omega_K =1}^{2^M} \sum_{m \in S_{\omega_k}} \frac{A_m}{\pi_m}
 p_{\omega_1} p_{\omega_2} \cdots p_{\omega_K} \nonumber \\
 &= \sum_{\omega_k=1}^{2^M} \sum_{m \in S_{\omega_k}} \frac{A_m}{\pi_m}
 p_{\omega_k} = A, \label{eq:expAht}
\end{align}
where the second to last identity follows from \eqref{eq:sump} and the last identity from \eqref{eq:expectationA_step1}. Again using that $\mathcal{A}_h(\boldsymbol \omega,t)$ only depends on $\omega_k$ for $t \in [t_{k-1}, t_k)$, we also find that
\begin{align}
\mathbb{E}[\| \mathcal{A}_h(t) - A \|^2] &= \sum_{\boldsymbol \omega \in \Omega^K} \| \mathcal{A}_h(\boldsymbol \omega,t) - A \|^2 p(\boldsymbol \omega) \nonumber \\
&= \sum_{\omega_k=1}^{2^M} \left\| \sum_{m \in S_{\omega_k}} \frac{A_m}{\pi_m} - A \right\|^2 p_{\omega_k} = \mathrm{Var}[\mathcal{A}], \label{eq:varA}
\end{align}
where the last identity follows from the definition of $\mathrm{Var}[\mathcal{A}]$ in \eqref{eq:def_varA}. Note that \eqref{eq:varA} holds for every time instant $t$ and that $\mathbb{E}[\| \mathcal{A}_h(t) - A \|^2]$ therefore does not depend on the considered time instant $t$. 

Another random element is the solution $x_h : \Omega^K \rightarrow L^2(0,T; \mathbb{R}^N)$ in \eqref{eq:dyn_xtilde}. 
We will frequently use that $\lvert x_h(\boldsymbol \omega,t) \rvert$ can be bounded as follows. First of all, observe that
\begin{equation}
\frac{d}{dt} \lvert x_h(\boldsymbol \omega,t) \rvert^2 = 2\langle x_h(\boldsymbol \omega,t), \mathcal{A}_h(\boldsymbol \omega,t) x_h(\boldsymbol \omega,t) + Bu(t) \rangle \leq 2 \lvert x_h(\boldsymbol \omega,t) \rvert  \lvert B u(t)\rvert, 
\end{equation}
where is was used that $\langle x, \mathcal{A}_h(\boldsymbol \omega,t) x \rangle \leq 0$ for all $x\in \mathbb{R}^N$ and $\boldsymbol \omega \in \Omega^K$ because of Assumption \ref{ass:dissipativity}. 
Now observe that
\begin{equation}
\frac{d}{dt} \lvert x_h(\boldsymbol \omega,t) \rvert = \frac{1}{2 \lvert x_h(\boldsymbol \omega,t)\rvert} \frac{d}{dt} \lvert x_h(\boldsymbol \omega,t) \rvert^2 \leq \lvert Bu(t) \rvert,
\end{equation}
from which we conclude that
\begin{equation}
\lvert x_h(\boldsymbol \omega) \rvert_{L^\infty(0,T;\mathbb{R}^N)} \leq \lvert x_0 \rvert + \lvert Bu \rvert_{L^1(0,T;\mathbb{R}^N)}. \label{eq:boundxh}
\end{equation}
For $x(t)$, a similar derivation shows that
\begin{equation}
\lvert x \rvert_{L^\infty(0,T;\mathbb{R}^)} \leq  \lvert x_0 \rvert + \lvert Bu \rvert_{L^1(0,T;\mathbb{R}^N)}. \label{eq:boundx}
\end{equation}

We will also consider situations in which we apply an input $u_h(\boldsymbol \omega,t)$ to the dynamics \eqref{eq:dyn_x} and \eqref{eq:dyn_xtilde} that depends on $\boldsymbol \omega$. The resulting solutions are then both random elements $x(\boldsymbol \omega,t)$ and $x_h(\boldsymbol \omega,t)$ which satisfy
\begin{align}
\dot{x}(\boldsymbol \omega,t) &= Ax(\boldsymbol \omega,t) + Bu_h(\boldsymbol \omega,t), \qquad \qquad & x(\boldsymbol \omega,0) = x_0, \label{eq:dyn_x_om} \\
\dot{x}_h(\boldsymbol \omega,t) &= \mathcal{A}_h(\boldsymbol \omega,t)x_h(\boldsymbol \omega,t) + Bu_h(\boldsymbol \omega,t), \qquad \qquad & x_h(\boldsymbol \omega,0) = x_0, \label{eq:dyn_xtilde_om}
\end{align}
In this case we can obtain estimates similar to \eqref{eq:boundxh} and \eqref{eq:boundx} with $u$ and $x$ replaced by $u_h(\boldsymbol \omega)$ and $x(\boldsymbol \omega)$, respectively. 

The third important random element in this paper is the optimal control $u^*_h(\boldsymbol \omega, \cdot)$ that minimizes $J_h(\boldsymbol \omega, \cdot)$ in \eqref{eq:Jtilde}. The coercivity of the functional $J_h(\boldsymbol \omega,\cdot)$ allows us to bound $\lvert u^*_h(\boldsymbol \omega) \rvert_{L^2(0,T;\mathbb{R}^q)}$ as follows. Denote the smallest eigenvalue of the matrix $R$ by $\alpha > 0$, then
\begin{equation}
\frac{\alpha}{2} \lvert u^*_h(\boldsymbol \omega) \rvert^2_{L^2(0,T;\mathbb{R}^q)} \leq \frac{1}{2}\int_0^T u_h^*(t)^\top R u_h^*(t) \ \mathrm{d}t \leq J_h(\boldsymbol \omega, u^*_h(\boldsymbol \omega)) \leq J_h(\boldsymbol \omega, 0), \label{eq:bound_uhopt_step1}
\end{equation}
where the last inequality follows because $u^*_h(\boldsymbol \omega)$ is the minimizer of $J_h(\boldsymbol \omega,\cdot)$. Next, observe that
\begin{align}
J_h(\boldsymbol \omega, 0) &\leq \frac{1}{2}\int_0^T (x_h(\boldsymbol \omega,t)-x_d(t))^\top Q (x_h(\boldsymbol \omega,t)-x_d(t)) \ \mathrm{d}t \nonumber \\
&\leq \tfrac{1}{2}\| Q \| \left(  \lvert x_h(\boldsymbol \omega) \rvert_{L^2(0,T;\mathbb{R}^N)} + \lvert x_d \rvert_{L^2(0,T;\mathbb{R}^N)} \right)^2 \nonumber \\
&\leq \tfrac{1}{2}\| Q \| \left( T  \lvert x_0 \rvert + \lvert x_d \rvert_{L^2(0,T;\mathbb{R}^N)} \right)^2 = C_{[x_0, Q, x_d,T]}, \label{eq:bound_uhopt_step2}
\end{align}
where $x_h(\boldsymbol \omega,t)$ denotes the solution of \eqref{eq:dyn_xtilde} with $u(t) = 0$ and the last inequality follows from \eqref{eq:boundxh}. Looking back at \eqref{eq:bound_uhopt_step1}, we find
\begin{equation}
\lvert u^*_h(\boldsymbol \omega) \rvert^2_{L^2(0,T;\mathbb{R}^N)} \leq C_{[x_0, Q, R, x_d, T]}. \label{eq:bound_uhopt}
\end{equation}

Finally, we repeat some standard definitions from the theory of the convex optimization, see, e.g., \cite{minoux1986}. A functional $J : V \rightarrow \mathbb{R}$ on a normed vector space $V$ is $\alpha$-convex if there exists an $\alpha \geq 0$ such that for all $u, v \in V$ and $\theta \in [0, 1]$
\begin{equation}
J((1 -\theta)u + \theta v) \leq (1 -\theta)J(u) + \theta J(v) - \tfrac{\alpha}{2} \theta (1-\theta) \lvert u -v \rvert_V^2. \label{eq:def_convex}
\end{equation}
One can easily verify that the functional $J_h(\boldsymbol \omega, \cdot)$ is $\alpha$-convex (for all $\boldsymbol \omega \in \Omega^K$) when we take $\alpha$ as the smallest eigenvalue of the positive definite matrix $R$. 
The G\^ateaux-derivative of $J$ at the point $u$ in the direction $v$ is denoted by $\delta J(u; v)$, i.e.\
\begin{equation}
\delta J(u; v) := \lim_{h \rightarrow 0} \frac{J(u+hv) - J(u)}{h}.
\end{equation}
By subtracting $J(u)$ from both sides of \eqref{eq:def_convex}, dividing the resulting inequality by $\theta$, and then taking the limit $\theta \rightarrow 0$, we find the well-known inequality
\begin{equation}
 J(v) \geq J(u) + \delta J(u; v - u) + \tfrac{\alpha}{2}  \lvert v-u \rvert_V^2.  \label{eq:gateaux}
\end{equation}

\subsection{The forward dynamics with a deterministic input} \label{ssec:forward_dynamics_deterministic}

In this subsection, we consider a deterministic $u(t)$ and derive a bound for the error
\begin{equation}
e_h(\boldsymbol \omega,t) := x_h(\boldsymbol \omega,t) - x(t), \label{eq:def_e}
\end{equation} 
where $x_h(\boldsymbol \omega,t)$ and $x(t)$ are the solutions of \eqref{eq:dyn_xtilde} and \eqref{eq:dyn_x} resulting from the same input $u(t)$, respectively. 

\begin{remark} \label{rem:bias_x}
It is important to stress that $x_h(t)$ is not an unbiased estimator for $x(t)$, i.e., we do \emph{not} have $\mathbb{E}[e_h(t)] = \mathbb{E}[x_h(t)] - x(t) = 0$.  This can for example be observed when we write the error dynamics as
\begin{align}
\dot{e}_h(\boldsymbol \omega,t) &=  \mathcal{A}_h(\boldsymbol \omega,t)x_h(\boldsymbol \omega,t) + Bu(t) - Ax(\boldsymbol \omega,t) - Bu(t) \nonumber \\
&= Ae_h(\boldsymbol \omega,t) + (\mathcal{A}_h(\boldsymbol \omega,t) - A)x_h(\boldsymbol \omega,t), \qquad\qquad e_h(\boldsymbol \omega,0) = 0,\label{eq:dyn_e2}
\end{align}
where we have substituted $x(\boldsymbol \omega,t) =x_h(\boldsymbol \omega,t) - e_h(\boldsymbol \omega,t)$. 
Taking the expected value in \eqref{eq:dyn_e2} we find
\begin{equation}
\frac{d}{dt}\mathbb{E}[e_h(t)] = A \mathbb{E}[e_h(t)] +  \mathbb{E}[(\mathcal{A}_h(t) - A)x_h(t)], \qquad\qquad \mathbb{E}[e_h(0)] = 0. \label{eq:dyn_emean}
\end{equation}
However, \eqref{eq:dyn_emean} does not imply that $\mathbb{E}[e_h(t)] = 0$ for all $t$ because generally 
\begin{equation}
\mathbb{E}[(\mathcal{A}_h(t)-A)x_h(t)] \neq \mathbb{E}[\mathcal{A}_h(t) - A] \mathbb{E} [x_h(t)] = 0, 
\end{equation}
where the equality follows because $\mathbb{E}[\mathcal{A}_h(t)] = A$, see \eqref{eq:expAht}. 
This would be the case when $\mathcal{A}_h(\boldsymbol \omega,t)$ and $x_h(\boldsymbol \omega,t)$ are independent, but they are correlated by the dynamics \eqref{eq:dyn_xtilde}.   Note, however, that at the beginning of each time interval $[t_{k-1}, t_k)$, the value of $\mathcal{A}_h(\boldsymbol \omega,t)$ changes and that $\mathcal{A}_h(\boldsymbol \omega,t_{k-1})$ is independent of the values of $\mathcal{A}_h(\boldsymbol \omega,t)$ for $t < t_{k-1}$ so that
\begin{equation}
\mathbb{E}[(\mathcal{A}_h(t_{k-1})-A)x_h(t_{k-1})] = \mathbb{E}[\mathcal{A}_h(t_{k-1}) - A] \mathbb{E} [x_h(t_{k-1})] = 0, 
\end{equation}
where the second identity again follows because $\mathbb{E}[\mathcal{A}_h(t)] = A$, see \eqref{eq:expAht}. 
This observation is crucial to obtain the main result of this subsection. 
\end{remark}

The main result in this subsection is the following. 
\setcounter{theorem}{0}
\begin{theorem} \label{thm:este}
Assume that the input $u(t)$ in \eqref{eq:dyn_xtilde} is deterministic and equal to the input $u(t)$ in \eqref{eq:dyn_x} and that Assumptions \ref{ass:dissipativity} and \ref{ass:positivity} hold, then
\begin{equation}
\mathbb{E}[\lvert e_h(t) \rvert^2] \leq h \mathrm{Var}[\mathcal{A}] (\| A\|t^2 + 2t) (\lvert x_0 \rvert + \lvert Bu \rvert_{L^1(0,T;\mathbb{R}^N}))^2. \label{eq:thm_este}
\end{equation}
\end{theorem}

\begin{proof} Observe that
\begin{align}
\dot{e}_h(\boldsymbol \omega,t) &=  \mathcal{A}_h(\boldsymbol \omega,t)x_h(\boldsymbol \omega,t) + Bu(t) - Ax(\boldsymbol \omega,t) - Bu(t) \nonumber \\
&= \mathcal{A}_h(\boldsymbol \omega,t ) e_h(\boldsymbol \omega,t) + (\mathcal{A}_h(\boldsymbol \omega,t) - A)x(t), \qquad\qquad e_h(\boldsymbol \omega,0) = 0,\label{eq:dyn_e1}
\end{align}
where the last equation follows after substituting $x_h(\boldsymbol \omega,t) = x(\boldsymbol \omega,t) + e_h(\boldsymbol \omega,t)$. 

Fix $t \in [0, T]$ and let $k \leq K$ be such that $t \in [t_{k-1}, t_k)$. 

Consider an arbitrary time instant $s \in [0,t)$ and let $\ell \in \{1,2, \ldots , k \}$ be such that $s \in [t_{\ell-1}, t_\ell)$. Then \eqref{eq:dyn_e1} shows that
\begin{align}
\frac{d}{ds} &\lvert e_h(\boldsymbol \omega,s)\rvert^2 = 2 \langle e_h(\boldsymbol \omega,s), \mathcal{A}_h(\boldsymbol \omega,s) e_h(\boldsymbol \omega,s) \rangle + 2 \langle e_h(\boldsymbol \omega,s), (\mathcal{A}_h(\boldsymbol \omega,s) - A) x(s) \rangle \nonumber \\
&\qquad = 2 \langle e_h(\boldsymbol \omega,s), \mathcal{A}_h(\boldsymbol \omega,s) e_h(\boldsymbol \omega,s) \rangle  + 2 \langle e_h(\boldsymbol \omega,t_{\ell-1}), (\mathcal{A}_h(\boldsymbol \omega,s) - A) x(s) \rangle \nonumber \\
&\qquad\qquad\qquad\qquad \qquad\qquad\qquad  + 2 \langle \Delta e_h(\boldsymbol \omega,s), (\mathcal{A}_h(\boldsymbol \omega,s) - A) x(s) \rangle, \label{eq:lem_este_step0}
\end{align}
where, in the second equality, we have introduced 
\begin{equation}
\Delta e_h(\boldsymbol \omega,s) := e_h(\boldsymbol \omega,s) - e_h(\boldsymbol \omega,t_{\ell-1}). \label{eq:lem_este_deltae_def}
\end{equation}

The first term on the RHS of \eqref{eq:lem_este_step0} is nonpositive due to Assumption \ref{ass:dissipativity}. We thus find after taking the expected value in \eqref{eq:lem_este_step0} that
\begin{align}
\frac{d}{ds} \mathbb{E}[ \lvert e_h(s) \rvert^2 ] &\leq  2 \mathbb{E}[ \langle e_h(t_{\ell-1}), (\mathcal{A}_h(s) - A) x(s) \rangle ] \nonumber \\
& \qquad \qquad + 2\mathbb{E}[\langle \Delta e_h(s), (\mathcal{A}_h(s) - A) x(s) \rangle]. \label{eq:lem_este_step1}
\end{align}

For the first term on the RHS of \eqref{eq:lem_este_step1}, observe that $e_h(\boldsymbol \omega,t_{\ell-1}) = x_h(\boldsymbol \omega,t_{\ell-1}) - x(t_{\ell-1})$ only depends on $\omega_1, \ldots \omega_{\ell-1}$, so that
\begin{align}
& \mathbb{E}[ \langle e_h(t_{\ell-1}), (\mathcal{A}_h(s) - A) x(s) \rangle ]
= \sum_{\boldsymbol \omega \in \Omega^K} \langle e_h(\boldsymbol \omega,t_{\ell-1}), (\mathcal{A}_h(\boldsymbol \omega,s) - A) x(s) \rangle p(\boldsymbol \omega) \nonumber \\
&= \sum_{\omega_1=1}^{2^M} \cdots \sum_{\omega_{\ell-1}=1}^{2^M} \sum_{\omega_\ell=1}^{2^M} \bigg\langle e_h(\boldsymbol \omega, t_{\ell-1}) , \bigg( \sum_{m \in S_{\omega_\ell}} \frac{A_m}{\pi_m} - A \bigg) x(s)
\bigg\rangle p_{\omega_1} \cdots p_{\omega_{\ell-1}}p_{\omega_\ell} \nonumber \\
&= \sum_{\omega_1=1}^{2^M} \cdots \sum_{\omega_{\ell-1}=1}^{2^M} \bigg\langle e_h(\boldsymbol \omega, t_{\ell-1}), \bigg( \sum_{\omega_\ell=1}^{2^M} \sum_{m \in S_{\omega_\ell}} \frac{A_m}{\pi_m}p_{\omega_\ell} - A \bigg) x(s)
\bigg\rangle p_{\omega_1} \cdots p_{\omega_{\ell-1}} \nonumber  \\
&= 0,
\label{eq:lem_este_step2}
\end{align}
where the second identity uses \eqref{eq:sump}, the third identity follows from \eqref{eq:sump} and the fact that $e_h(\boldsymbol \omega, t)$ does not depend on $\omega_\ell$, and the last identity follows because \eqref{eq:expectationA_step1} shows that the factor between round brackets vanishes. 

For the second term on the RHS of \eqref{eq:lem_este_step1}, we use that
\begin{align}
\mathbb{E}[& \langle \Delta e_h(s), (\mathcal{A}_h(s) - A) x(s) \rangle] \leq \mathbb{E}[\lvert \Delta e_h(s)\rvert \|\mathcal{A}_h(s) - A \| \lvert x(s) \rvert] \nonumber \\
&\leq \sqrt{\mathbb{E}[\lvert \Delta e_h(s) \rvert^2] \mathbb{E}[\| \mathcal{A}_h(s) - A\|^2 \lvert x(s) \rvert^2]} = \sqrt{\mathbb{E}[\lvert \Delta e_h(s) \rvert^2]} \sqrt{\mathrm{Var}[\mathcal{A}]} \lvert x(s) \rvert \nonumber \\
&\leq \sqrt{\mathbb{E}[\lvert \Delta e_h(s) \rvert^2]} \sqrt{\mathrm{Var}[\mathcal{A}]} (\lvert x_0 \rvert + \lvert Bu \rvert_{L^1(0,T;\mathbb{R}^N)}),
\label{eq:lem_este_step3}
\end{align}
where the first identity follows from the Cauchy-Schwartz inequality in $\mathbb{R}^N$, the second inequality from Cauchy-Schwartz inequality in the probability space, and the last inequality follows from \eqref{eq:boundx}.

We now claim that
\begin{equation}
\mathbb{E}[\lvert \Delta e_h(s) \rvert^2] \leq h^2 \mathrm{Var}[\mathcal{A}] (\| A\|s + 1)^2(\lvert x_0 \rvert + \lvert Bu \rvert_{L^1(0,T;\mathbb{R}^N)})^2. \label{eq:lem_este_claim}
\end{equation}
We will prove \eqref{eq:lem_este_claim} at the end of the proof. 
Inserting the claim \eqref{eq:lem_este_claim} into \eqref{eq:lem_este_step3}, we find
\begin{align}
\mathbb{E}[& \langle \Delta e_h(s), (\mathcal{A}_h(s) - A) x(s) \rangle] \leq h \mathrm{Var}[\mathcal{A}] (\| A\|s + 1) (\lvert x_0 \rvert +  \lvert Bu \rvert_{L^1(0,T;\mathbb{R}^N)})^2. 
\label{eq:lem_este_step4}
\end{align}
Inserting \eqref{eq:lem_este_step2} and \eqref{eq:lem_este_step4} into \eqref{eq:lem_este_step1} shows that
\begin{equation}
\frac{d}{ds} \mathbb{E}[ \lvert e_h(s)\rvert^2 ] \leq 2 h \mathrm{Var}[\mathcal{A}](\| A\|s + 1) ( \lvert x_0 \rvert + \lvert Bu_d \rvert_{L^1(0,T;\mathbb{R}^N)})^2. \label{eq:lem_este_step5}
\end{equation}
Integrating \eqref{eq:lem_este_step5} from $s = 0$ to $s = t$ using that $e_h(\omega,0) = 0$ now shows that
\begin{equation}
\mathbb{E}[ \lvert e_h(t) \rvert ^2 ] \leq h \mathrm{Var}[\mathcal{A}] (\| A\|t^2 + 2t) (\lvert x_0 \rvert + \lvert Bu_d \rvert_{L^1(0,T;\mathbb{R}^N)})^2,
\end{equation}
which is the desired estimate \eqref{eq:thm_este}. 

It thus remains to show that \eqref{eq:lem_este_claim} holds. Recall that, for $\tau \in [t_{\ell-1}, s)$, \eqref{eq:lem_este_deltae_def} shows that $\Delta e_h(\boldsymbol \omega,\tau) = e_h(\boldsymbol \omega,\tau) - e_h(\boldsymbol \omega,t_{\ell-1})$. Using \eqref{eq:dyn_e2}, we thus see that $\Delta e_h(\boldsymbol \omega,\tau)$ is the solution of the ODE
\begin{equation}
\tfrac{d}{d\tau}\Delta e_h(\boldsymbol \omega,\tau) = \dot{e}_h(\boldsymbol \omega,t) = Ae_h(\boldsymbol \omega,t) + (\mathcal{A}_h(\boldsymbol \omega,t) - A)x_h(\boldsymbol \omega,t),
\end{equation}
with initial condition $\Delta e_h(\boldsymbol \omega,t_{\ell-1}) = 0$. We therefore also have that
\begin{equation}
\frac{d}{d\tau} \lvert \Delta e_h(\boldsymbol \omega,\tau) \rvert= \frac{\langle \Delta e_h(\boldsymbol \omega,\tau), \dot{e}_h(\boldsymbol \omega,\tau) \rangle}{\lvert \Delta e_h(\boldsymbol \omega,\tau) \rvert} \leq \lvert Ae_h(\boldsymbol \omega,\tau) \rvert + \lvert(\mathcal{A}_h(\boldsymbol \omega,\tau) - A)x_h(\boldsymbol \omega,\tau) \rvert. \label{eq:lem_deltae_step1}
\end{equation}
Using that $\Delta e_h(\boldsymbol \omega,t_{\ell-1}) = 0$, integrating \eqref{eq:lem_deltae_step1} from $\tau = t_{\ell-1}$  to $\tau = s$ yields
\begin{align}
\lvert \Delta e_h(\boldsymbol \omega,s) \rvert \leq \int_{t_{\ell-1}}^s \left( \| A \| \lvert e_h(\boldsymbol \omega,\tau) \rvert + \lvert (\mathcal{A}_h(\boldsymbol \omega,\tau) - A)x_h(\boldsymbol \omega,\tau) \rvert \right) \ \mathrm{d}\tau. \label{eq:lem_deltae_step2}
\end{align}
To bound $e_h(\boldsymbol \omega,\tau)$, we apply the variation of constants formula to the error dynamics in \eqref{eq:dyn_e2} and obtain
\begin{align} 
\lvert e_h(\boldsymbol \omega,\tau) \rvert &= \left\lvert \int_0^\tau e^{A(\tau-\sigma)} (\mathcal{A}_h(\boldsymbol \omega,\sigma) - A ) x_h(\boldsymbol \omega,\sigma) \ \mathrm{d}\sigma \right\rvert \nonumber \\
&\leq \int_0^\tau \| \mathcal{A}_h(\boldsymbol \omega,\sigma) - A \| \ \mathrm{d}\sigma \ ( \lvert x_0 \rvert + \lvert Bu \rvert_{L^1(0,T;\mathbb{R}^N)}), \label{eq:lem_deltae_step3}
\end{align}
where we have used the bound for $x_h(\boldsymbol \omega,\sigma)$ in \eqref{eq:boundxh} and that $\| e^{A(\tau-\sigma)}\| \leq 1$ because Assumption \ref{ass:dissipativity} implies that $A$ is dissipative. 
Using this result in \eqref{eq:lem_deltae_step2}, we find
\begin{equation}
\lvert \Delta e_h(\boldsymbol \omega,s) \rvert \leq \int_{t_{\ell-1}}^s g(\boldsymbol \omega,\tau) \ \mathrm{d}\tau\ (\rvert x_0 \lvert + \rvert Bu \lvert_{L^1(0,T;\mathbb{R}^N)}), \label{eq:lem_deltae_step4}
\end{equation}
where we have again used the bound on $x_h(\boldsymbol \omega,t)$ in \eqref{eq:boundxh} for the second term in \eqref{eq:lem_deltae_step2} and introduced
\begin{equation}
g(\boldsymbol \omega,\tau) := \left( \| A \| \int_0^\tau \|\mathcal{A}_h(\boldsymbol \omega,\sigma) - A \| \ \mathrm{d}\sigma + \|\mathcal{A}_h(\boldsymbol \omega,\tau) - A \| \right). \label{eq:lem_deltae_step5}
\end{equation}
Squaring both sides in \eqref{eq:lem_deltae_step4} and taking the expectation, we find
\begin{align}
\mathbb{E}[ \lvert \Delta e_h(s) \rvert^2] &\leq \mathbb{E}\left[ \left( \int_{t_{\ell-1}}^s g(\tau) \ \mathrm{d}\tau \right)^2  \right] (\lvert x_0 \rvert + \lvert Bu \rvert_{L^1(0,T;\mathbb{R}^N)})^2 \nonumber \\
&\leq (s-t_{\ell-1}) \int_{t_{\ell-1}}^s \mathbb{E}[(g(\tau))^2] \ \mathrm{d}\tau\ (\lvert x_0 \rvert + \lvert Bu \rvert_{L^1(0,T;\mathbb{R}^N)})^2, \label{eq:lem_deltae_step6}
\end{align}
where the second inequality follows from the Cauchy-Schwartz inequality in $L^2(t_{\ell-1},s)$. 
Now observe that \eqref{eq:lem_deltae_step5} shows that
\begin{align}
\mathbb{E}[(g(\tau))^2] &= \| A \|^2 \int_0^\tau \int_0^\tau \mathbb{E}[ \|\mathcal{A}_h(\sigma) - A \|\|\mathcal{A}_h(\sigma') - A \| ] \ \mathrm{d}\sigma \ \mathrm{d}\sigma' \nonumber \\ 
& \quad + 2 \| A \| \int_0^\tau \mathbb{E}[ \|\mathcal{A}_h(\sigma) - A \| \|\mathcal{A}_h(\tau) - A \| ] \ \mathrm{d}\sigma + \mathbb{E}[\|\mathcal{A}_h(\tau) - A \|^2]. \label{eq:lem_deltae_step7}
\end{align}
Because $\mathbb{E}[\| \mathcal{A}_h(t) - A \|^2] = \mathrm{Var}[\mathcal{A}]$ for all $t$, we also have that
\begin{equation}
\mathbb{E}[ \|\mathcal{A}_h(\sigma) - A \| \|\mathcal{A}_h(\tau) - A \| ] \leq \sqrt{ \mathbb{E}[ \|\mathcal{A}_h(\sigma) - A \|^2] \mathbb{E}[\|\mathcal{A}_h(\tau) - A \|^2]} = \mathrm{Var}[\mathcal{A}]. \label{eq:lem_deltae_step8}
\end{equation}
Using this result in \eqref{eq:lem_deltae_step7}, we obtain
\begin{equation}
\mathbb{E}[(g(\tau))^2] \leq \mathrm{Var}[\mathcal{A}] (\| A \|\tau + 1)^2.  \label{eq:lem_deltae_step9}
\end{equation}
Using this result again in \eqref{eq:lem_deltae_step6}, also using that $s-t_{\ell-1} \leq h$ and $\tau \leq s$, we find the claimed inequality \eqref{eq:lem_este_claim}. 
\end{proof}

Some remarks regarding Theorem \ref{thm:este} are in order. 

\begin{remark} \label{rem:infdim}
The error estimate in Theorem \ref{thm:este} involves the operator norm of the matrix $A$. This suggests that the expected error $\mathbb{E}[ \lvert e_h(t) \rvert^2]$ grows when we are considering better approximations $A$ of an unbounded operator, which for example happens when we consider a discretization of a PDE and refine the spatial grid. However, Figure \ref{fig:heat1D_convX} in Section \ref{sec:examples} indicates that $\mathbb{E}[\lvert e_h(t) \rvert] \leq C \sqrt{h \mathrm{Var}[\mathcal{A}]}$ for a constant $C$ that does not increase (but even seems to decrease) when the spatial grid is refined. 

A first step in understanding the infinite-dimensional case better is taken in Appendix \ref{app:commutative}, where we prove that
\begin{equation}
\mathbb{E}[\lvert e_h(t) \rvert^2] \leq 2ht \mathrm{Var}_W[\mathcal{A}] \lvert W^{-1}x_0 \rvert. 
\end{equation}
under the additional assumptions that $u(t) \equiv 0$ and that all matrices $A_m$ commute pairwise. Here, $W$ is any invertible matrix and $\mathrm{Var}_W[\mathcal{A}]$ is the weighted variance introduced in Remark \ref{rem:def_varAW}. Observe that the operator norm $\| A \|$ does not appear in this estimate. The result from Appendix \ref{app:commutative} extends naturally to an infinite dimensional setting in which all operators $A_m$ have the same domain $D(A_m) = D(A)$.

Recall from Remark \ref{rem:def_varAW} that a typical choice for $W$ is $W = (A - \lambda I)^{-1}$ for some $\lambda$ in the resolvent of $A$. For $\lvert W^{-1}x_0 \rvert$ to be bounded, we thus require that $x_0 \in D(A)$, where $D(A)$ denotes the domain of the operator $A$. In an infinite dimensional setting we thus need an additional smoothness assumption on the initial condition $x_0$. Such conditions are typical for (deterministic) splitting algorithms, see e.g.\ \cite{hansen2008, ignat2011}. Further details can be found in Appendix  \ref{app:commutative}.
\end{remark}

\begin{remark}
The error estimate in Theorem \ref{thm:este} is derived based on the error dynamics \eqref{eq:dyn_e1}. Considering the error dynamics \eqref{eq:dyn_e2} leads to a less clean proof because instead of the 3 terms on the RHS of \eqref{eq:lem_este_step0}, we then get 4 terms
\begin{align}
\frac{d}{ds} \lvert &e_h(\boldsymbol \omega,s) \rvert^2 = 2 \langle e_h(\boldsymbol \omega,s), A e_h(\boldsymbol \omega,s) \rangle + 2 \langle e_h(\boldsymbol \omega,s), (\mathcal{A}_h(\boldsymbol \omega,s) - A) x_h(\boldsymbol \omega,s) \rangle \nonumber \\
&= 2 \langle e_h(\boldsymbol \omega,s), A e_h(\boldsymbol \omega,s) \rangle + 2 \langle e_h(\boldsymbol \omega,t_{\ell-1}), (\mathcal{A}_h(\boldsymbol \omega,s) - A) x_h(\boldsymbol \omega,t_{\ell-1}) \rangle \nonumber \\
& \qquad\qquad\qquad + 2 \langle \Delta e_h(\boldsymbol \omega,s), (\mathcal{A}_h(\boldsymbol \omega,s) - A) x_h(\boldsymbol \omega,s) \rangle \nonumber \\ & \qquad\qquad\qquad + 2 \langle e_h(\boldsymbol \omega,s),  (\mathcal{A}_h(\boldsymbol \omega,s) - A) \Delta x_h(\boldsymbol \omega,s) \rangle, \label{eq:lem_este_rem1}
\end{align}
where $\Delta e_h(\boldsymbol \omega,s) := e_h(\boldsymbol \omega,s) - e_h(\boldsymbol \omega,t_{\ell-1})$ and $\Delta x_h(\boldsymbol \omega,s) := x_h(\boldsymbol \omega,s) - x_h(\boldsymbol \omega,t_{\ell-1})$. This approach is closer to proofs for interacting particle systems in \cite{jin2020}. 

Note that the fourth term in \eqref{eq:lem_este_rem1} is needed because $x_h(\boldsymbol \omega,s)$ is correlated to $\mathcal{A}_h(\boldsymbol \omega,s)$ for $s\in (t_{\ell-1},t_\ell)$. Because $x(s)$ is not correlated to $\mathcal{A}_h(\boldsymbol \omega,s)$, it was not necessary to introduce such a term in \eqref{eq:lem_este_step0}. The proof of Theorem \ref{thm:este} based on the error dynamics \eqref{eq:dyn_e1} presented above is thus simpler than a proof based on \eqref{eq:dyn_e2}.
\end{remark}

\begin{remark} \label{rem:dissipativity}
When we look back at the proof of Theorem \ref{thm:este}, we see that Assumption \ref{ass:dissipativity} is only used to assure that the matrices $A$ and $\mathcal{A}_h(\boldsymbol \omega,t)$ are dissipative (for all $\boldsymbol \omega$ with $p(\boldsymbol \omega) > 0$ and all $t \in [0,T]$). When Assumption \ref{ass:dissipativity} is not satisfied, there must exist a constant $a > 0$ such that $\hat{A} = A-aI$ and $\hat{\mathcal{A}}_h(\boldsymbol \omega,t) = \mathcal{A}_h(\boldsymbol \omega,t) - a I$ are dissipative (for all $\boldsymbol \omega$ with $p(\boldsymbol \omega) > 0$ and all $t \in [0,T]$). Because $\mathbb{E}[\mathcal{A}_h(t)] = A$, it follows that $\mathbb{E}[\hat{\mathcal{A}}_h(t)] = \mathbb{E}[\mathcal{A}_h(t)] - aI = A - aI = \hat{A}$ and $\mathrm{Var}[\| \hat{\mathcal{A}}_h(t) - \hat{A} \|^2] = \mathrm{Var}[\mathcal{A}]$. When we let $\hat{x}(t)$ and $\hat{x}_h(\boldsymbol \omega,t)$ denote the solutions generated by $\hat{A}$ and $\hat{\mathcal{A}}_h(\boldsymbol \omega,t)$, respectively, we can now prove in a similar way as in Theorem \ref{thm:este} that the error $\hat{e}_h(\boldsymbol \omega,t) = \hat{x}_h(\boldsymbol \omega,t) - \hat{x}(t)$ can be bounded as
\begin{equation}
\mathbb{E}[\lvert \hat{e}_h(t) \rvert^2] \leq h \mathrm{Var}[\mathcal{A}](\| \hat{A} \|t^2+2t)(\lvert x_0 \rvert + \lvert Bu \rvert_{L^1(0,T;\mathbb{R}^N)})^2. \label{eq:rem_dissipativity1}
\end{equation}
Because $x(t) = e^{at}\hat{x}(t)$ and $x_h(\boldsymbol \omega,t) = e^{at}\hat{x}_h(\boldsymbol \omega,t)$, also
\begin{equation}
e_h(\boldsymbol \omega,t) = x_h(\boldsymbol \omega,t) - x(t) = e^{at} \hat{x}_h(\boldsymbol \omega,t) - e^{at} \hat{x}(t) = e^{at} \hat{e}_h(\boldsymbol \omega,t). 
\end{equation}
Taking the expectation and using \eqref{eq:rem_dissipativity1}, we find
\begin{equation}
\mathbb{E}[\lvert e_h(t) \rvert^2] \leq h e^{at} \mathrm{Var}[\mathcal{A}](\| \hat{A} \|t^2+2t)(\lvert x_0 \rvert + \lvert Bu \rvert_{L^1(0,T;\mathbb{R}^N)})^2. 
\end{equation}
The error estimate now grows exponentially in time. 
\end{remark}

\subsection{The forward dynamics with a stochastic input}
\label{ssec:forward_dynamics_stochastic}
In this subsection, we prove a result similar to Theorem \ref{thm:este} for inputs $u_h(\boldsymbol \omega,t)$ that are stochastic, i.e., which depend on $\boldsymbol \omega$. We thus want to bound the error
\begin{equation}
e_h(\boldsymbol \omega,t) = x_h(\boldsymbol \omega,t) - x(\boldsymbol \omega,t),
\end{equation}
where $x_h(\boldsymbol \omega,t)$ and $x(\boldsymbol \omega,t)$ are the solutions of \eqref{eq:dyn_xtilde_om} and \eqref{eq:dyn_x_om}, respectively. 

To this end, we consider the semi-group $e^{At}$ generated by the matrix $A$ and the evolution operator $S_h(\boldsymbol \omega,t,s)$ associated to $\mathcal{A}_h(\boldsymbol \omega,t)$. The evolution operator $S_h(\boldsymbol \omega,t,s)$ is defined by property that for all vectors $x_s \in \mathbb{R}^N$ (and all $t \geq s$), $S_h(\boldsymbol \omega,t,s)x_s$ is equal to the solution $y_h(\boldsymbol \omega,t)$ of
\begin{equation}
\dot{y}_h(\boldsymbol \omega,t) = \mathcal{A}_h(\boldsymbol \omega,t) y_h(\omega,t), \qquad \qquad y_h(\boldsymbol \omega,s) = x_s. \label{eq:dyn_ytilde}
\end{equation}

\begin{figure}
\includegraphics[width=0.9\textwidth]{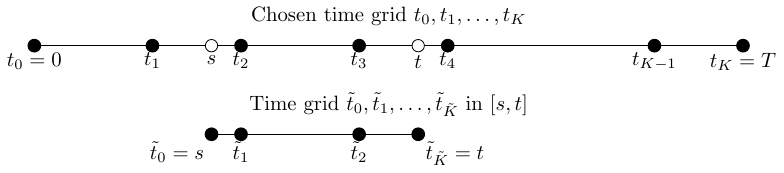}
\caption{The relation between the chosen time grid $t_0, t_1, \ldots, t_K$ and the time grid $\tilde{t}_0, \tilde{t}_1, \ldots , \tilde{t}_{\tilde{K}}$ used in Remark \ref{rem:Shnew}. In the displayed example, $\ell = 2$, $k = 4$, and $\tilde{K} = 3$. }
\label{fig:timegrids}
\end{figure}

\begin{remark} \label{rem:Shnew}
An explicit formula for the evolution operator $S_h(\boldsymbol \omega, t, s)$ can be obtained as follows. 
Let $0 \leq s \leq t \leq T$ and let $\ell, k \in \{1,2, \ldots, K \}$ be selected such that
\begin{equation}
s \in [t_{\ell-1}, t_\ell), \qquad \qquad t \in [t_{k-1}, t_k). 
\end{equation}
By restricting the given time grid $0 = t_0 < t_1 < t_2 < \ldots < t_{K-1} < t_K = T$ to the interval $[s,t]$, we obtain a grid with $\tilde{K} = k-\ell+1$ grid points
%\begin{equation}
%s < 
%t_\ell < 
%t_{\ell+1} < 
%\ldots < 
%t_{k-1}  < 
%t. 
%\end{equation}
\begin{equation}
\tilde{t}_0 := s < 
\tilde{t}_1 := t_\ell < 
\tilde{t}_2 := t_{\ell+1} < 
\ldots < 
\tilde{t}_{\tilde{K}-1} := t_{k-1}  < 
\tilde{t}_{\tilde{K}} := t. 
\end{equation}
The construction of the time grid $\tilde{t}_0, \tilde{t}_1, \ldots \tilde{t}_{\tilde{K}}$ is illustrated in Figure \ref{fig:timegrids}. We also denote $\tilde{h}_p := \tilde{t}_p - \tilde{t}_{p-1}$ (for $p \in \{1,2, \ldots, \tilde{K} \}$) and introduce (for each $\omega \in \{1,2, \ldots, 2^M \}$)
\begin{equation}
\mathcal{A}_\omega := \sum_{m \in S_\omega} \frac{A_m}{\pi_m}.
\end{equation}
Because $\mathcal{A}_h(\boldsymbol \omega,\tau) = \mathcal{A}_{\omega_p}$ is constant for $\tau \in [\tilde{t}_{p-1}, \tilde{t}_p)$, it is now easy to see that
\begin{equation}
S_h(\boldsymbol \omega, t,s) = e^{\mathcal{A}_{\omega_k} \tilde{h}_{\tilde{K}}} \cdots e^{\mathcal{A}_{\omega_{\ell+1}} \tilde{h}_2}  e^{\mathcal{A}_{\omega_\ell} \tilde{h}_1} =  \prod_{p=1}^{\tilde{K}} e^{\mathcal{A}_{\omega_p+\ell-1}\tilde{h}_p}.\label{eq:Udef}
\end{equation}
Under Assumption \ref{ass:dissipativity}, all matrices $\mathcal{A}_{\omega_p}$ are dissipative and \eqref{eq:Udef} shows that 
\begin{equation}
\| S_h(\boldsymbol \omega, t,s) \| \leq 1.
\end{equation}
\end{remark}

Using the variation of constants formula, the solutions of $x_h(\boldsymbol \omega,t)$ and $x(\boldsymbol \omega,t)$ can expressed as
\begin{align}
x_h(\boldsymbol \omega,t) &= S_h(\boldsymbol \omega,t,0) x_0 + \int_0^t S_h(\boldsymbol \omega,t,s) Bu_h(\boldsymbol \omega,s) \ \mathrm{d}s, \label{eq:xtildesol_om_varconst} \\
x(\boldsymbol \omega, t) &= e^{At} x_0 + \int_0^t e^{A(t-s)} Bu_h(\boldsymbol \omega,s) \ \mathrm{d}s. \label{eq:xsol_om_varconst}
\end{align}
Subtracting \eqref{eq:xsol_om_varconst} from \eqref{eq:xtildesol_om_varconst} we find the following expression for the error $e_h(\boldsymbol \omega,t)$
\begin{equation}
e_h(\boldsymbol \omega,t) = E_h(\boldsymbol \omega, t,0) x_0 + \int_0^t E_h(\boldsymbol \omega,t,s)Bu_h(\boldsymbol \omega,s) \ \mathrm{d}s, \label{eq:esol_om_varconst}
\end{equation}
where $E_h(\boldsymbol \omega,t,s) = S_h(\boldsymbol \omega,t,s) - e^{A(t-s)}$. The following corollary of Theorem \ref{thm:este} shows that we can bound $E_h(\boldsymbol \omega,t,s) = S_h(\boldsymbol \omega,t,s) - e^{A(t-s)}$. 

\begin{corollary} \label{corr:este1}
Under Assumptions \ref{ass:dissipativity} and \ref{ass:positivity}, we have that
\begin{equation}
\mathbb{E}[\| S_h(t,s) - e^{A(t-s)} \|^2] \leq (\|A \|T^2 + 2T) h  \mathrm{Var}[\mathcal{A}], \label{eq:cor_este1}
\end{equation}
for all $0 \leq s \leq t \leq T$. 
\end{corollary}

\begin{proof} Fix $s \in [0,T]$ and an initial condition $x_s \in \mathbb{R}^N$.

Define $y(t) = e^{A(t-s)}x_s$ and let $y_h(\boldsymbol \omega,t)$ be the solution of \eqref{eq:dyn_ytilde}, both for $t \in [s,T]$. We then apply Theorem \ref{thm:este} with $u(t) \equiv 0$ to the time-shifted solutions $\tilde{y}(\tilde{t}) = y(\tilde{t}+s)$ and $\tilde{y}_h(\boldsymbol \omega, \tilde{t}) = y_h(\boldsymbol \omega,\tilde{t}+s)$ and the time-shifted matrix $\tilde{\mathcal{A}}_h(\boldsymbol \omega,\tilde{t}) = \mathcal{A}_h(\boldsymbol \omega,\tilde{t}+s)$ defined on $\tilde{t} \in [0, T-s]$. We thus conclude that (writing $\tilde{t} = t-s$)
\begin{equation}
\mathbb{E}[\lvert y_h(t) - y(t) \rvert^2] = \mathbb{E}[\lvert \tilde{y}_h(\tilde{t}) - \tilde{y}(\tilde{t})\rvert^2] \leq h \mathrm{Var}[\mathcal{A}] (\| A\| \tilde{t}^2 + 2\tilde{t}) \lvert x_s \rvert^2.
\end{equation}
Noting that, by definition, $y(t) = e^{A(t-s)}x_s$ and $y_h(\boldsymbol \omega,t) = S_h(\boldsymbol \omega,t,s)x_s$, we find that (for $x_s \neq 0$)
\begin{equation}
\mathbb{E}\left[ \frac{\lvert (S_h(\boldsymbol \omega,t,s) - e^{A(t-s)})x_s \rvert^2}{\lvert x_s \rvert^2} \right] \leq h \mathrm{Var}[\mathcal{A}](\| A\| T^2 + 2T),
\end{equation}
where it was used that $\tilde{t} = t-s \leq T$. 
The result now follows from the definition of the operator-norm. 
\end{proof}

\begin{remark}
In Appendix \ref{app:commutative}, we prove a result similar to Corollary \ref{corr:este1} under the additional assumption that all matrices $A_m$ commute pairwise. The result in Appendix \ref{app:commutative} extends naturally to an infinite dimensional setting under the additional assumption that the domains of the operators $A_m$ are the same. This is not the case for Corollary \ref{corr:este1} because the operator norm $\| A \|$ appears in \eqref{eq:cor_este1}. 
\end{remark}

We are now ready for the main result of this subsection. 

\begin{theorem} \label{thm:este2} 
Consider any control $u_h : \Omega^K \rightarrow L^2(0,T;\mathbb{R}^q)$. 
Assume that Assumptions \ref{ass:dissipativity} and \ref{ass:positivity} are satisfied and let $U$ be such that \begin{equation}
\lvert Bu_h(\boldsymbol \omega) \rvert_{L^2(0,T;\mathbb{R}^q)} \leq U, \label{eq:boundU}
\end{equation} 
for all $\boldsymbol \omega \in \Omega^K$, 
then
\begin{equation}
\mathbb{E}[\lvert e_h(t) \rvert^2] \leq (\| A \|T^2 + 2T) h \mathrm{Var}[\mathcal{A}] \left( \lvert x_0 \rvert + U \sqrt{T}  \right)^2.
\end{equation}
\end{theorem}

\begin{proof} 
Using the triangle inequality in \eqref{eq:esol_om_varconst}, we find
\begin{align}
\lvert e_h(\boldsymbol \omega,t) \rvert &\leq \|E_h(\boldsymbol \omega,t,0) \| \lvert x_0 \rvert + \int_0^t \| E_h(\boldsymbol \omega,t,s) \| \lvert Bu_h(\boldsymbol \omega,s) \rvert \ \mathrm{d}s \nonumber \\
&\leq \|E_h(\boldsymbol \omega,t,0) \| \lvert x_0 \rvert + \sqrt{\int_0^t \| E_h(\boldsymbol \omega,t,s) \|^2 \ \mathrm{d}s} \lvert Bu_h(\boldsymbol \omega) \rvert_{L^2(0,T;\mathbb{R}^q)}, \label{eq:thm_este2_step0}
\end{align}
where the second inequality follows from the Cauchy-Schwarz inequality in $L^2(0,t)$. 
Squaring both sides and using the bound \eqref{eq:boundU}, we find
\begin{multline}
\lvert e_h(\boldsymbol \omega,t) \rvert^2 \leq \|E_h(\boldsymbol \omega,t,0) \|^2 \lvert x_0 \rvert^2 + U^2 \int_0^t \| E_h(\boldsymbol \omega,t,s) \|^2 \ \mathrm{d}s \\
+ 2 U \lvert x_0 \rvert \|E_h(\boldsymbol \omega,t,0) \| \sqrt{\int_0^t \| E_h(\boldsymbol \omega,t,s) \|^2 \ \mathrm{d}s}. \label{eq:thm_este2_step1}
\end{multline}
In order to use the bound from Corollary \ref{corr:este1} to estimate the last term, note that we can use the Cauchy-Schwartz inequality in the probability space to find
\begin{equation}
\mathbb{E}\left[\|E_h(t,0) \| \sqrt{\int_0^t \| E_h(t,s) \|^2 \ \mathrm{d}s} \right] \leq \sqrt{\mathbb{E}[\|E_h(t,0) \|^2 ] \int_0^t \mathbb{E}[ \| E_h(t,s) \|^2] \ \mathrm{d}s}
\end{equation}
Taking the expected value in \eqref{eq:thm_este2_step1} and using that the bound on $\mathbb{E}[\| E_h(t,s) \|^2]$ from Corollary \ref{corr:este1} does not depend on $t$ and $s$, we find
\begin{equation}
\mathbb{E}[\lvert e_h(t) \rvert^2] \leq (\lvert x_0 \rvert + U \sqrt{t})^2 (\| A \|T^2 + 2T) h \mathrm{Var}[\mathcal{A}], 
\end{equation}
which gives the desired estimate. 
\end{proof}

\begin{remark}
Because $\Omega^K$ is finite, we can always find a constant $U$ such that \eqref{eq:boundU} is satisfied for a given $u_h: \Omega^K \rightarrow L^2(0,T;\mathbb{R}^q)$. However, when we consider a family of temporal grids for which $h \rightarrow 0$, the constant $U$ may depend on $h$ (depending on the considered family of controls $u_h(\boldsymbol \omega,t)$). Fortunately, we only need to apply Theorem \ref{thm:este2} with $u_h(\boldsymbol \omega,t) = u_h^*(\boldsymbol \omega,t)$, where $u^*_h(\boldsymbol \omega,t)$ is the control that minimizes the cost functional $J_h(\boldsymbol \omega,\cdot)$ in \eqref{eq:Jtilde}. For this control, the coercivity of the cost functional $J_h(\boldsymbol \omega, \cdot)$ implies that the constant $U$ can be chosen independent of the considered temporal grid, see \eqref{eq:bound_uhopt}.  
\end{remark}

\begin{remark}
Note that the estimate in Theorem \ref{thm:este} depends on the $L^1$-norm of the control but that estimate in Theorem  \ref{thm:este2} depends through \eqref{eq:boundU} on the $L^2$-norm. Setting $u_h(\boldsymbol \omega,t) = u(t)$ in Theorem \ref{thm:este2} therefore does not give the estimate in Theorem \ref{thm:este}. This underlines the additional difficulty posed by stochastic controls. 
\end{remark}

\subsection{A no-gap condition} \label{ssec:nogap}
With the results regarding forward dynamics from the previous two subsections, we are now ready to address the optimal control problem. The main result of this subsection is the no-gap condition in Theorem \ref{thm:nogap}. To prove this result, we need the following technical lemma. 

\begin{lemma} \label{lem:dJuh} 
Consider any control $u_h : \Omega^K \rightarrow L^2(0,T;\mathbb{R}^q)$. Assume that Assumptions \ref{ass:dissipativity} and \ref{ass:positivity} hold and let $U > 0$ be such that \eqref{eq:boundU} is satisfied. Then
\begin{equation}
\mathbb{E}[\lvert J_h(u_h) - J(u_h) \rvert] \leq C_{[A,x_0, Q, x_d, T, U]} \left( \sqrt{h \mathrm{Var}[\mathcal{A}]} +  h \mathrm{Var}[\mathcal{A}] \right). \label{eq:lem_dJuh}
\end{equation}
\end{lemma}

\begin{proof}
Let $x(\boldsymbol \omega,t)$ and $x_h(\boldsymbol \omega,t)$ be the solutions of \eqref{eq:dyn_x_om} and \eqref{eq:dyn_xtilde_om} for the considered control $u_h(\boldsymbol \omega,t)$. For brevity, we write $\tilde{x}(\boldsymbol \omega,t) = x(\boldsymbol \omega,t) - x_d(t)$ and $\tilde{x}_h(\boldsymbol \omega,t) = x_h(\boldsymbol \omega,t) - x_d(t)$. By definition of the cost functionals $J(\cdot)$ and $J_h(\boldsymbol \omega,\cdot)$ in \eqref{eq:J} and \eqref{eq:Jtilde}, we have
\begin{align}
&J_h(\omega,u_h(\boldsymbol \omega)) - J(u_h(\boldsymbol \omega)) = \tfrac{1}{2}\int_0^T \left( \tilde{x}_h(\boldsymbol \omega,t)^\top Q \tilde{x}_h(\boldsymbol \omega,t) - \tilde{x}(\boldsymbol \omega,t)^\top Q \tilde{x}(\boldsymbol \omega,t) \right) \ \mathrm{d}t \nonumber \\
&= \int_0^T  \tilde{x}(\boldsymbol \omega,t)^\top Q (\tilde{x}_h(\boldsymbol \omega,t) - \tilde{x}(\boldsymbol \omega,t)) \ \mathrm{d}t \nonumber \\
&\qquad\qquad\qquad +\tfrac{1}{2} \int_0^T(\tilde{x}_h(\boldsymbol \omega,t) - \tilde{x}(\boldsymbol \omega,t))^\top Q (\tilde{x}_h(\boldsymbol \omega,t) - \tilde{x}(\boldsymbol \omega,t))  \ \mathrm{d}t \nonumber \\
&= \int_0^T \left( \tilde{x}(\boldsymbol \omega,t)^\top Q e_h(\boldsymbol \omega,t) +\tfrac{1}{2} e_h(\boldsymbol \omega,t)^\top Q e_h(\boldsymbol \omega,t) \right) \ \mathrm{d}t,
\end{align}
where the last identity follows because $e_h(\boldsymbol \omega,t) = x_h(\boldsymbol \omega,t) - x(t)=  \tilde{x}_h(\boldsymbol \omega,t) - \tilde{x}(t)$. 
Taking the absolute value and estimating the RHS, we find
\begin{align}
& \lvert J_h(\boldsymbol \omega,u_h)  - J(u_h(\boldsymbol \omega)) \rvert \leq \| Q \| \int_0^T \left( \lvert \tilde{x}(\boldsymbol \omega,t) \rvert \lvert e_h(\boldsymbol \omega,t)\rvert  + \tfrac{1}{2} \lvert e_h(\boldsymbol \omega,t)\rvert ^2 \right) \ \mathrm{d}t \nonumber \\
&\quad \leq \| Q \| \left( \lvert \tilde{x}(\boldsymbol \omega) \rvert_{L^2(0,T;\mathbb{R}^N)} \lvert e_h(\boldsymbol \omega) \rvert_{L^2(0,T;\mathbb{R}^N)} + \tfrac{1}{2}\lvert e_h(\boldsymbol \omega) \rvert_{L^2(0,T;\mathbb{R}^N)}^2 \right).
\end{align}
Taking the expectation and using the Cauchy-Schwartz inequality, we find that
\begin{multline}
\mathbb{E}[\lvert J_h(u_h) - J(u_h) \rvert ] \leq \\
\| Q \| \left( \sqrt{\mathbb{E}[\lvert \tilde{x} \rvert_{L^2(0,T;\mathbb{R}^N)}^2]} \sqrt{\mathbb{E}[\lvert e_h \rvert_{L^2(0,T;\mathbb{R}^N)}^2]} + \tfrac{1}{2}\mathbb{E}[\lvert e_h \rvert_{L^2(0,T;\mathbb{R}^N)}^2] \right). \label{eq:lem_dJuh_step1}
\end{multline}
Using the estimate from Theorem \ref{thm:este2}, we find
\begin{equation}
\mathbb{E}[\lvert e_h \rvert_{L^2(0,T;\mathbb{R}^N)}^2] =  \int_0^T \mathbb{E}[ \lvert e_h(t) \rvert^2] \ \mathrm{d}t \leq h \mathrm{Var}[\mathcal{A}] C_{[A,x_0,T,U]}. \label{eq:lem_dJuh_step2}
\end{equation}
Because $\tilde{x}(\boldsymbol \omega,t) = x(\boldsymbol \omega, t) - x_d(t)$, \eqref{eq:boundx} shows that
\begin{equation}
\lvert \tilde{x}(\boldsymbol \omega) \rvert^2_{L^2(0,T;\mathbb{R}^N)} \leq ( \sqrt{T}(\lvert x_0 \rvert + \lvert Bu_h(\boldsymbol \omega) \rvert_{L^1(0,T;\mathbb{R}^N)})  + \lvert x_d \rvert_{L^2(0,T;\mathbb{R}^N)} )^2.\label{eq:lem_dJuh_step3}
\end{equation}
Because $\lvert Bu_h(\boldsymbol \omega) \rvert_{L^1(0,T;\mathbb{R}^N)} \leq \sqrt{T} \lvert Bu_h(\boldsymbol \omega) \rvert_{L^2(0,T;\mathbb{R}^N)} \leq \sqrt{T}U$, we see from \eqref{eq:lem_dJuh_step3} that $\mathbb{E}[\lvert \tilde{x} \rvert_{L^2(0,T;\mathbb{R}^N)}^2] \leq  C_{[x_0, x_d, T, U]}$. The result now follows by inserting this estimate and \eqref{eq:lem_dJuh_step2} into \eqref{eq:lem_dJuh_step1}. 
\end{proof}

We are now ready to prove the main result of this section which can be considered as a no-gap condition for the RBM optimal control problem. 

\begin{theorem} \label{thm:nogap}
Let $u^*(t)$ be the (deterministic) control that minimizes the cost functional $J(u)$ in \eqref{eq:J} and let $u_h^*(\boldsymbol \omega,t)$ be the control that minimizes the cost functional $J_h(\boldsymbol \omega, u)$ in \eqref{eq:Jtilde}. Then
\begin{equation}
\mathbb{E}[\lvert J_h(u^*_h) - J(u^*)\rvert] \leq C_{[A,B,x_0,Q,R,x_d,T]} \left(  \sqrt{h \mathrm{Var}[\mathcal{A}]} + h \mathrm{Var}[\mathcal{A}] \right). \label{eq:thm_nogap}
\end{equation}
\end{theorem}

\begin{proof}
We have that
\begin{multline}
J(u^*) \leq J(u_h^*(\boldsymbol \omega)) = J_h(\boldsymbol \omega, u_h^*(\boldsymbol \omega)) + \delta(\boldsymbol \omega) \\
 \leq J_h(\boldsymbol \omega,u^*) + \delta(\boldsymbol \omega) = J(u^*) + \delta(\boldsymbol \omega) + \varepsilon(\boldsymbol \omega), \label{eq:thm_nogap_step1}
\end{multline}
where $\delta(\boldsymbol \omega) = J(u^*_h(\boldsymbol \omega)) - J_h(\boldsymbol \omega, u^*_h(\boldsymbol \omega))$ and $\varepsilon(\boldsymbol \omega) = J_h(\boldsymbol \omega,u^*) - J(u^*)$. Note that the first inequality follows because $u^*$ is the minimizer of $J$ and the second inequality because $u_h^*(\boldsymbol \omega)$ is the minimizer of $J_h(\boldsymbol \omega, \cdot)$. Subtracting $J(u^*) + \delta(\boldsymbol \omega)$ from the first, third, and fifth expressions in \eqref{eq:thm_nogap_step1}, shows that
\begin{equation}
-\delta(\boldsymbol \omega) \leq 
J_h(\boldsymbol \omega,u^*_h(\boldsymbol \omega)) - J(u^*) \leq \varepsilon(\boldsymbol \omega). 
\end{equation}
Taking the absolute value, we find
\begin{equation}
\lvert J_h(\boldsymbol \omega,u^*_h(\boldsymbol \omega)) - J(u^*) \rvert \leq \max\{\lvert \delta(\boldsymbol \omega) \rvert, \lvert \varepsilon(\boldsymbol \omega) \rvert \} \leq \lvert \delta(\boldsymbol \omega) \rvert + \lvert \varepsilon(\boldsymbol \omega) \rvert.
\end{equation}
Therefore also
\begin{equation}
\mathbb{E}[\lvert J_h(u^*_h) - J(u^*) \rvert] \leq \mathbb{E}[\lvert \delta \rvert] + \mathbb{E}[\lvert \varepsilon \rvert]. \label{eq:thm_nogap_step2}
\end{equation}

Lemma \ref{lem:dJuh} can now be used to find bounds for $\mathbb{E}[\lvert \delta \rvert] = \mathbb{E}[\lvert J_h(u_h^*) - J(u_h^*) \rvert]$ and $\mathbb{E}[\lvert \varepsilon \rvert] = \mathbb{E}[\lvert J_h(u^*)-J(u^*) \rvert]$. 

For the bound on $\mathbb{E}[\lvert \delta \rvert]$, we use that \eqref{eq:bound_uhopt} shows that there exists a constant such that $\lvert Bu^*_h(\boldsymbol \omega)\rvert_{L^2(0,T;\mathbb{R}^N)} \leq C_{[B,x_0,Q,R,x_d,T]}$ so that \eqref{eq:boundU} is satisfied with a constant $U$ that does not depend on the used temporal grid $t_0, t_1, \ldots, t_K$. Lemma \ref{lem:dJuh} thus implies that
\begin{equation}
\mathbb{E}[\lvert \delta \rvert ] \leq C_{[A,B,x_0,Q,R,x_d,T]} \left(  \sqrt{h \mathrm{Var}[\mathcal{A}]} + h \mathrm{Var}[\mathcal{A}] \right). \label{eq:thm_nogap_step6}
\end{equation}

For the bound on $\mathbb{E}[\lvert \varepsilon \rvert]$, we can simply take $U = \lvert Bu^*(t) \rvert_{L^2(0,T;\mathbb{R}^N)}$, which is a constant that only depends on the parameters $A,B,x_0,Q,R,x_d,T$ that define the deterministic problem \eqref{eq:dyn_x}--\eqref{eq:J}. Lemma \ref{lem:dJuh} thus also shows that
\begin{equation}
\mathbb{E}[\lvert \varepsilon \rvert] \leq C_{[A,B,x_0,Q,R,x_d,T]} \left(  \sqrt{h \mathrm{Var}[\mathcal{A}]} + h \mathrm{Var}[\mathcal{A}] \right). \label{eq:thm_nogap_step3}
\end{equation}
Inserting \eqref{eq:thm_nogap_step6} and \eqref{eq:thm_nogap_step3} into \eqref{eq:thm_nogap_step2} we find \eqref{eq:thm_nogap}. 
\end{proof}

\subsection{Convergence in the controls} \label{ssec:controls}
In the last stage of our analysis of the RBM-optimal control problem, we bound the expected difference between the optimal control $u^*_h$ that minimizes $J_h$ in \eqref{eq:Jtilde} and the optimal control $u^*$ for the original problem. The proof is based on the strong convexity of the functional $J_h$ in \eqref{eq:Jtilde}. 

To prove the main result, we need the following lemma which bounds the difference between the G\^ateaux derivative of $J_h$ and the G\^ateaux derivative of $J$ in expectation. 

\begin{lemma} \label{lem:gateaux}
For any deterministic control $u \in L^2(0,T;\mathbb{R}^q)$ and any stochastic perturbation $v_h : \Omega^K \rightarrow L^2(0,T;\mathbb{R}^q)$,
\begin{equation}
\mathbb{E} [\lvert \delta J_h(u; v_h) - \delta J(u;v_h) \rvert] \leq C_{[A,B,x_0,Q,x_d,T,u]} \sqrt{h \mathrm{Var}[\mathcal{A}]} \sqrt{\mathbb{E}[\lvert v_h \rvert_{L^2(0,T;\mathbb{R}^q)}^2]}. \label{eq:lem_gateaux}
\end{equation}
\end{lemma}

\begin{proof}
Let $x(t)$ and $x_h(\boldsymbol \omega,t)$ be the solutions of \eqref{eq:dyn_x} and \eqref{eq:dyn_xtilde}, respectively.  
Furthermore, denote
\begin{equation}
y(\boldsymbol \omega, t) = \int_0^t e^{A(t-s)}Bv_h(\boldsymbol \omega,s) \ \mathrm{d}s, \qquad
y_h(\boldsymbol \omega, t) = \int_0^t S_h(\boldsymbol \omega,t,s)Bv_h(\boldsymbol \omega,s) \ \mathrm{d}s. \label{eq:lem_gateaux_step0}
\end{equation}

Directly from the definition of the G\^ateaux derivative, we find that
\begin{align}
\delta J(u,v_h(\boldsymbol \omega)) &= \int_0^T \left( \tilde{x}(t)^\top Q y(\boldsymbol \omega,t) + u(t)^\top R v_h(\boldsymbol \omega,t) \right) \ \mathrm{d}t, \label{eq:lem_gateaux_step1} \\
\delta J_h(\boldsymbol \omega, u,v_h(\boldsymbol \omega)) &= \int_0^T \left( \tilde{x}_h(\boldsymbol \omega,t)^\top Q y_h(\boldsymbol \omega,t) + u(t)^\top R v_h(\boldsymbol \omega,t) \right) \ \mathrm{d}t, \label{eq:lem_gateaux_step2}
\end{align}
where we write $\tilde{x}(t) = x(t) - x_d(t)$ and $\tilde{x}_h(\boldsymbol \omega,t) = x_h(\boldsymbol \omega,t) - x_d(t)$. 

Subtracting \eqref{eq:lem_gateaux_step1} from \eqref{eq:lem_gateaux_step2}, we find
\begin{align}
\delta J_h(&\boldsymbol \omega,u,v_h(\boldsymbol \omega)) - \delta J(u,v_h(\boldsymbol \omega)) \nonumber \\
&= \int_0^T \left( \tilde{x}_h(\boldsymbol \omega,t)^\top Q y_h(\boldsymbol \omega,t) - \tilde{x}(t)^\top Q y(\boldsymbol \omega,t) \right) \ \mathrm{d}t \nonumber \\
&= \int_0^T \left( \tilde{x}_h(\boldsymbol \omega,t)^\top Q (y_h(\boldsymbol \omega,t) - y(\boldsymbol \omega,t)) + (\tilde{x}_h(\boldsymbol \omega,t) - \tilde{x}(t))^\top Q y(\boldsymbol \omega,t) \right) \ \mathrm{d}t \nonumber \\
&= \int_0^T \left( \tilde{x}_h(\boldsymbol \omega,t)^\top Q f_h(\boldsymbol \omega,t) + e_h(\boldsymbol \omega,t)^\top Q y(\boldsymbol \omega,t) \right) \ \mathrm{d}t, \label{eq:lem_gateaux_step3}
\end{align}
where $e_h(\boldsymbol \omega,t) = x_h(\boldsymbol \omega,t) - x(t) = \tilde{x}_h(\boldsymbol \omega,t) - \tilde{x}(t)$ and $f_h(\boldsymbol \omega,t) = y_h(\boldsymbol \omega,t) - y(\boldsymbol \omega,t)$.
Taking the absolute value, we find
\begin{align}
\lvert \delta J_h(&\boldsymbol \omega,u,v_h(\boldsymbol \omega)) - \delta J(u,v_h(\boldsymbol \omega)) \rvert \nonumber \\
&\leq \| Q \|\int_0^T \left( \lvert \tilde{x}_h(\boldsymbol \omega,t)\rvert \lvert f_h(\boldsymbol \omega,t) \rvert + \lvert e_h(\boldsymbol \omega,t)\rvert \lvert y(\boldsymbol \omega,t) \rvert \right) \ \mathrm{d}t.  \label{eq:lem_gateaux_step4}
\end{align}

Using \eqref{eq:boundxh}, we find the following bound for $\tilde{x}_h(\boldsymbol \omega,t) = x_h(\boldsymbol \omega, t) - x_d(t)$
\begin{equation}
\lvert \tilde{x}_h(\boldsymbol \omega,t) \rvert \leq \lvert x_h(\boldsymbol \omega,t)\rvert + \lvert x_d(t)\rvert \leq \lvert x_0 \rvert + \lvert B u \rvert_{L^1(0,T;\mathbb{R}^N)} + \lvert x_d(t) \rvert.
\end{equation}
We thus have $\lvert \tilde{x}_h(\boldsymbol \omega,t) \rvert \leq  C_{[B,x_0,x_d,T,u]}$ for all $\boldsymbol \omega \in \Omega^K$. 

Taking the expectation in \eqref{eq:lem_gateaux_step4} using this result shows that
\begin{align}
\mathbb{E}[&\lvert\delta J_h(u,v_h) - \delta J(u,v_h) \rvert ] \nonumber \\
&\leq \| Q \|\int_0^T \left( C_{[B,x_0,x_d,T,u]} \mathbb{E}[\lvert f_h(t) \rvert] - \sqrt{\mathbb{E}[\lvert e_h(t) \rvert^2]} \sqrt{\mathbb{E}[ \lvert y(t) \rvert^2]} \right) \ \mathrm{d}t, \label{eq:lem_gateaux_step5}
\end{align}
where the second term on the RHS follows from the Cauchy-Schwartz inequality. 

Again using the notation $E_h(\boldsymbol \omega,t,s) := S_h(\boldsymbol \omega,t,s) - e^{A(t-s)}$, \eqref{eq:lem_gateaux_step0} shows that 
\begin{equation}
f_h(\boldsymbol \omega,t) = y_h(\boldsymbol \omega,t) - y(\boldsymbol \omega,t) = \int_0^t E_h(\boldsymbol \omega,t,s) B v_h(\boldsymbol \omega,s) \ \mathrm{d}s.
\end{equation}
Therefore, 
\begin{align}
\mathbb{E}[\lvert f_h(t) \rvert] &\leq \int_0^t \mathbb{E}[\|E_h(t,s) \| \lvert Bv_h(s) \rvert] \ \mathrm{d}s \nonumber \\
&\leq \int_0^t \sqrt{\mathbb{E}[\|E_h(t,s) \|^2]} \sqrt{\mathbb{E}[ \lvert Bv_h(s) \rvert^2]} \ \mathrm{d}s \nonumber \\
&\leq C_{[A,T]}\sqrt{h \mathrm{Var}[\mathcal{A}]} \int_0^t \sqrt{\mathbb{E}[\lvert Bv_h(s)\rvert^2]} \ \mathrm{d}s \nonumber \\
&\leq C_{[A,T]} \sqrt{h \mathrm{Var}[\mathcal{A}]} \sqrt{t} \sqrt{\int_0^t \mathbb{E}[\lvert Bv_h(s) \rvert^2] \ \mathrm{d}s} \nonumber \\
&\leq C_{[A,T]} \sqrt{h \mathrm{Var}[\mathcal{A}]}  \sqrt{\mathbb{E}[\lvert Bv_h\rvert^2_{L^2(0,T;\mathbb{R}^N)}]}, \label{eq:lem_gateaux_step6}
\end{align}
where the second inequality follows from the Cauchy-Schwartz inequality in the probability space, the third inequality from Corollary \ref{corr:este1}, and the third inequality from the Cauchy-Schwartz inequality in $L^2(0,t)$. 

Because the control $u(t)$ is deterministic, Theorem \ref{thm:este} shows that
\begin{equation}
\mathbb{E}[\lvert e_h(t) \rvert^2] \leq h \mathrm{Var}[\mathcal{A}] C_{[A,B,x_0,T,u]}. \label{eq:lem_gateaux_step7}
\end{equation}
Finally, note 
\begin{multline}
\lvert y(\boldsymbol \omega,t) \rvert^2 = \left( \int_0^t \| e^{A(t-s)}\| \lvert Bv_h(\boldsymbol \omega,s) \rvert \ \mathrm{d}s \right)^2 \\ 
\leq \int_0^t \| e^{A(t-s)}\|^2 \ \mathrm{d}s \int_0^t \lvert Bv_h(\boldsymbol \omega,s)r\vert^2 \ \mathrm{d}s \leq t\lvert Bv_h(\boldsymbol \omega)\rvert_{L^2(0,T;\mathbb{R}^N)}^2. \label{eq:lem_gateaux_step9}
\end{multline}
Therefore, also
\begin{equation}
\mathbb{E}[\lvert y(t) \rvert^2] \leq C_{[B,T]}\mathbb{E}[\lvert v_h \rvert_{L^2(0,T;\mathbb{R}^N)}^2] . \label{eq:lem_gateaux_step8}
\end{equation}
Inserting \eqref{eq:lem_gateaux_step6}, \eqref{eq:lem_gateaux_step7}, and \eqref{eq:lem_gateaux_step8} into \eqref{eq:lem_gateaux_step5} completes the proof. 
\end{proof}

We are now ready to prove the convergence result for the optimal controls. 

\begin{theorem} \label{thm:controls}
Suppose that the functional $J_h(\boldsymbol \omega,\cdot)$ in \eqref{eq:Jtilde} is $\alpha$-convex for all $\boldsymbol \omega \in \Omega^K$. Let $u_h^*(\boldsymbol \omega,t)$ be the minimizer of $J_h(\boldsymbol \omega,\cdot)$ in \eqref{eq:Jtilde} and $u^*(t)$ be the minimizer of $J$ in \eqref{eq:J}, then
\begin{equation}
\alpha^2 \mathbb{E}[\lvert u_h^* - u^*\rvert_{L^2(0,T; \mathbb{R}^q)}^2] \leq C_{[A,B,x_0,Q,R,x_d,T]} h \mathrm{Var}[\mathcal{A}]. \label{eq:thm_controls}
\end{equation}
\end{theorem} 

\begin{proof} 
We apply \eqref{eq:gateaux} with $J(\cdot) = J_h(\boldsymbol \omega, \cdot)$, $v = u^*_h(\boldsymbol \omega)$, and $u = u^*$ to find 
\begin{equation}
J_h(\omega,u^*_h(\boldsymbol \omega)) \geq J_h(\boldsymbol \omega,u^*) + \delta J_h(\boldsymbol \omega, u^*; u_h^*(\boldsymbol \omega) - u^*) + \tfrac{\alpha}{2} \lvert u_h^*(\boldsymbol \omega) - u^* \rvert^2_{L^2(0,T;\mathbb{R}^q)}.
\end{equation}
Because $u^*_h(\boldsymbol \omega)$ is the minimizer of $J_h(\boldsymbol \omega, \cdot)$, $J_h(\boldsymbol \omega,u_h^*(\boldsymbol \omega)) \leq J_h(\boldsymbol \omega,u^*)$ and
\begin{equation}
0 \geq \delta J_h(\boldsymbol \omega,u^*; u_h^*(\boldsymbol \omega) - u^*) + \tfrac{\alpha}{2} \lvert u_h^*(\boldsymbol \omega) - u^* \rvert^2_{L^2(0,T;\mathbb{R}^q)}.
\end{equation}
Bringing $\delta J_h$ to the other side, taking the absolute value and then the expectation, yields
\begin{equation}
\frac{\alpha}{2}\mathbb{E}[\lvert u_h^*  - u^*\rvert^2_{L^2(0,T;\mathbb{R}^q)}] \leq \mathbb{E}[\lvert \delta J_h(u^*; u_h^* - u^*) \rvert]. 
\end{equation}
Since $u^*$ is the minimizer of $J$,  $\delta J(u^*, v) = 0$ for all perturbation $v \in L^2(0,T;\mathbb{R}^q)$. In particular, we have that $\delta J(u^*, u_h^*(\boldsymbol \omega) - u^*) = 0$ for all $\boldsymbol \omega \in \Omega^K$ so that also
\begin{equation}
\frac{\alpha}{2}\mathbb{E}[\lvert u_h^* - u^* \rvert^2_{L^2(0,T;\mathbb{R}^q)}] \leq \mathbb{E}[\lvert \delta J_h(u^*; u_h^* - u^*) - \delta J(u^*; u_h^* - u^* )\rvert ]. 
\end{equation}
We now apply Lemma \ref{lem:gateaux} to the RHS with $u = u^*$ and $v_h(\boldsymbol \omega) = u^*_h(\boldsymbol \omega) - u^*$, which shows that
\begin{equation}
\frac{\alpha}{2}\mathbb{E}[\lvert u_h^*  - u^* \rvert^2_{L^2(0,T;\mathbb{R}^q)}] \leq C_{[B,x_0,Q,x_d,T,u^*]} \sqrt{h \mathrm{Var}[\mathcal{A}]} \sqrt{\mathbb{E}[\lvert u_h^* - u^* \rvert^2_{L^2(0,T;\mathbb{R}^q)}]}. \label{eq:thm_convu_step1}
\end{equation}
Next, we divide \eqref{eq:thm_convu_step1} by $\tfrac{1}{2} \sqrt{\mathbb{E}[\lvert u_h^*-u^* \rvert^2_{L^2(0,T;\mathbb{R}^q)}]}$ to find
\begin{equation}
\alpha \sqrt{\mathbb{E}[\lvert u_h^* - u^* \rvert^2_{L^2(0,T;\mathbb{R}^q)}]} \leq C_{[A,B,x_0,Q,x_d,T,u^*]} \sqrt{h \mathrm{Var}[\mathcal{A}]}. \label{eq:thm_convu_step2}
\end{equation}
Squaring both sides we arrive at
\begin{equation}
\alpha^2 \mathbb{E}[\lvert u_h^* - u^* \rvert^2_{L^2(0,T;\mathbb{R}^q)}] \leq C_{[A,B,x_0,Q,x_d,T,u^*]} h \mathrm{Var}[\mathcal{A}]. 
\end{equation}
The result follows because the optimal control $u^*(t)$ only depends on the parameters $A,B,x_0,Q,R,x_d$, and $T$ that define the original problem \eqref{eq:dyn_x}--\eqref{eq:J}. 
\end{proof}

We now point out two corollaries of Theorem \ref{thm:controls} that are important when we use the control $u_h^*(\boldsymbol \omega ,t)$ (optimized for the RBM-dynamics) to control the original dynamics. For the first corollary, we introduce the notation 
\begin{align}
x^*_h(\boldsymbol \omega, t) &= e^{At}x_0 + \int_0^t e^{A(t-s)}Bu_h^*(\boldsymbol \omega, s) \ \mathrm{d}s, \label{eq:def_xstarh} \\
x^*(t) &= e^{At}x_0 + \int_0^t e^{A(t-s)}Bu^*(s) \ \mathrm{d}s, \label{eq:def_xstar}
\end{align}
i.e., $x^*_h(\boldsymbol \omega, t)$ is the solution of the original dynamics \eqref{eq:dyn_x} resulting from the control $u_h^*(\boldsymbol \omega, t)$ optimized for the RBM-dynamics and $x^*(t)$ is the solution of the original dynamics \eqref{eq:dyn_x} resulting from the optimal control $u^*(t)$. 
\begin{corollary} \label{corr:xstar}
Suppose that the functional $J_h(\boldsymbol \omega,\cdot)$ in \eqref{eq:Jtilde} is $\alpha$-convex for all $\boldsymbol \omega \in \Omega^K$ and let $x_h^*(\boldsymbol \omega, t)$ and $x^*(t)$ be as in \eqref{eq:def_xstarh} and \eqref{eq:def_xstar}, respectively. Then
\begin{equation}
\alpha^2 \mathbb{E}[\lvert x^*_h(t) - x^*(t) \rvert^2] \leq C_{[A,B,x_0,Q,R,x_d,T]} h \mathrm{Var}[\mathcal{A}],
\end{equation}
for all $t \in [0,T]$. 
\end{corollary}
\begin{proof}
Note that
\begin{equation}
x^*_h(\boldsymbol \omega, t) - x^*(t) = \int_0^t e^{A(t-s)}B(u_h^*(\boldsymbol \omega, s)- u^*(s)) \ \mathrm{d}s.
\end{equation}
Therefore also
\begin{align}
&\lvert x^*_h(\boldsymbol \omega, t) - x^*(t) \rvert \leq \int_0^t \| e^{A(t-s)}\| \| B\| \lvert u_h^*(\boldsymbol \omega, s)- u^*(s) \rvert \ \mathrm{d}s \nonumber \\ &\leq \| B \| \lvert u_h^*(\boldsymbol \omega)- u^* \rvert_{L^1(0,T; \mathbb{R}^q)} \leq \| B \| \sqrt{T} \sqrt{\lvert u_h^*(\boldsymbol \omega)- u^*\rvert_{L^2(0,T; \mathbb{R}^q)}},
\end{align}
where the second inequality uses that $\|e^{At}\| \leq 1$ in view of Assumption \ref{ass:dissipativity}. The result now follows after squaring this inequality, taking the expectation, and using \eqref{eq:thm_controls}. 
\end{proof}

\begin{corollary} \label{corr:suboptimality}
Suppose that the cost functional $J_h(\boldsymbol \omega,\cdot)$ is $\alpha$-convex for all $\boldsymbol \omega \in \Omega^K$. Let $u^*(t)$ be the (deterministic) control that minimizes the cost functional $J(u)$ in \eqref{eq:J} and let $u_h^*(\boldsymbol \omega,t)$ be the control that minimizes the cost functional $J_h(\boldsymbol \omega, u)$ in \eqref{eq:Jtilde}. Then
\begin{equation}
\alpha^2 \mathbb{E}[\lvert J(u^*_h) - J(u^*) \rvert] \leq C_{[A,B,x_0,Q,R,x_d,T]} h \mathrm{Var}[\mathcal{A}]. \label{eq:corr_suboptimality}
\end{equation}
\end{corollary}

\begin{proof}
Denote $v_h(\boldsymbol \omega,t) := u_h^*(\boldsymbol \omega,t) - u^*(t)$ and $y(\boldsymbol \omega,t) := \int_0^t e^{A(t-s)} Bv_h(\boldsymbol \omega,s) \ \mathrm{d}s$. Because the considered functional is quadratic,
\begin{align}
J(u^*_h(\boldsymbol \omega)) - J(u^*) &= J(u^* + v_h(\boldsymbol \omega)) - J(u^*) \nonumber \\
&= \delta J(u^*, v_h(\boldsymbol \omega)) + \delta^2 J(v_h(\boldsymbol \omega), v_h(\boldsymbol \omega)), \label{eq:corr_suboptimality_step1}
\end{align}
where the Hessian $\delta^2 J(v_h(\boldsymbol \omega), v_h(\boldsymbol \omega))$ is given by
\begin{align}
\delta^2 J(v_h(\boldsymbol \omega), v_h(\boldsymbol \omega)) &= \frac{1}{2}\int_0^T \left( y(\boldsymbol \omega,t)^\top Q y(\boldsymbol \omega,t) + v_h(\boldsymbol \omega,t)^\top R v_h(\boldsymbol \omega,t) \right) \ \mathrm{d}t. 
\end{align}
Because $u^*$ is the minimizer of $J(\cdot)$, $\delta J(u^*, v) = 0$ for all $v \in L^2(0,T; \mathbb{R}^q)$. The first term on the RHS of \eqref{eq:corr_suboptimality_step1} thus vanishes. Also observe that
\begin{equation}
\delta^2 J(v_h(\boldsymbol \omega), v_h(\boldsymbol \omega)) \leq \tfrac{1}{2}\| Q \| \lvert y(\boldsymbol \omega) \rvert_{L^2(0,T;\mathbb{R}^N)}^2 + \tfrac{1}{2}\| R \| \lvert v_h(\boldsymbol \omega) \rvert_{L^2(0,T;\mathbb{R}^q)}^2.
\end{equation}
A similar estimate as \eqref{eq:lem_gateaux_step9} shows that $\lvert y(\boldsymbol \omega)\rvert_{L^2(0,T;\mathbb{R}^N)}^2 \leq C_{[B,T]}\lvert v_h(\boldsymbol \omega)\rvert_{L^2(0,T;\mathbb{R}^q)}^2$. Combining these results in \eqref{eq:corr_suboptimality_step1}, we conclude
\begin{align}
\lvert J(u^*_h(\boldsymbol \omega)) - J(u^*) \rvert &\leq J(u^*_h(\boldsymbol \omega)) - J(u^*) \nonumber \\
&\leq \delta^2 J(v_h(\boldsymbol \omega), v_h(\boldsymbol \omega)) \leq C_{[B,Q,R,T]}\lvert v_h(\boldsymbol \omega) \rvert_{L^2(0,T;\mathbb{R}^q)}^2. 
\end{align}
The result now follows after taking the expectation and using the result from Theorem \ref{thm:controls} to bound $\mathbb{E}[\lvert v_h \rvert_{L^2(0,T;\mathbb{R}^q})^2] = \mathbb{E}[\lvert u^*_h - u^* \rvert_{L^2(0,T;\mathbb{R}^q)}^2]$. 
\end{proof}

\section{Numerical results} \label{sec:examples}
In this section, we apply our proposed method to three medium to large scale linear dynamical systems that are obtained after spatial discretization of a linear PDE. 

\subsection{A discretized 1D heat equation} \label{ssec:example1}
We consider a controlled heat equation on the 1-D spatial domain $[-L, L]$, 
\begin{align}
&y_t(t,\xi) = y_{\xi\xi}(t,\xi) + \chi_{[-L/3,0]}(\xi) u(t), \qquad \xi \in [-L,L], \label{eq:heat1D_PDE} \\
&y_\xi(t,-L) = y_\xi(t,L) = 0, \qquad \qquad y(0,\xi) = e^{-\xi^2} + \xi^2e^{-L^2},  \label{eq:heat1D_IC}
\end{align}
where $\chi_{[-L/3,0]}(\xi)$ denotes the characteristic function for the interval $[-L/3,0]$. We want to compute the optimal control $u^*(t)$ that minimizes
\begin{equation}
\mathcal{J}(u) = \frac{100}{2}\int_0^T \int_{-L}^0 y(t,\xi)^2 \ \mathrm{d}\xi \ \mathrm{d}t + \frac{1}{2}\int_0^T u(t)^2 \ \mathrm{d}t. \label{eq:heat1D_cost}
\end{equation}
The spatial discretization of the dynamics \eqref{eq:heat1D_PDE}--\eqref{eq:heat1D_IC} is made by finite differences and the cost functional in \eqref{eq:heat1D_cost} is discretized by the trapezoid rule. We choose a uniform spatial grid with $N = 61$ grid points $\xi_i = (i-1)\Delta \xi - L$ ($i \in \{ 1,2, \ldots , N \}$), where $\Delta \xi = 2L/(N-1)$ is the grid spacing, and obtain a system of the form \eqref{eq:dyn_x}. 

The resulting $A$-matrix is of the form
\begin{equation}
A = \frac{1}{\Delta \xi^2}\begin{bmatrix}
-2 &  2 & 0 & \cdots & 0 & 0 & 0 \\
 1 & -2 & 1 &        & 0 & 0 & 0 \\
 0 & 1 & -2 &        & 0 & 0 & 0 \\
 \vdots &  &  & \ddots & & & \vdots \\
 0 & 0 & 0 &  & -2 & 1 & 0 \\
 0 & 0 & 0 &  & 1 & -2 & 1 \\
 0 & 0 & 0 & \cdots & 0 & 2 & -2
\end{bmatrix}. \label{eq:heat1D_A}
\end{equation}
Observe that $A$ can be written as 
\begin{equation}
A = \sum_{i=1}^{n} \tilde{A}_i, \label{eq:heat1D_sumAi}
\end{equation}
where the $n := N-1 = 60$ matrices $\tilde{A}_i \in \mathbb{R}^{N \times N}$ are zero except for the entries
\begin{align*}
\begin{bmatrix}
[\tilde{A}_1]_{11} & [\tilde{A}_1]_{12} \\
[\tilde{A}_1]_{21} & [\tilde{A}_1]_{22}
\end{bmatrix}%tilde{A}_1]_{1:2,1:2} 
&= \begin{bmatrix}
-2 & 2 \\
1 & -1
\end{bmatrix}, \\
\begin{bmatrix}
[\tilde{A}_i]_{ii} & [\tilde{A}_i]_{i,i+1} \\
[\tilde{A}_i]_{i+1,i} & [\tilde{A}_i]_{i+1,i+1}
\end{bmatrix}%[\tilde{A}_i]_{i:i+1,i:i+1} 
&= \begin{bmatrix}
-1 & 1 \\
1 & -1
\end{bmatrix}, & 2 \leq i \leq n-1, \\
\begin{bmatrix}
[\tilde{A}_n]_{nn} & [\tilde{A}_n]_{n,n+1} \\
[\tilde{A}_n]_{n+1,n} & [\tilde{A}_n]_{n+1,n+1}
\end{bmatrix} %[\tilde{A}_n]_{n:n+1,n:n+1} 
&= \begin{bmatrix}
-1 & 1 \\
2 & -2
\end{bmatrix}.
\end{align*}
%\begin{align}
%&[\tilde{A}_1]_{11} = -2, \qquad [\tilde{A}_1]_{22} = -1, \qquad &&[\tilde{A}_1]_{12} = 2, \qquad [\tilde{A}_1]_{21} = 1, \\
%&[\tilde{A}_i]_{ii} = [\tilde{A}_i]_{i+1,i+1} = -1, \qquad &&[\tilde{A}_i]_{i,i+1} = [\tilde{A}_i]_{i+1,i} = 1, \\
%&[\tilde{A}_n]_{nn} = -1, \quad [\tilde{A}_n]_{n+1,n+1} = -2, \quad &&[\tilde{A}_n]_{n,n+1} = 1, \quad [\tilde{A}_n]_{n+1,n} = 2.
%\end{align}
One can easily verify that the matrices $\tilde{A}_i$ are dissipative. We now define the $M$ submatrices $A_m$ (for $M = 1,2,3,4$) as
\begin{equation}
A_m = \sum_{i=i_{m-1}+1}^{i_m} \tilde{A}_i, \label{eq:heat1D_defAm}
\end{equation}
where $i_m = nm/M$. Because of \eqref{eq:heat1D_sumAi}, it is easy to see that the submatrices $A_m$ satisfy \eqref{eq:Aform}. Because the submatrices $\tilde{A}_i$ are dissipative, the submatrices $A_m$ in \eqref{eq:heat1D_defAm} are dissipative and Assumption \ref{ass:dissipativity} is satisfied. 
\begin{example} \label{ex:decomposition}
For $M = 2$ and $N = 61$, we obtain the splitting of the $A$-matrix in \eqref{eq:heat1D_A} as $A = A_1 + A_2$, with
\begin{equation}
A_1 = \begin{bmatrix}
A_{11} & 0_{31\times 30} \\
0_{30\times 31} & 0_{30\times 30}
\end{bmatrix}, \qquad \qquad A_2 = \begin{bmatrix}
0_{30\times 30} & 0_{30\times 31} \\
0_{31\times 30} & A_{22}
\end{bmatrix},
\end{equation}
where $A_{11}$ and $A_{22}$ are the $31 \times 31$-matrices 
\begin{align}
A_{11} &= \frac{1}{\Delta \xi^2}\begin{bmatrix}
-2 &  2 & 0 & \cdots & 0 & 0 & 0 \\
 1 & -2 & 1 &        & 0 & 0 & 0 \\
 0 & 1 & -2 &        & 0 & 0 & 0 \\
 \vdots &  &  & \ddots & & & \vdots \\
 0 & 0 & 0 &  & -2 & 1 & 0 \\
 0 & 0 & 0 &  & 1 & -2 & 1 \\
 0 & 0 & 0 & \cdots & 0 & 1 & -1
\end{bmatrix}, \\
A_{22} &= \frac{1}{\Delta \xi^2}\begin{bmatrix}
-1 &  1 & 0 & \cdots & 0 & 0 & 0 \\
 1 & -2 & 1 &        & 0 & 0 & 0 \\
 0 & 1 & -2 &        & 0 & 0 & 0 \\
 \vdots &  &  & \ddots & & & \vdots \\
 0 & 0 & 0 &  & -2 & 1 & 0 \\
 0 & 0 & 0 &  & 1 & -2 & 1 \\
 0 & 0 & 0 & \cdots & 0 & 2 & -2
\end{bmatrix}. 
\end{align}
\end{example}

%\begin{remark}
%We can interpret the matrices $A_{1}$ and $A_{2}$ in Example \ref{ex:decomposition} as a finite difference discretizations of the operators
%\begin{equation}
%\mathfrak{A}_1 = 
%\end{equation}
%\end{remark}

We will present numerical results for four cases:
\begin{description}
\item[Case i] We decompose $A$ into $M = 2$ submatrices and assign a probability $\tfrac{1}{2}$ to the subsets $\{ 1 \}$ and $ \{ 2 \}$  and a probability $0$ to the subsets $\emptyset$ and $\{1,2 \}$. 
\item[Case ii] We decompose $A$ into $M = 3$ submatrices and assign a probability $\tfrac{1}{3}$ to the subsets $\{ 1 \}$, $ \{ 2 \}$, and $\{ 3 \}$ and a probability $0$ to the other subsets of $\{1,2,3 \}$. 
\item[Case iii] We decompose $A$ into $M = 4$ submatrices and assign a probability $\tfrac{1}{4}$ to the subsets $\{ 1 \}$, $ \{ 2 \}$, $\{ 3 \}$, and $\{ 4 \}$ and a probability $0$ to the other subsets of $\{1,2,3,4\}$. 
\item[Case iv] We decompose $A$ into $M = 4$ submatrices and assign a probability $\tfrac{1}{2}$ to the subsets $\{ 1,3 \}$ and $ \{ 2,4 \}$ and a probability $0$ to the other subsets of $\{1,2,3,4\}$. %We thus effectively decompose the domain $[-L, L]$ into two disconnected subdomains $[-L, -L/2] \cup [0,L/2]$ and $[-L/2,0] \cup [L/2,L]$. 
\end{description}
In all 4 cases, we fix $N = 61$, $L = \tfrac{3}{2}$, and $T = \tfrac{1}{2}$. 

We use a uniform grid $0=t_0 < t_1 < \ldots < t_{K-1} < t_K = T$ with a uniform grid spacing $h$. We will present results for $h = 2^{-5}$, $2^{-7}$, $2^{-9}$, $2^{-11}$, $2^{-13}$, and $2^{-15}$. For each of the $K = T/h$ time intervals $[t_{k-1},t_k)$, we select an index $\omega_k$ according to the probabilities specified in Cases i--iv above. The state $x_h(\boldsymbol \omega,t)$ that satisfies \eqref{eq:dyn_xtilde} is computed using a single Crank-Nicholson step in each time interval $[t_{k-1}, t_k)$. We use precomputed LU-factorizations of the matrices $I - \tfrac{h}{2}\sum_{m\in S_{\omega}}\tfrac{A_m}{\pi_m}$ (for subsets $S_{\omega}$ with a nonzero probability $p_{\omega}$) that need to be inverted frequently. 

The optimal control $u_h^*(\boldsymbol \omega,t)$ that minimizes $J_h(\boldsymbol \omega,u)$ in \eqref{eq:Jtilde} is computed with a gradient-descent algorithm. The gradient is computed using the adjoint state $\varphi_h(\boldsymbol \omega,t)$, see Remark \ref{rem:phih}. The time discretization for the adjoint state equation \eqref{eq:dyn_phih} is done using the scheme proposed in \cite{apel2012} that leads to discretely consistent gradients. 
The iterates $u^k$ are computed as $u^{k+1} = u^k-\beta\nabla J_h(\boldsymbol \omega,u^k)$. The step size $\beta$ is chosen such that $J_h(\boldsymbol \omega, u^k - \beta \nabla J_h(\boldsymbol \omega,u^k))$ is minimal. The algorithm is terminated when the relative change in $J_h(\boldsymbol \omega,u)$ is below $10^{-6}$. 

\begin{figure}[h!]
\begin{subfigure}{0.45\textwidth}
\includegraphics[scale=0.78]{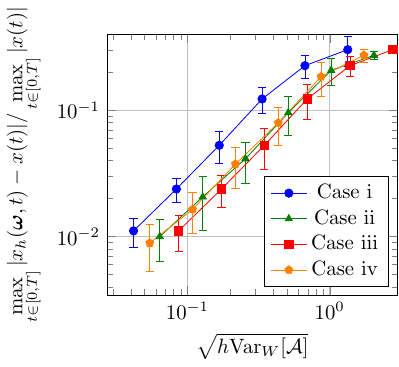}
\caption{Error in $x_h(t)$}
\label{fig:heat1D_errX}
\end{subfigure} ~
\begin{subfigure}{0.45\textwidth}
\includegraphics[scale=0.78]{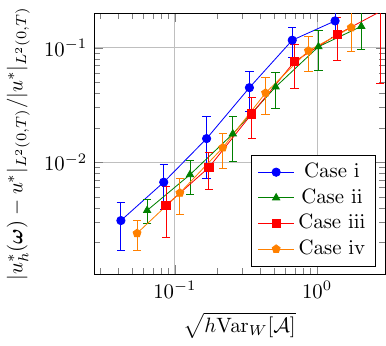}
\caption{Error in $u_h^*(t)$}
\label{fig:heat1D_errU}
\end{subfigure}\\
\begin{subfigure}{0.45\textwidth}
\includegraphics[scale=0.78]{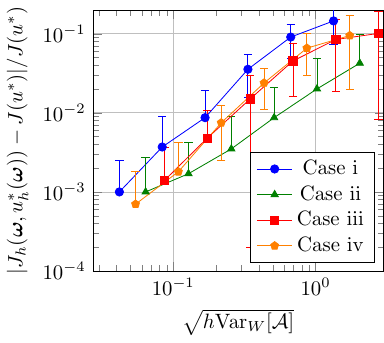}
\caption{Error in $J_h(u^*_h)$}
\label{fig:heat1D_errJ}
\end{subfigure} ~
\begin{subfigure}{0.45\textwidth}
\includegraphics[scale=0.78]{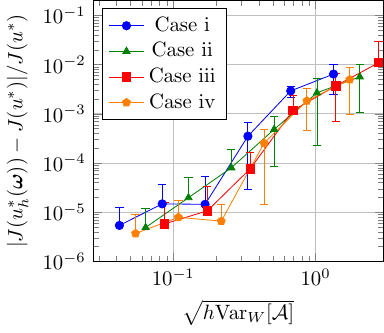}
\caption{Error in $J(u^*_h)$}
\label{fig:heat1D_errJ2}
\end{subfigure}\\
\begin{subfigure}{0.45\textwidth}
\includegraphics[scale=0.78]{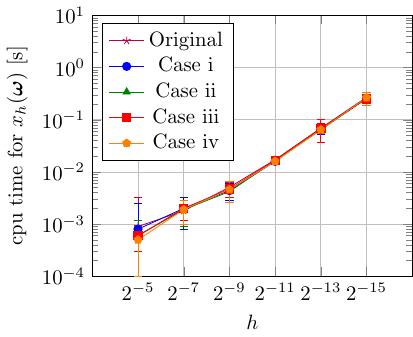}
\caption{Computational time for $x_h(t)$}
\label{fig:heat1D_durX}
\end{subfigure} ~
\begin{subfigure}{0.45\textwidth}
\includegraphics[scale=0.78]{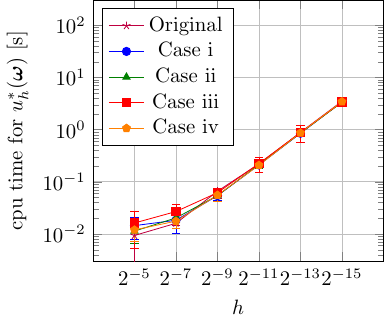}
\caption{Computational time for $u_h^*(t)$}
\label{fig:heat1D_durU}
\end{subfigure}
\caption{Simulation results for the discretized 1D heat equation}
\label{fig:heat1D}
\end{figure}

The results for the four considered cases are displayed in Figure \ref{fig:heat1D}. Because the obtained results depend on the randomly selected indices stored in $\boldsymbol \omega$, each marker in the subfigures in Figure \ref{fig:heat1D} represents the average error or duration over 25 random realizations of $\boldsymbol \omega$. The errorbars represent the $2\sigma$-confidence interval estimated from these 25 realizations. The errors are computed w.r.t.\ the solutions $x(t)$ and $u^*(t)$ that are computed on the same time grid as the corresponding solutions $x_h(\boldsymbol \omega,t)$ and $u^*_h(\boldsymbol \omega,t)$. The displayed errors therefore do not reflect the errors due to the temporal (or spatial) discretization but capture only the error introduced by the proposed randomized splitting method. 

\begin{figure*}
\begin{subfigure}{0.45\textwidth}
\includegraphics[scale=0.78]{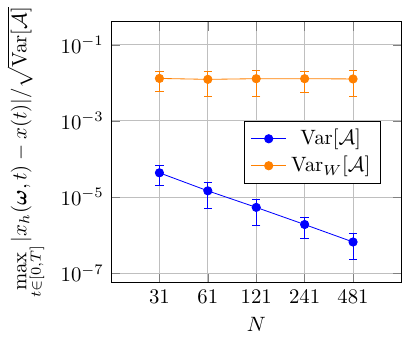}
\caption{Error in $x_h(t)$}
\label{fig:heat1D_convX}
\end{subfigure} ~
\begin{subfigure}{0.45\textwidth}
\includegraphics[scale=0.78]{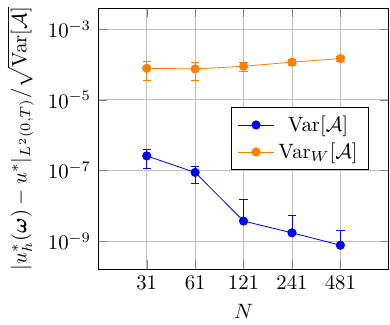}
\caption{Error in $u_h^*(t)$}
\label{fig:heat1D_convU}
\end{subfigure}
\caption{The errors in the forward dynamics $x_h(\omega,t)$ and the optimal control $u_h^*(\omega,t)$ divided by $\mathrm{Var}[\mathcal{A}]$ and $\mathrm{Var}_W[\mathcal{A}]$ (with $W = (A - 0.1 I)^{-1}$) as a function of the number of nodes $N$. The results are presented for case i, so $A$ is decomposed in $M = 2$ parts. }
\label{fig:heat1Dconv}
\end{figure*}

Because the matrices $A$ and $A_m$ represent approximations of unbounded operators, the variance $\mathrm{Var}[\mathcal{A}]$ defined in \eqref{eq:def_varA} will grow unbounded when the mesh is refined. This is also reflected by the large values of $\mathrm{Var}[\mathcal{A}]$ given in Table \ref{tab:heat1D_varA}. 
It is therefore more natural to consider the variance $\mathrm{Var}_W[\mathcal{A}]$ in \eqref{eq:def_varAW} weighted by a matrix of the form $W = (A - \lambda I)^{-1}$. The values of $\mathrm{Var}_W[\mathcal{A}]$ are indeed much smaller than the values of $\mathrm{Var}[\mathcal{A}]$ in Table \ref{tab:heat1D_varA}. The results at the end of this subsection (in Figure \ref{fig:heat1Dconv}) also indicate that the weighted variance $\mathrm{Var}_W[\mathcal{A}]$ reflects the behavior of the error better when the mesh is refined. 

The error estimates in Theorems \ref{thm:este}, \ref{thm:nogap}, and \ref{thm:controls} and in Corollary \ref{corr:suboptimality} are proportional to $h\mathrm{Var}[\mathcal{A}]$. We therefore plot the errors in Figures \ref{fig:heat1D_errX}--\ref{fig:heat1D_errJ2} against $\sqrt{h\mathrm{Var}_W[\mathcal{A}]}$ (with $W = (A - 0.1 I)^{-1}$) and expect that the errors for the different cases will be (approximately) on one line.

\begin{table}
\caption{Values of $\mathrm{Var}[\mathcal{A}]$ and $\mathrm{Var}_W[\mathcal{A}]$ for $W = (A - \lambda I)^{-1}$ with $\lambda = 0.1$}
\label{tab:heat1D_varA}
\centering
\begin{tabular}{|c|c|c|c|c|} \hline
   & Case i & Case ii & Case iii & Case iv \\ \hline
$\mathrm{Var}[\mathcal{A}]$ 
   & $4.16 \cdot 10^7$ & $1.65 \cdot 10^8$ & $3.68 \cdot 10^8$ & $4.16 \cdot 10^7$ \\ \hline
$\mathrm{Var}_W[\mathcal{A}]$ 
   & $57.32$ & $133.91$ & $246.54$ & $96.68$ \\ \hline
\end{tabular}
\end{table}

Figure \ref{fig:heat1D_errX} shows the difference $\lvert x_h(\boldsymbol \omega,t) - x(t)\rvert$ between the solutions $x(t)$ and $x_h(\omega,t)$ of \eqref{eq:dyn_x} and \eqref{eq:dyn_xtilde} with $u(t) = 0$. Recall that the markers in this figure indicate the average error observed over 25 realizations of $\boldsymbol \omega$, and are thus estimates for $\mathbb{E}[\max_{t \in [0,T]}\lvert x_h(t) - x(t) \rvert]$. Because $\mathbb{E}[\lvert x_h(t) - x(t) \rvert] \leq \sqrt{\mathbb{E}[ \lvert x_h(t) - x(t)\rvert^2]}$, we expect (based on the bound in Theorem \ref{thm:este}) that the errors in Figure \ref{fig:heat1D_errX} are proportional to $\sqrt{h\mathrm{Var}_W[\mathcal{A}]}$. This is indeed confirmed by Figure \ref{fig:heat1D_errX}.

\begin{figure}
\begin{subfigure}{0.45\textwidth}
\includegraphics[scale=0.65]{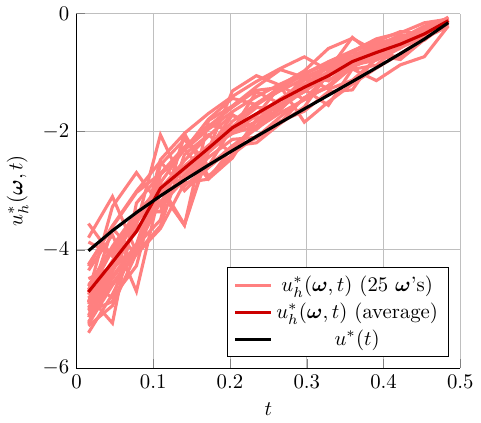}
\caption{$h = 2^{-5}$}
\label{fig:heat1D_controls1}
\end{subfigure} ~
\begin{subfigure}{0.45\textwidth}
\includegraphics[scale=0.65]{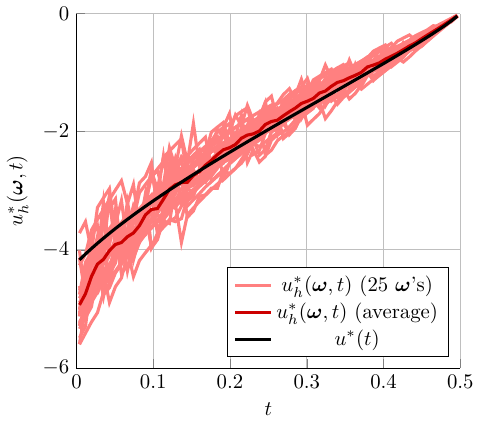}
\caption{$h = 2^{-7}$}
\label{fig:heat1D_controls2}
\end{subfigure}\\
\begin{subfigure}{0.45\textwidth}
\includegraphics[scale=0.65]{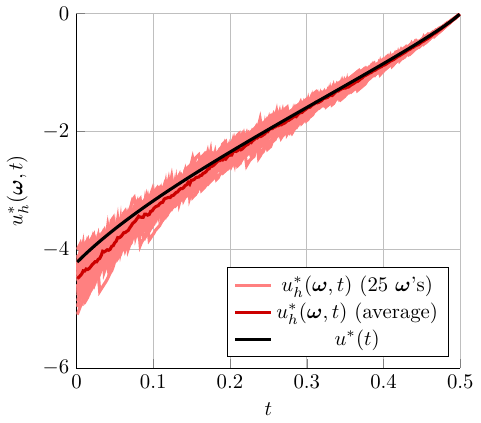}
\caption{$h = 2^{-9}$}
\label{fig:heat1D_controls3}
\end{subfigure} ~
\begin{subfigure}{0.45\textwidth}
\includegraphics[scale=0.65]{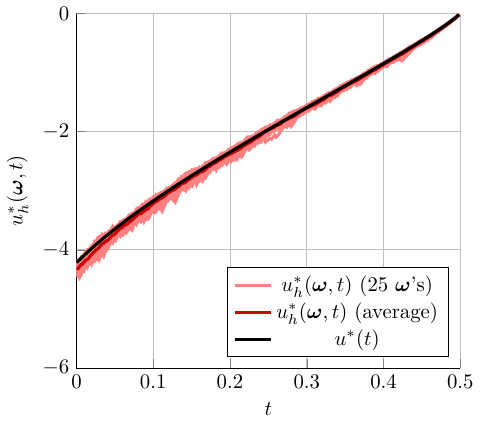}
\caption{$h = 2^{-11}$}
\label{fig:heat1D_controls4}
\end{subfigure}\\
\begin{subfigure}{0.45\textwidth}
\includegraphics[scale=0.65]{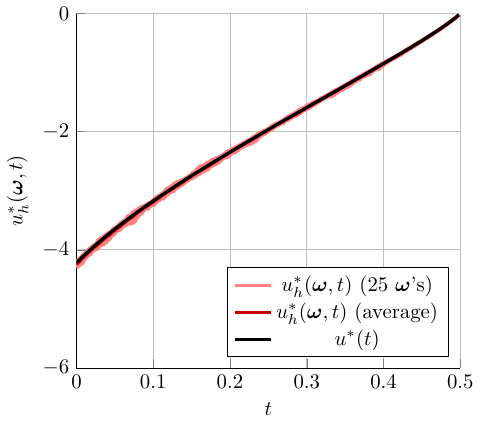}
\caption{$h = 2^{-13}$}
\label{fig:heat1D_controls5}
\end{subfigure} ~
\begin{subfigure}{0.45\textwidth}
\includegraphics[scale=0.65]{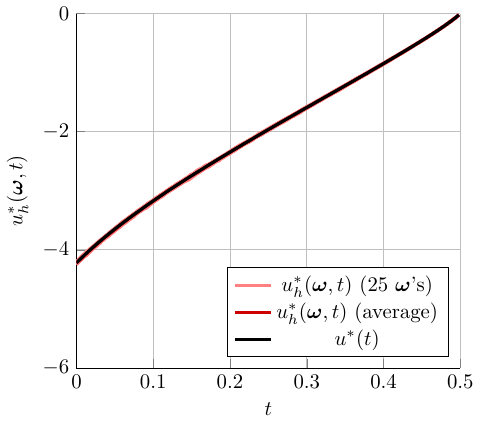}
\caption{$h = 2^{-15}$}
\label{fig:heat1D_controls6}
\end{subfigure}
\caption{The optimal controls computed for the 1D heat equation for different time steps $h$. The controls $u^*_h(\boldsymbol \omega,t)$ computed with the proposed randomized time-splitting method are shown for 25 realizations of $\boldsymbol \omega$ and compared to the optimal control $u^*(t)$ for the original system. }
\label{fig:heat1D_controls}
\end{figure}

Figure \ref{fig:heat1D_errU} shows the difference $\lvert u_h^* - u^* \rvert_{L^2(0,T)}$ between the optimal controls $u^*(t)$ and $u_h^*(\boldsymbol \omega,t)$ that minimize \eqref{eq:J} and \eqref{eq:Jtilde}, respectively. Based on the estimate in Theorem \ref{thm:controls}, we again expect that the observed errors are proportional to $\sqrt{h\mathrm{Var}_W[\mathcal{A}]}$. This is indeed the case  and the proportionality constants for the different cases are again (approximately) equal, which is also expected based on the error estimate in Theorem \ref{thm:controls}. 

The convergence in the optimal controls in Figure \ref{fig:heat1D_errU} is also illustrated in Figure \ref{fig:heat1D_controls}. This figure shows the optimal controls $u_h^*(\boldsymbol \omega,t)$ obtained for 25 randomly selected realizations of $\boldsymbol \omega \in \Omega^K$ (light red) for the six considered grid spacings $h$ of the temporal grid. The figure also shows the average of the 25 optimal controls $u^*_h(\boldsymbol \omega,t)$ (dark red) and the optimal control $u^*(t)$ for the original system (black). Figure \ref{fig:heat1D_controls} indeed shows that the optimal controls $u_h^*\boldsymbol \omega,t)$ get closer to the optimal control $u^*(t)$ when the spacing of the temporal grid $h$ is reduced. Especially in Figures \ref{fig:heat1D_controls1} and \ref{fig:heat1D_controls2}, it is also clear that the average of the 25 optimal controls $u^*_h(\boldsymbol \omega,t)$ (dark red) is not equal to the optimal control $u^*(t)$ for the original system (black). This indicates that $\mathbb{E}[u_h^*] \neq u^*$, see also Remark \ref{rem:bias_x}. This means that $u^*_h$ is a biased estimator for $u^*$ and averaging several realizations of $u^*(\boldsymbol \omega,t)$ can only improve the approximation of $u^*(t)$ to a limited extend. Note, however, that
\begin{equation}
\lvert \mathbb{E}[u_h^*] - u^* \rvert = \lvert \mathbb{E}[u_h^* - u^*] \rvert \leq \mathbb{E}[\lvert u^*_h - u^* \rvert] \leq \sqrt{\mathbb{E}[\lvert u^*_h - u^* \rvert^2]},
\end{equation}
so that Theorem \ref{thm:controls} shows that $\mathbb{E}[u_h^*] \rightarrow u^*$ at a rate of $\sqrt{h \mathrm{Var}[\mathcal{A}]}$. An analysis of the numerical results (that is not presented in Figure \ref{fig:heat1D}) also indicates that the average of the 25 realizations of $u_h^*(\boldsymbol \omega,t)$ converges to $u^*(t)$ at this rate.

Figures \ref{fig:heat1D_errJ} and \ref{fig:heat1D_errJ2} illustrate the convergence of $J_h(\boldsymbol \omega,u^*_h(\boldsymbol \omega))$ and $J(u^*_h(\boldsymbol \omega))$ to $J(u^*)$. Figure \ref{fig:heat1D_errJ} illustrates the error estimate in Theorem \ref{thm:nogap} and shows that the optimality gap $\lvert J_h(\omega,u^*_h(\boldsymbol \omega)) - J(u^*) \rvert$ is indeed proportional to $\sqrt{h\mathrm{Var}_W[\mathcal{A}]}$. The difference between the different cases is more visible than in Figures \ref{fig:heat1D_errX} and \ref{fig:heat1D_errU}. Figure \ref{fig:heat1D_errJ2} illustrates the error estimate in Corollary \ref{corr:suboptimality}, which shows that the suboptimality of the RBM-control $\lvert J(u_h^*(\boldsymbol \omega)) - J(u^*) \rvert$ is proportional to $h\mathrm{Var}_W[\mathcal{A}]$. The convergence rate is now twice as high as in the previous cases and the relative error stabilizes around $10^{-5}$, which seems to be related to the tolerance of $10^{-6}$ used in the computation of the optimal controls. 

Figures \ref{fig:heat1D_durX} and \ref{fig:heat1D_durU} show the computational times for (one realization of) $x_h(\boldsymbol \omega,t)$ and $u^*_h(\boldsymbol \omega,t)$ in Cases i--iv and the computational time for the original problem (labeled `Original'). Note that the results have been generated on temporal grids with different grid spacings $h$ and that the computational time generally increases when the more time steps are used, i.e.\ when $h$ is smaller. The figures indicate that $x_h(\boldsymbol \omega,t)$ and $u^*_h(\boldsymbol \omega,t)$ are not computed faster than the solutions $x(t)$ and $u^*(t)$ of the original problem. The proposed method does thus not lead to any reduction in computational time in this example. It seems that we cannot observe any reduction in computational time for this example because the original $A$-matrix is quite small ($N = 61$) and sparse ($A$ is tridiagonal). The examples in the following two subsections indicate that a reduction in computational cost is obtained when the state dimension $N$ is significantly higher or when $A$ has significantly more nonzero off-diagonal elements. 

To conclude this example, we study the dependence of our results on the number of grid points $N$. This gives us some indication whether the RBM can also be applied to infinite dimensional problems. In particular, the results give us some indication whether the proposed randomized splitting also works for the underlying PDE problem \eqref{eq:heat1D_PDE}--\eqref{eq:heat1D_cost}. As we also noted in Remarks \ref{rem:def_varAW} and \ref{rem:infdim}, the main concerns are related to operator norm of $A$, that appears in $\mathrm{Var}[\mathcal{A}]$ and in the estimate in Theorem \ref{thm:este}, which grows unbounded when the mesh is refined. These concerns also motivated the introduction of the weighted variance $\mathrm{Var}_W[\mathcal{A}]$, see Remark \ref{rem:def_varAW}. 

When the estimate in Theorem \ref{thm:este} indeed depends on $\| A \|$, the error $\lvert x_h(\boldsymbol \omega,t) - x(t) \rvert$ divided by $\mathrm{Var}[\mathcal{A}]$ should grow when $N$ is increased. Figure \ref{fig:heat1D_convX} shows that this is not the case, but that this ratio actually decreases when $N$ is increased. However, when we divided the errors by $\mathrm{Var}_W[\mathcal{A}]$, the result seems to be independent of the mesh size. Figure \ref{fig:heat1D_convU} shows that the same trend is observed for the errors in the optimal control. 

The numerical results in Figure \ref{fig:heat1Dconv} match well with the result from Appendix \ref{app:commutative}, where we prove an error estimate proportional to $\mathrm{Var}_W[\mathcal{A}]$ under the additional assumption that all matrices $A_m$ commute. This result also extends to an infinite-dimensional setting when the domains the operators $A_m$ coincide. However, in the setting considered here, the matrices $A_m$ do not commute and are not approximations of operators with the same domains. Proving the convergence of the proposed randomized time splitting method for the underlying PDE problem \eqref{eq:heat1D_PDE}--\eqref{eq:heat1D_cost} with the proposed randomized time splitting method is a challenging topic for future research.

\subsection{A discretized 3D heat equation} \label{ssec:example2}
We now consider a heat equation on the a 3-D spatial domain $V = [-L, L]^3$, 
\begin{align}
&y_t(t,\boldsymbol \xi) = \Delta y(t, \boldsymbol \xi), \qquad\qquad\qquad & \boldsymbol \xi \in [-L,L]^3, \label{eq:heat3D_PDE} \\
&\nabla y(t,\boldsymbol \xi) \cdot \mathbf{n} = u(t), \qquad \qquad\qquad & \boldsymbol \xi \in S_{\mathrm{top}}, \\
&\nabla y(t,\boldsymbol \xi) \cdot \mathbf{n} = 0, \qquad \qquad\qquad & \boldsymbol \xi \in \partial V \backslash S_{\mathrm{top}}, \\
&y(0,\boldsymbol \xi) = e^{-\lvert \boldsymbol \xi \rvert^2/(8L^2)},  \label{eq:heat3D_IC}
\end{align}
where $\nabla$ and $\Delta$ are the gradient and Laplacian operators  w.r.t.\ $\boldsymbol \xi$, $\mathbf{n}$ is the outward pointing normal, and $S_{\mathrm{top}}$ denotes the top surface $S_{\mathrm{top}} = \{ (\xi_1, \xi_2, \xi_3) \in [-L,L]^3 \mid \xi_3 = L \}$. The control $u(t)$ can be considered as a uniform heat load on the top surface. We want to compute the control $u^*(t)$ that minimizes
\begin{equation}
J = 1000 \int_0^T \iint_{S_{\mathrm{side}}} (y(t,\boldsymbol \xi))^2 \ \mathrm{d}\boldsymbol \xi \ \mathrm{d}t + \int_0^T (u(t))^2 \ \mathrm{d}t,\label{eq:heat3D_cost}
\end{equation}
where $S_{\mathrm{side}}= \{ (\xi_1, \xi_2, \xi_3) \in [-L,L]^3 \mid \xi_1 = -L \}$. We fix $L = 0.75$ and $T = 2$. 

The spatial discretization of \eqref{eq:heat3D_PDE}--\eqref{eq:heat3D_cost} is made by finite differences using $16 \times 16 \times 16$ grid points the $\xi_1$-, $\xi_2$-, and $\xi_3$-directions. This leads to a model of the form \eqref{eq:dyn_x}--\eqref{eq:J} with $N = 16^3 = 4096$ states. The resulting $A$-matrix is again dissipative. We create the decomposition of $A$ into submatrices $A_m$ by observing that $A$ is diagonally dominant. In particular, we have that
\begin{equation}
[A]_{ii} = -\sum_{\substack{j=1\\j \neq i}}^N [A]_{ij},
\end{equation}
where the off-diagonal elements $[A]_{ij}$ ($j \neq i$) are positive and the diagonal elements $[A]_{ii}$ are negative. 
By associating a matrix $\tilde{A}_{ij} \in \mathbb{R}^{N\times N}$ to each pair $(i,j)$ with $j > i$, we obtain a decomposition of $A$ as
\begin{equation}
A = \sum_{\substack{j=1\\j > i}}^N \tilde{A}_{ij},
\end{equation}
where the matrices $\tilde{A}_{ij}$ ($j > i$) are zero except for the entries
\begin{equation}
\begin{bmatrix}
[\tilde{A}_{ij}]_{ii} & [\tilde{A}_{ij}]_{ij} \\
[\tilde{A}_{ij}]_{ji} & [\tilde{A}_{ij}]_{jj}
\end{bmatrix} = [A]_{ij} \begin{bmatrix}
-1 & 1\\ 
1 & -1
\end{bmatrix}
\end{equation}
Because the off-diagonal elements $[A]_{ij} \geq 0$ ($j \neq i$), it is easy to verify that all the matrices $\tilde{A}_{ij}$ are dissipative. Also note that the matrix $A$ contains many zero off-diagonal elements, so that many of the matrices $\tilde{A}_{ij}$ are zero. There are only $3(16-1)16^2 = 11,520$ nonzero off-diagonal elements and thus only $11,520$ nonzero matrices $\tilde{A}_{ij}$. The $11,520$ nonzero matrices $\tilde{A}_{ij}$ are randomly divided into $M$ groups of (approximately) equal size. The matrices $A_m$ in \eqref{eq:Aform} are formed by summing the matrices $\tilde{A}_{ij}$ in each group. 

\begin{figure}[h!]
\begin{subfigure}{0.45\textwidth}
\includegraphics[scale=0.78]{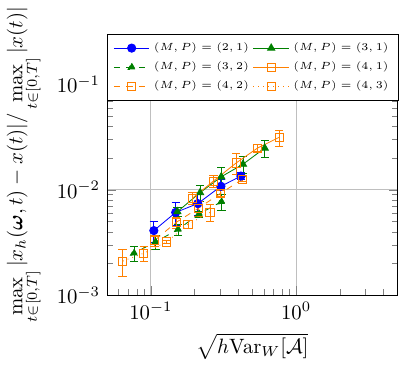}
\caption{Error in $x_h(t)$}
\label{fig:heat3D_errX}
\end{subfigure} ~
\begin{subfigure}{0.45\textwidth}
\includegraphics[scale=0.78]{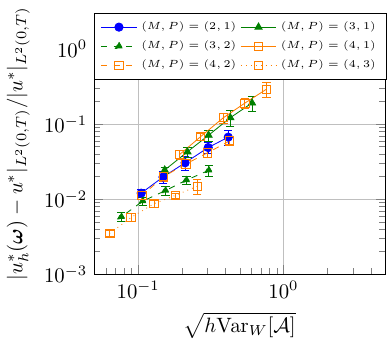}
\caption{Error in $u_h^*(t)$}
\label{fig:heat3D_errU}
\end{subfigure}\\
\begin{subfigure}{0.45\textwidth}
\includegraphics[scale=0.78]{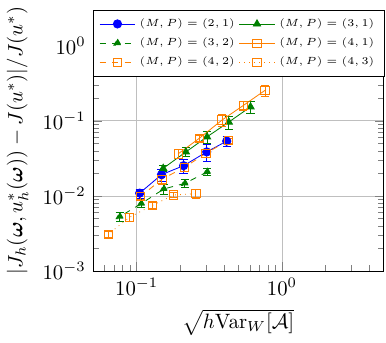}
\caption{Error in $J_h(u^*_h)$}
\label{fig:heat3D_errJ}
\end{subfigure} ~
\begin{subfigure}{0.45\textwidth}
\includegraphics[scale=0.78]{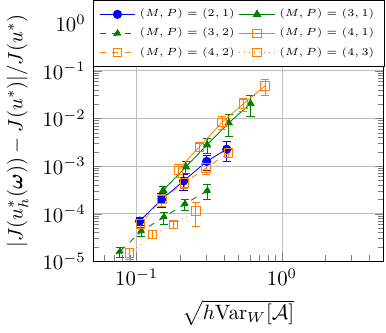}
\caption{Error in $J(u^*_h)$}
\label{fig:heat3D_errJ2}
\end{subfigure}\\
\begin{subfigure}{0.45\textwidth}
\includegraphics[scale=0.78]{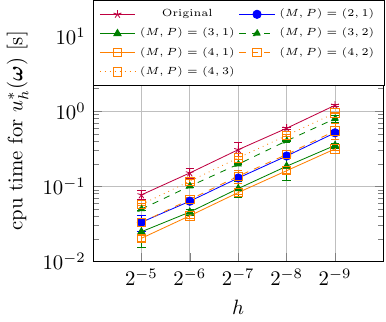}
\caption{Computational time for $x_h(t)$}
\label{fig:heat3D_durX}
\end{subfigure} ~
\begin{subfigure}{0.45\textwidth}
 \includegraphics[scale=0.78]{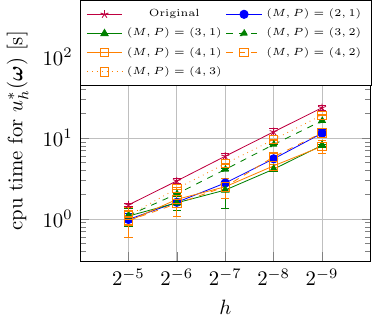}
\caption{Computational time for $u_h^*(t)$}
\label{fig:heat3D_durU}
\end{subfigure}
\caption{Results for the discretized 3D heat equation}
\label{fig:heat3D}
\end{figure}

We again consider uniform time grids with a grid spacing $h$. In each time interval $[t_{k-1}, t_k)$, we randomly use $P$ of the $M$ submatrices simultaneously. In our formalism, we thus assign a probability $1/{M \choose P}$ to each of the ${M \choose P}$ subsets of $\{1,2, \ldots ,M \}$ of size $P$. The states $x_h(\boldsymbol \omega,t)$ and the optimal controls $u^*_h(\boldsymbol \omega,t)$ are computed in the same way as for the example in the previous subsection. 

The obtained results are presented in Figure \ref{fig:heat3D}. The average errors (indicated by the markers) and the 2$\sigma$-confidence intervals (indicated by the error bars) are now estimated based on 10 realizations of $\boldsymbol \omega$. Figures \ref{fig:heat3D_errX}--\ref{fig:heat3D_errJ2} again show the convergence rates expected based on our theoretical results, just as in Figures \ref{fig:heat1D_errX}--\ref{fig:heat1D_errJ2} for the example in the previous subsection. We also observe that the errors are smaller when larger parts of $A$ are used simultaneously, i.e., when $P/M$ is larger. 

Figures \ref{fig:heat3D_durX} and \ref{fig:heat3D_durU} also show a computational advantage of the proposed method. Naturally, the computational advantage increases when the matrix $\mathcal{A}_h(\boldsymbol \omega,t)$ is more sparse, i.e., when $P/M$ is smaller. This situation is significantly different from the 1D heat equation considered in the previous subsection. For that example, the proposed method did not lead to any computational advantage. Apart from the larger state dimension $N$ in the 3D example, this difference seems to be related to the more `dense interconnection structure' of the 3D problem (in which every node is typically connected to 6 neighboring nodes) compared to the 1D problem (in which every node is connected to two neighboring nodes). This idea will be explored further in the next subsection in which we consider a model with an even denser interconnection structure. 

\subsection{A FE discretization of the fractional Laplacian} \label{ssec:example3}
We consider a controlled fractional heat equation on the a 1-D spatial domain $\xi \in [-L, L]$, 
\begin{align}
&y_t(t,\xi) = -(-d_\xi^2)^sy(t,\xi) + \chi_{[-L/3,0]}(\xi) u_1(t) + \chi_{[L/3,2L/3]}(\xi) u_2(t), \label{eq:fl1D_PDE} \\
&y(t,-L) = y(t,L) = 0, \qquad \qquad y(0,\xi) = e^{-\beta^2\xi^2} - e^{-\beta^2L^2},  \label{eq:fl1D_IC}
\end{align}
with the fractional power $s \in (0,1)$. We fix $s = 0.7$, $L = 5$, and $\beta = 0.4$. Note that the control $u(t) = [u_1(t), u_2(t)]^\top$ now has two components.  Our aim is to compute the optimal control $u^*(t) = [u^*_1(t), u^*_2(t)]^\top$ that minimizes
\begin{equation}
\mathcal{J}(u) = \frac{100}{2}\int_0^T \int_{-L}^L y(t,\xi)^2 \ \mathrm{d}\xi \ \mathrm{d}t + \frac{1}{2}\int_0^T \left( u_1(t)^2 + u_2^2(t) \right) \ \mathrm{d}t. \label{eq:fl1D_cost}
\end{equation}
A Finite Element (FE) discretization of \eqref{eq:fl1D_PDE}--\eqref{eq:fl1D_IC} with $N+1$ linear elements of equal length takes the form
\begin{equation}
E \dot{x}(t) = Ax(t) + Bu(t), \qquad \qquad x(0) = x_0, \label{eq:fl1D_disc}
\end{equation}
where the state $x(t)$ evolves in $\mathbb{R}^N$. 
Note that \eqref{eq:fl1D_disc} now also contains the symmetric and positive definite mass matrix $E$ and is thus not exactly of the form \eqref{eq:dyn_x}, but that the proposed method also applies to systems of this form. An explicit expression for the stiffness matrix $A$ can be found in \cite{biccari2019}. Because the fractional Laplacian is a nonlocal operator, all elements of $A$ are nonzero. From the expressions for the coefficients of $A$ in \cite{biccari2019} we can verify that $A$ is symmetric and diagonally dominant, i.e.\
\begin{equation}
-[A]_{ii} > \sum_{\substack{j=1 \\ j \neq i}}^N \lvert [A]_{ij} \rvert. 
\end{equation}
We can now write 
\begin{equation}
A = \sum_{\substack{j=1 \\ j \geq i}}^N  \tilde{A}_{ij} = \sum_{\substack{j=1 \\ j > i}}^N  \tilde{A}_{ij} + \sum_{i=1}^N \tilde{A}_{ii}, 
\end{equation}
where the matrices $A_{ij} \in \mathbb{R}^{N \times N}$ ($j \geq i$) are zero except for the coefficients
\begin{equation}
\begin{bmatrix}
[\tilde{A}_{ij}]_{ii} & [\tilde{A}_{ij}]_{ij} \\
[\tilde{A}_{ij}]_{ji} & [\tilde{A}_{ij}]_{jj}
\end{bmatrix} = \begin{bmatrix}
-\lvert [A]_{ij} \rvert & [A]_{ij} \\ 
[A]_{ij} & -\lvert [A]_{ij} \rvert
\end{bmatrix}, \quad
[A_{ii}]_{ii} = [A]_{ii} + \sum_{\substack{j=1 \\ j \neq i}}^N \lvert  [A]_{ij} \rvert.
\end{equation}
Again, it is easy to verify that the matrices $A_{ij}$ ($j \geq i$) are dissipative. 

\begin{figure}[h!]
\begin{subfigure}{0.45\textwidth}
\includegraphics[scale=0.78]{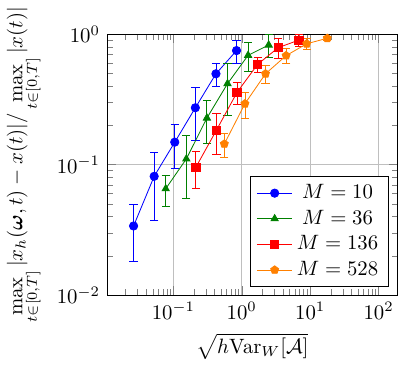}
\caption{Error in $x_h(t)$}
\label{fig:fl1D_errX}
\end{subfigure} ~
\begin{subfigure}{0.45\textwidth}
\includegraphics[scale=0.78]{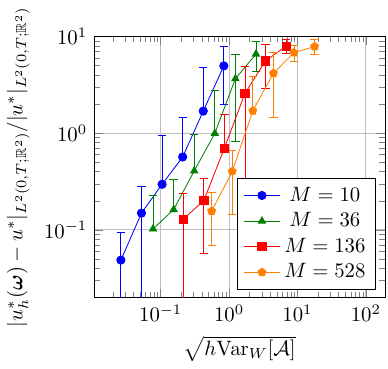}
\caption{Error in $u_h^*(t)$}
\label{fig:fl1D_errU}
\end{subfigure}\\
\begin{subfigure}{0.45\textwidth}
\includegraphics[scale=0.78]{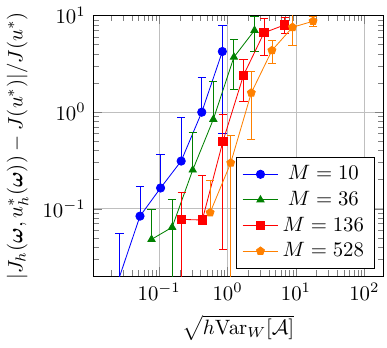}
\caption{Error in $J_h(u^*_h)$}
\label{fig:fl1D_errJ}
\end{subfigure} ~
\begin{subfigure}{0.45\textwidth}
\includegraphics[scale=0.78]{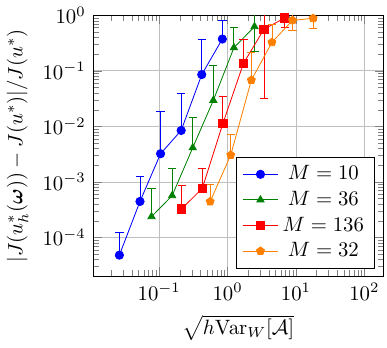}
\caption{Error in $J(u^*_h)$}
\label{fig:fl1D_errJ2}
\end{subfigure}\\
\begin{subfigure}{0.45\textwidth}
\includegraphics[scale=0.75]{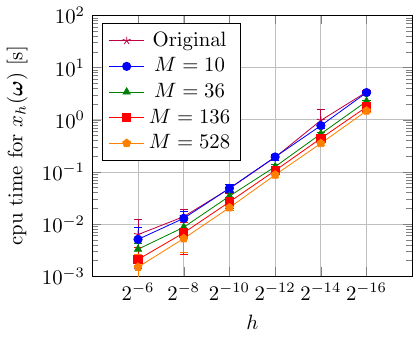}
\caption{Computational time for $x_h(t)$}
\label{fig:fl1D_durX}
\end{subfigure} ~
\begin{subfigure}{0.45\textwidth}
\includegraphics[scale=0.75]{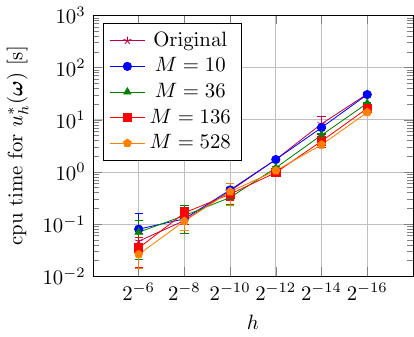}
\caption{Computational time for $u_h^*(t)$}
\label{fig:fl1D_durU}
\end{subfigure}
\caption{Results for the discretized 1D fractional heat equation with $s=0.7$}
\label{fig:fl1D}
\end{figure}

Now assume that $N$ is divisable by some number $P$. We then decompose $A$ into $M = P(P+1)/2$ submatrices $A_m$ as in \eqref{eq:Aform} by setting
\begin{equation}
A_{m(p,q)} = \sum_{i=i_{p-1}+1}^{i_p} \sum_{j=i_{q-1}+1}^{i_q} \tilde{A}_{ij}, \qquad \qquad q \geq p \in \{ 1,2, \ldots , P\},
\end{equation}
where $i_p = pN/P$ and $m(p,q)$ is a bijection 
\begin{equation}
m: \{ (p,q) \in \{1,2, \ldots , P \}^2 \mid q \geq p \} \rightarrow \{1,2, \ldots, P(P+1)/2 \}.
\end{equation}
We thus effectively decompose $A$ into $N/P \times N/P$ blocks, but we treat the diagonal in such a way that all submatrices $A_m$ are dissipative. 

We only use one of the matrices $A_m$ in each time interval $[t_{k-1},t_k)$ and thus assign uniform probabilities $2/(P(P+1))$ to each of the $M = P(P+1)/2$ subsets of $\{1,2, \ldots M \}$ of size 1. 

The results obtained for $N = 96$ are shown in Figure \ref{fig:fl1D}. The markers and the error bars in this figure again indicate the average and $2\sigma$-confidence interval estimated from 10 realizations of $\boldsymbol \omega$. Results are presented for for $P=4$, $8$, $16$, and $32$, which correspond to values of $M=10$, $36$, $136$, and $528$, respectively. Note that the number of submatrices $M$ is now much larger than in the previous two examples, and that also $h\mathrm{Var}[\mathcal{A}]$ and the relative errors are larger than in the previous examples. Figures \ref{fig:fl1D_errU} and \ref{fig:fl1D_errJ} even show relative errors that exceed 100\%. However, we still observe the convergence rates predicted by the theoretical results in Section \ref{sec:convergence} in Figures \ref{fig:fl1D_errX}--\ref{fig:fl1D_errJ2}. In particular, the convergence rate in Figure \ref{fig:fl1D_errJ2} is again twice as high as in the other figures. 

When we inspect the computational times in Figures \ref{fig:fl1D_durX} and \ref{fig:fl1D_durU}, we see that increasing $M$ decreases the computational time. In particular, solutions for $M = 528$ are typically computed 2-3 times faster than the solutions for the original dynamics. We expect that the computational advantage of the proposed method increases further when we increase the state dimension $N$.

\section{Conclusions and discussions} \label{sec:conc}
\subsection{Conclusions}

We have proposed a general framework for randomized time-splitting in LQ optimal control problems. It has been shown  that the dynamics, the minimal values of the cost functional, and the optimal control obtained with the proposed randomized time-splitting method converge in expectation to their analogues in the original problem when the grid spacing of the time grid goes to zero. The convergence rates in our theoretical results are also observed in three numerical examples. 

In two of the three considered examples, the proposed method leads to a typical reduction in computational cost of a factor 2-3. Only in the first example of a heat equation on a 1-D spatial domain, no reduction in computational cost could be observed. This seems to be the case because the matrix $A$ is not very large and already very sparse in this example.

\subsection{Extension to unbounded operators} \label{ssec:infdim}
We have considered finite-dimensional systems in this paper, but the numerical examples in Section \ref{sec:examples} are all obtained after spatial discretization of an infinite-dimensional system. A natural question is therefore whether our results can be extended to an infinite-dimensional setting. We already touched on this question in Remarks \ref{rem:def_varAW} and \ref{rem:infdim} and in Appendix \ref{app:commutative}. In particular, at the end of Appendix \ref{app:commutative} we indicate how results can be extended to an infinite dimensional setting under the (strong) additional assumptions that all operators $A_m$ commute and have the same domain $D(A_m)$.

It should be noted that the assumption that $D(A_m) = D(A)$ is very strong and will not be satisfied in many applications. A prototypical example is the splitting of an advection diffusion problem with zero Dirichlet boundary conditions (represented by $A$) in an advective part (represented by $A_1$) and a diffusive part (represented by $A_2$). Functions in $D(A_2)$ can then satisfy the zero Dirichlet boundary conditions on the whole boundary, but the functions in $D(A_1)$ only satisfy the zero Dirichlet boundary conditions on the parts of the boundary where the velocity field is pointing inward. The analysis of the RBM becomes much more subtle in these kind of situations. The numerical results in Figure \ref{fig:heat1Dconv} also seem to indicate that the proposed randomized time splitting method converges under weaker assumptions than the ones in Appendix \ref{app:commutative}. 

The technical difficulties encountered when weakening these assumptions are related to the difficulties in deterministic operator splitting with unbounded operators. These date back to the paper \cite{trotter1959} by Trotter, and have been an active field of research since then, see, e.g., \cite{kato1978, lapidus1981, neidhardt1998, hansen2008, osterman2013}. As the large literature on this topic indicates, determining the necessary conditions for the convergence of the proposed stochastic operator splitting method with unbounded operators is an interesting but challenging topic for future research. 

\subsection{Extension to nonlinear dynamics}
Another important topic for future research is the extension of our results for the linear quadratic optimal control problem to problems with nonquadratic cost functions constrained by nonlinear dynamics. This extension is particularly interesting because of the connections between the training of certain types of Deep Neural Networks (DNNs) and optimal control, see, e.g., \cite{E2017, benning2019, esteve2021, esteve2021sparse, ruizbalet2021}, and is also important for the control of interacting particles systems, see \cite{ko2020}. 

In the most general setting, we would replace the linear dynamics \eqref{eq:dyn_x} by the nonlinear dynamics
\begin{equation}
\dot{x}(t) = f(x(t), u(t)), \qquad \qquad x(0) = x_0, 
\end{equation}
where $f : \mathbb{R}^N \times \mathbb{R}^q \rightarrow \mathbb{R}^N$ is Lipschitz in the first variable $x$. 
As an analogue of \eqref{eq:Aform}, we then write (for $x \in \mathbb{R}^N$ and $u \in \mathbb{R}^q$)
\begin{equation}
f(x,u) = \sum_{m=1}^M f_m(x,u),
\end{equation}
for certain Lipschitz continuous functions $f_m : \mathbb{R}^N \times \mathbb{R}^q \rightarrow \mathbb{R}^N$. Similarly as in this paper, we choose a time grid $0 = t_0 < t_1 < t_2 < \ldots < t_K = T$, enumerate the subsets $S_1, S_2, \ldots ,S_{2^M}$ of $\{1,2, \ldots , M \}$ and assign probabilities $p_1, p_2, \ldots , p_{2^M}$ to them, and randomly select a $K$-tuple $\boldsymbol \omega = (\omega_1, \omega_2, \ldots , \omega_K)$ of indices $\omega_k \in \{1,2,\ldots 2^M \}$ according to the selected probabilities. We then  consider the (typically simpler) dynamics 
\begin{equation}
\dot{x}_h(\boldsymbol \omega, t) = \sum_{m \in S_{\omega_k}} \frac{f_m(x_h(\boldsymbol \omega,t), u_h(\boldsymbol \omega, t))}{\pi_m}, \qquad \qquad t \in [t_{k-1},t_k).
\end{equation}

Extending Theorem \ref{thm:este} (which considers the forward dynamics with a deterministic control $u_h(\boldsymbol \omega,t) = u(t)$) to such a nonlinear setting seems possible along the lines of the results for interacting-particle systems in \cite{jin2020}. The main difficulty is in Theorem \ref{thm:este2} where we use the variation of constants formula to obtain an estimate for a stochastic control $u_h(\boldsymbol \omega,t)$ (which depends on the randomly selected indices in $\boldsymbol \omega$). The variation of constants formula can be extended to a nonlinear setting, see, e.g., \cite{brauer1966}, but this leads to several additional complications which we aim to address in a future work. 

When an analogue of Theorem \ref{thm:este2} for nonlinear dynamics can be obtained, a bound on $\mathbb{E}[\lvert J_h(u_h) - J(u_h) \rvert]$ as in Lemma \ref{lem:dJuh} should follow relatively easily from a Lipschitz condition on the integrand in the considered cost function. An analogue of the no-gap condition, i.e., a bound on $\mathbb{E}[\lvert J(u_h^*) - J(u^*) \rvert]$, can then be obtained using classical arguments from the calculus of variations and the bound on $\mathbb{E}[\lvert J_h(u_h) - J(u_h) \rvert]$, similarly as for the linear-quadratic case in Theorem \ref{thm:nogap}. 

With these results, the suboptimality gap $\mathbb{E}[\lvert J_h(u^*_h) - J(u^*) \rvert]$ be bounded using the analogues of Lemma \ref{lem:dJuh} and Theorem \ref{thm:nogap} as follows. We start by noting that the triangle inequality shows that
\begin{equation}
\lvert J(u^*_h(\boldsymbol \omega)) - J(u^*) \rvert \leq \lvert J(u^*_h(\boldsymbol \omega)) - J_h(\boldsymbol \omega, u^*_h(\boldsymbol \omega)) \rvert + \lvert J_h(\boldsymbol \omega, u^*_h(\boldsymbol \omega)) - J(u^*) \rvert.
\end{equation}
Taking the expectation in this inequality, we see that the first term on the RHS can be bounded using (the analogue of) Lemma \ref{lem:dJuh} and the second term on the RHS can be bounded using (the analogue of) Theorem \ref{thm:nogap}. We thus obtain a bound on $\mathbb{E}[\lvert J_h(u^*_h) - J(u^*) \rvert]$ that is of order $\sqrt{h}$. It is interesting to observe that this rate is slower than the rate of order $h$ found for the linear-quadratic case in Corollary \ref{corr:suboptimality}. This difference seems to occur because Corollary \ref{corr:suboptimality} relies on the strict convexity of the functional, which is lost in a setting in which the dynamics are nonlinear. 

\subsection{Combination with model predictive control}
As suggested in \cite{ko2020}, it is natural to combine the proposed randomized time-splitting method with an MPC strategy. The resulting algorithm is essentially a receding horizon strategy, see, e.g., \cite{reble2012, amzi2016, amzi2018}, but we now use the proposed stochastic time-splitting method to approximate the optimal controls that need to be computed in each step. An important element of such a receding horizon strategy is that the optimal control is computed based on the current state of the original dynamics \eqref{eq:dyn_x}. This creates a feedback mechanism that provides additional robustness against the errors introduced by the proposed stochastic time-splitting method. 

The receding horizon strategy introduces two additional parameters in the control algorithm: the prediction horizon $\hat{T}$ and the control horizon $\tau$.
When the prediction horizon $\hat{T}$ is too short, the difference between the controls computed on the prediction horizon $[0, \hat{T}]$ and the desired optimal control on $[0, \infty)$ will be large. Decreasing the control horizon $\tau$ strengthens the feedback mechanism of the MPC strategy, which will likely allow for larger errors in the proposed stochastic time-splitting method. This idea could be formalized further by deriving an explicit error estimate that demonstrates the interaction of the control horizon $\tau$ and $h\mathrm{Var}[\mathcal{A}]$ (which characterizes the accuracy of the proposed random time-splitting method).

\backmatter
\bmhead{Acknowledgments}

This project has received funding from the European Research Council (ERC) under the European Union’s Horizon 2020 research and innovation programme (grant agreement NO: 694126-DyCon),  the Alexander von Humboldt-Professorship program, the European Unions Horizon 2020 research and innovation programme under the Marie Sklodowska-Curie grant agreement No.765579-ConFlex and the Transregio 154 Project ‘‘Mathematical Modelling, Simulation and Optimization Using the Example of Gas Networks’’, project C08, of the German DFG,   the Grant MTM2017-92996-C2-1-R COSNET of MINECO (Spain) and  the Elkartek grant KK-2020/00091 CONVADP of the Basque government.

\begin{appendices}

\section{Interacting particle systems in the proposed framework} \label{app:interacting_particles}

In this appendix, we explain the connection of our framework to the previously proposed RBMs for interacting particle systems in \cite{jin2020, jin2020convergence,li2020,ko2020}. We consider a (linearized first-order) system of $N$ interacting particles
\begin{equation}
\dot{x}_i(t) = \frac{1}{N-1} \sum_{\substack{j=1\\j\neq i}}^N a_{ij}(x_j(t) - x_i(t)), \quad x_i(0) = x_{0,i}, \quad i \in \{1,2,\ldots N \}, \label{eq:ip_dyn_x_ip}
\end{equation}
where the $a_{ij} \in \mathbb{R}$ ($j \neq i$) are constants. To simplify the following exposition, we assume that the number of particles $N$ is divisible by some number $P > 1$.

We discuss here one particular RBM called RBM-1 in \cite{jin2020}, but other variants can be treated similarly. We first choose a time grid $0 = t_0 < t_1 < t_2 < \ldots < t_{K_1} < t_K = T$ in the time interval $[0,T]$.  In each time interval $[t_{k-1}, t_k)$, we then choose a random partition of the index set $\{1,2, \ldots, n \}$ into disjoint subsets $\mathcal{B}_r^k$ (also called batches) of size $P$ ($r \in \{1,2, \ldots , N/P \}$). We consider only the interactions between particles that are in the same batch. To formalize this idea, note that, in each time interval $[t_{k-1}, t_k)$, every particle $i$ is contained in precisely one batch $\mathcal{B}^k_{r(i,k)}$. We thus consider the dynamics
\begin{equation}
\dot{x}_{\mathrm{RBM},i}(t) = \frac{1}{P-1}\sum_{\substack{j\in \mathcal{B}^k_{r(i,k)} \\ j \neq i}} a_{ij}(x_{\mathrm{RBM},j}(t) - x_{\mathrm{RBM},i}(t)), \qquad x_i(0) = x_{0,i}. 
\label{eq:ip_dyn_x_rbm_ip}
\end{equation}
To connect this idea to our framework, we write \eqref{eq:ip_dyn_x_ip} in matrix form
\begin{equation}
\dot{x}(t) = A x(t), \qquad \qquad x(0) = x_0, \qquad \qquad 
 A = \frac{1}{N-1} \sum_{\substack{i,j = 1\\i \neq j}}^N \tilde{A}_{ij}, 
\end{equation}
where $x(t) = [x_1(t), x_2(t), \ldots x_N(t)]^\top$ and $x_0 = [x_{0,1}, x_{0,2}, \ldots , x_{0,N}]$ and the entries of the matrices $\tilde{A}_{ij}$ ($j \neq i$) are zero except for the entries
\begin{equation}
\begin{bmatrix}
[\tilde{A}_{ij}]_{ij} & [\tilde{A}_{ij}]_{ii}
\end{bmatrix} = a_{ij} \begin{bmatrix}
1 & -1
\end{bmatrix}. 
\end{equation}
Also the RBM-dynamics \eqref{eq:ip_dyn_x_rbm_ip} can be written in matrix form as
\begin{equation}
\dot{x}_{\mathrm{RBM}}(t) = \mathcal{A}_{\mathrm{RBM}}(t) x_{\mathrm{RBM}}(t), \qquad \qquad x_{\mathrm{RBM}}(0) = x_0,
\end{equation}
where
\begin{equation}
\mathcal{A}_{\mathrm{RBM}}(t) = \frac{1}{P-1} \sum_{r=1}^{N/P}\sum_{\{i,j\} \subseteq \mathcal{B}_r^k} \tilde{A}_{ij}, \qquad \qquad t \in [t_{k-1}, t_k). \label{eq:ip_ARBMt}
\end{equation}
Note that the probability that two distinct indices $i$ and $j$ are in the same batch (i.e., the probability that $j \neq i$ is in the batch $\mathcal{B}^k_{r(i,k)}$) is $(P-1)/(N-1)$ because there are $P-1$ of the $N-1$ places in $\mathcal{B}^k_{r(i,k)}$ remaining after the index $i$ has been fixed. This factor is also visible in the definitions of $A$ and $\mathcal{A}_{\mathrm{RBM}}(t)$. 

To make the connection to our proposed framework, we enumerate the $M = N(N-1)$ interaction matrices $A_{ij}$, i.e., we choose a bijection 
\begin{equation}
\mathfrak{m} : \{ (i,j) \in \{1,2, \ldots, N \}^2 \mid i \neq j \} \rightarrow \{1,2, \ldots, N(N-1) \},
\end{equation}
and set 
\begin{equation}
A_{\mathfrak{m}(i,j)} := \frac{1}{N-1} \tilde{A}_{ij}.
\end{equation}
We then need to assign probabilities $p_\omega$ to the $2^M$ subsets $S_\omega$ of $\{1,2, \ldots , M \}$. Naturally, we only assign nonzero probabilities to subsets $S_\omega$ that correspond to a partition $\dot{\cup}_r \mathcal{B}_r= \{1,2, \ldots, N \}$, i.e.\ sets of the form
\begin{equation}
S_\omega = \{ \mathfrak{m}(i,j) \mid \exists_{i,j,r}\ \mathrm{such\ that}\ i \neq j\ \mathrm{and}\ \{ i,j\} \subseteq \mathcal{B}_r \}. \label{eq:ip_Sl_form}
\end{equation}
Standard combinatorics shows that there are
\begin{equation}
\mathcal{N} = \frac{N!}{(P!)^{N/P} \left( N/P \right)!},
\end{equation}
distinct partitions of $N$ indices into $N/P$ subsets of size $P$. We assign a probability $p_\omega = 1 /\mathcal{N}$ to each of the subsets of the form \eqref{eq:ip_Sl_form}. 

It remains to compute the probabilities $\pi_m = \pi_{\mathfrak{m}(i,j)}$ defined in \eqref{eq:pim}, i.e.\ to determine how many of the subsets $S_\omega$ of the form \eqref{eq:ip_Sl_form} contain $m = \mathfrak{m}(i,j)$. 
When a certain batch $\mathcal{B}_{r^*}$ contains $i$ and $j$ ($j \neq i$) there are $N-2 \choose P-2$ ways to fill the remaining positions in $\mathcal{B}_{r^*}$ with $P-2$ of the $N-2$ remaining indices. Once the indices in $\mathcal{B}_{r^*}$ are fixed, there are
\begin{equation}
\mathcal{M} = \frac{(N-P)!}{(P!)^{N/P-1} \left( N/P - 1\right)!},
\end{equation}
ways to distribute the remaining $N-P$ indices into $N/P-1$ subsets of size $P$. We thus conclude that
\begin{equation}
\pi_m = \frac{{N-2 \choose P-2} \mathcal{M}}{\mathcal{N}}
\end{equation}
Using the formulas for $\mathcal{N}$ and $\mathcal{M}$, it can be verified that
\begin{equation}
\pi_m = \frac{P-1}{N-1}.
\end{equation}
It is now easy to verify that the definition of $\mathcal{A}_h(\boldsymbol \omega,t)$ in \eqref{eq:defAt} is equivalent to the definition of $\mathcal{A}_{\mathrm{RBM}}(t)$ in \eqref{eq:ip_ARBMt}. 

\section{An alternative for Corollary \ref{corr:este1}}
\label{app:commutative}

In this appendix, we will prove a result similar to Corollary \ref{corr:este1} under the additional assumption that all matrices commute. The proof is quite intuitive and gives an idea about how the results in this paper can be generalized to an infinite dimensional setting. 

The analysis in this appendix uses the following additional assumption. 
\begin{assumption} \label{ass:commutative}
Suppose that the matrices $A_1, A_2, \ldots , A_M$ all commute pairwise, i.e.\
\begin{equation}
A_{m}A_{m'} = A_{m'}A_m,
\end{equation}
for all $m, m' \in \{1,2, \ldots , M \}$. 
\end{assumption}

Also observe that for any two dissipative matrices $X,Y \in \mathbb{R}^{N \times N}$ and vector $x_0 \in \mathbb{R}^N$ we have that
\begin{align}
\lvert e^Xx_0 - e^Yx_0 \rvert
&= \left\lvert \int_0^1 \frac{d}{d\tau} e^{X \tau + Y(1-\tau)}x_0 \ \mathrm{d}\tau \right\rvert \nonumber \\
&\leq \int_0^1 \| e^{X \tau + Y(1-\tau)} \| \lvert (X-Y)x_0 \rvert  \mathrm{d}\tau
\leq \lvert (X - Y)x_0 \rvert, \label{eq:diff_mat_exp}
\end{align}
where it was used that $X\tau + Y(1-\tau)$ is dissipative for $\tau \in [0,1]$ because $X$ and $Y$ are dissipative by assumption. 

\begin{theorem} \label{thm:commutative}
Under Assumptions \ref{ass:dissipativity}, \ref{ass:positivity}, and \ref{ass:commutative}, we have that
\begin{equation}
\mathbb{E}[\| S_h(t,s)x_0 - e^{A(t-s)}x_0 \|^2] \leq 2 h(t-s)  \mathrm{Var}_W[\mathcal{A}] \lvert W^{-1}x_0 \rvert^2, \label{eq:thm_commutative}
\end{equation}
for all $0 \leq s \leq t \leq T$, all $x_0 \in \mathbb{R}^N$, and all invertible matrices $W$. 
\end{theorem}

\begin{proof} We use the notation from Remark \ref{rem:Shnew}, so $\ell$ and $k$ are such that $s \in [t_{\ell-1}, t_\ell)$ and $t \in [t_{k-1}, t_k)$, $\tilde{K} = k - \ell + 1$, and 
\begin{equation}
\tilde{t}_0 := s < 
\tilde{t}_1 := t_\ell < 
\tilde{t}_2 := t_{\ell+1} < 
\ldots < 
\tilde{t}_{\tilde{K}-1} := t_{k-1}  < 
\tilde{t}_{\tilde{K}} := t,  \label{eq:thm_comm_step1}
\end{equation}
see also Figure \ref{fig:timegrids} on page \pageref{fig:timegrids}. 
Furthermore, we denote $\tilde{h}_p := \tilde{t}_p - \tilde{t}_{p-1}$ for $p \in \{1,2, \ldots, \tilde{K} \}$ and denote $\mathcal{A}_\omega := \sum_{m \in S_\omega} A_m/\pi_m$ for $\omega \in \{1,2, \ldots ,2^M \}$. Note that $\mathcal{A}_h(\boldsymbol \omega,\tau) = \mathcal{A}_{\omega_p}$ for $\tau \in [\tilde{t}_{p-1}, \tilde{t}_p)$ and that $\mathcal{A}_\omega$ is dissipative for all $\omega \in \{1,2, \ldots , 2^M \}$ because of Assumption \ref{ass:dissipativity}. 

Because the matrices $\mathcal{A}_\omega$ (with $\omega \in \{1,2, \ldots, 2^M \}$) all commute pairwise due to Assumption \ref{ass:commutative}, the formula for $S_h(\boldsymbol \omega, t,s)$ in \eqref{eq:Udef} in Remark \ref{rem:Shnew} reduces to
\begin{equation}
S_h(\boldsymbol \omega, t,s)x_0 = \exp\left(\sum_{p=1}^{\tilde{K}} \mathcal{A}_{\omega_{p+\ell-1}}\tilde{h}_p \right)x_0. \label{eq:thm_comm_step2}
\end{equation}
Because Assumption \ref{ass:dissipativity} implies that the matrix in the exponent in the formula above and $A$ are both dissipative, \eqref{eq:diff_mat_exp} can be applied to find that
\begin{equation}
\lvert S_h(\boldsymbol \omega, t, s)x_0 - e^{A(t-s)}x_0 \rvert \leq \left\lvert \sum_{p=1}^{\tilde{K}} \left( \mathcal{A}_{\omega_{p+\ell-1}} - A \right)\tilde{h}_p x_0 \right\rvert, \label{eq:thm_comm_step3}
\end{equation}
where it was used that $\sum_{p=1}^{\tilde{K}} \tilde{h}_p = t-s$. Squaring this expression yields
\begin{multline}
\lvert S_h(\boldsymbol \omega, t, s)x_0 - e^{A(t-s)}x_0 \rvert^2 \leq \\ 
\sum_{p,p'=1}^{\tilde{K}} \tilde{h}_p \tilde{h}_{p'} \langle (\mathcal{A}_{\omega_{p+\ell-1}} - A) x_0, (\mathcal{A}_{\omega_{p'+\ell-1}} - A ) x_0 \rangle. \label{eq:thm_comm_step4}
\end{multline}
When we take the expected value, the terms with $p \neq p'$ disappear because 
\begin{align}
\mathbb{E}[\langle & (\mathcal{A}_{\omega_{p+\ell-1}} - A) x_0, (\mathcal{A}_{\omega_{p'+\ell-1}} - A ) x_0 \rangle] \nonumber \\
&= \sum_{\omega = 1}^{2^M}  \sum_{\omega' = 1}^{2^M} \langle (\mathcal{A}_{\omega} - A) x_0, (\mathcal{A}_{\omega'} - A ) x_0 \rangle p_{\omega} p_{\omega'} \nonumber \\
&= \bigg\langle  \sum_{\omega = 1}^{2^M} (\mathcal{A}_\omega - A)x_0, \sum_{\omega'=1}^{2^M} (\mathcal{A}_{\omega'} - A)x_0 \bigg\rangle = \langle 0, 0 \rangle = 0 \label{eq:thm_comm_step5}
\end{align}
where the first identity follows after writing $\omega = \omega_{p-\ell+1}$ and $\omega'= \omega_{p'-\ell+1}$, and the second to last identity from \eqref{eq:expectationA_step1} and \eqref{eq:sump}. Therefore, only the terms with $p = p'$ remain after taking the expected value of \eqref{eq:thm_comm_step4} and
\begin{align}
&\mathbb{E}[\lvert S_h(t, s)x_0 - e^{A(t-s)}x_0 \rvert^2] \nonumber \\
& \leq \sum_{\omega_\ell = 1}^{2^M} \sum_{\omega_{\ell+1} = 1}^{2^M}  \cdots \sum_{\omega_{\ell+\tilde{K}-1} = 1}^{2^M} \sum_{p=1}^{\tilde{K}} \tilde{h}_p^2 \lvert (\mathcal{A}_{\omega_{p+\ell-1}} - A) x_0 \rvert^2 p_{\omega_\ell} p_{\omega_{\ell+1}} \ldots p_{\omega_{\ell + \tilde{K}-1}} \nonumber \\
&= \sum_{p=1}^{\tilde{K}} \tilde{h}_p^2 \sum_{\omega=1}^{2^M} \lvert (\mathcal{A}_{\omega_p+\ell-1} - A)x_0 \rvert^2 p_\omega. \label{eq:thm_comm_step6}
\end{align}
The proof is completed with two straightforward observations. 
First of all, note that because $\tilde{h}_p \leq h$
\begin{equation}
\sum_{p=1}^{\tilde{K}} \tilde{h}_p^2 \leq \sum_{p=1}^{\tilde{K}} h \tilde{h}_p = h \sum_{p=1}^{\tilde{K}} \tilde{h}_p = h(t-s).  \label{eq:thm_comm_step7}
\end{equation}
Secondly, we have that
\begin{align}
\sum_{\omega=1}^{2^M} \lvert (\mathcal{A}_{\omega_p+\ell-1} - A)x_0 \rvert^2 p_\omega &= \sum_{\omega=1}^{2^M} \lvert (\mathcal{A}_{\omega_p+\ell-1} - A)WW^{-1}x_0 \rvert^2 p_\omega \nonumber \\
&\leq \sum_{\omega=1}^{2^M} \|(\mathcal{A}_{\omega_p+\ell-1} - A)W \|^2  \lvert W^{-1}x_0 \rvert^2 p_\omega.  \label{eq:thm_comm_step8}
\end{align}
The result follows after inserting \eqref{eq:thm_comm_step7} and \eqref{eq:thm_comm_step8} into \eqref{eq:thm_comm_step6}. 
\end{proof}

The proof of Theorem \ref{thm:commutative} extends naturally to an infinite dimensional setting as follows. Most of the definitions and notations from Section \ref{sec:method} remain unchanged, apart from the following. 
\begin{itemize}
\item The state and the control no longer evolve in the finite-dimensional spaces $\mathbb{R}^N$ and $\mathbb{R}^q$, but in the (potentially) infinite-dimensional Hilbert spaces $X$ and $U$, respectively. 

\item $A$ and $A_m$ (with $m \in \{1,2, \ldots, M \}$) now represent the generators of $C_0$-semigroups $e^{At}$ and $e^{A_mt}$ on the Hilbert space $X$ with domains $D(A)$ and $D(A_m)$, respectively. 

\item $B$ is now a bounded linear operator from $U$ to $X$. 
\end{itemize}

For simplicity we assume that the domains of the operators $A_m$ are all the same and equal to the domain of $A$, i.e.\ $D(A_m) = D(A)$. For a value of $\lambda$ in the resolvent set of $A$, the resolvent $W = (A - \lambda I)^{-1}$ is a bounded operator $X \rightarrow D(A) \subset X$ with (unbounded) inverse $A - \lambda I$ and one now easily verifies that $AW$ and $A_mW$ represent bounded operators on $X$, meaning that $\mathrm{Var}_W[\mathcal{A}]$ as introduced in Remark \ref{rem:def_varAW} is bounded. For $\lvert W^{-1}x_0 \rvert = \lvert (A - \lambda I)x_0 \rvert$ to be bounded, we require that $x_0 \in D(A)$. The proof of Theorem \ref{thm:commutative} can thus be applied in this setting with the additional assumption that $x_0 \in D(A)$. The proof remains effectively unchanged. 

Note that when we want to use Theorem \ref{thm:commutative} to obtain a result similar to Theorem \ref{thm:este2}, we also need a smoothness assumption on the input operator $B$. 
In particular, similarly as \eqref{eq:thm_este2_step0} in Theorem \ref{thm:este2}, we would then like to bound
\begin{equation}
 \int_0^t \left\lvert (S_h(\boldsymbol \omega,t,s) - e^{A(t-s)})Bu_h(\boldsymbol \omega,s) \right\rvert \ \mathrm{d}s, 
\end{equation}
which is only possible with Theorem \ref{thm:commutative} when $\lvert W^{-1}Bu_h(\boldsymbol \omega,s) \rvert$ is finite. To this end one would typically require that the range of $B$ is contained in $D(A)$. 

\end{appendices}

\bibliographystyle{abbrvnat}
\bibliography{references}% common bib file


\end{document}